\newcount\mgnf\newcount\tipi\newcount\tipoformule\newcount\greco

\tipi=2          
\tipoformule=0   


\global\newcount\numsec
\global\newcount\numfor
\global\newcount\numtheo
\global\advance\numtheo by 1

\def\senondefinito#1{\expandafter\ifx\csname#1\endcsname\relax}

\def\SIA #1,#2,#3 {\senondefinito{#1#2}%
\expandafter\xdef\csname #1#2\endcsname{#3}\else
\write16{???? ma #1,#2 e' gia' stato definito !!!!} \fi}

\def\etichetta(#1){(\veroparagrafo.\veraformula)%
\SIA e,#1,(\veroparagrafo.\veraformula) %
\global\advance\numfor by 1%
\write15{\string\FU (#1){\equ(#1)}}%
\write16{ EQ #1 ==> \equ(#1) }}

\def\letichetta(#1){\veroparagrafo.\verotheo
\SIA e,#1,{\veroparagrafo.\verotheo}
\global\advance\numtheo by 1
\write15{\string\FU (#1){\equ(#1)}}
\write16{ Sta \equ(#1) == #1 }}

\def\tetichetta(#1){\veroparagrafo.\veraformula 
\SIA e,#1,{(\veroparagrafo.\veraformula)}
\global\advance\numfor by 1
\write15{\string\FU (#1){\equ(#1)}}
\write16{ tag #1 ==> \equ(#1)}}

\def\FU(#1)#2{\SIA fu,#1,#2 }

\def\etichettaa(#1){(A\veroparagrafo.\veraformula)%
\SIA e,#1,(A\veroparagrafo.\veraformula) %
\global\advance\numfor by 1%
\write15{\string\FU (#1){\equ(#1)}}%
\write16{ EQ #1 ==> \equ(#1) }}

\def\BOZZA{
\def\alato(##1){%
 {\rlap{\kern-\hsize\kern-1.4truecm{$\scriptstyle##1$}}}}%
\def\aolado(##1){%
 {
{
 \rlap{\kern-1.4truecm{$\scriptstyle##1$}}}}}
 }

\def\alato(#1){}
\def\aolado(#1){}

\def\veroparagrafo{\number\numsec}
\def\veraformula{\number\numfor}
\def\verotheo{\number\numtheo}

\def\Eq(#1){\eqno{\etichetta(#1)\alato(#1)}}
\def\eq(#1){\etichetta(#1)\alato(#1)}
\def\leq(#1){\leqno{\aolado(#1)\etichetta(#1)}}
\def\teq(#1){\tag{\aolado(#1)\tetichetta(#1)\alato(#1)}}
\def\Eqa(#1){\eqno{\etichettaa(#1)\alato(#1)}}
\def\eqa(#1){\etichettaa(#1)\alato(#1)}
\def\eqv(#1){\senondefinito{fu#1}$\clubsuit$#1
\write16{#1 non e' (ancora) definito}%
\else\csname fu#1\endcsname\fi}
\def\equ(#1){\senondefinito{e#1}\eqv(#1)\else\csname e#1\endcsname\fi}

\def\Lemma(#1){\aolado(#1)Lemma \letichetta(#1)}%
\def\Theorem(#1){{\aolado(#1)Theorem \letichetta(#1)}}%
\def\Proposition(#1){\aolado(#1){Proposition \letichetta(#1)}}%
\def\Corollary(#1){{\aolado(#1)Corollary \letichetta(#1)}}%
\def\Remark(#1){{\noindent\aolado(#1){\bf Remark \letichetta(#1).}}}%
\def\Definition(#1){{\noindent\aolado(#1){\bf Definition 
\letichetta(#1)$\!\!$\hskip-1.6truemm}}}
\def\Example(#1){\aolado(#1) Example \letichetta(#1)$\!\!$\hskip-1.6truemm}

\def\include#1{
\openin13=#1.aux \ifeof13 \relax \else
\input #1.aux \closein13 \fi}

\openin14=\jobname.aux \ifeof14 \relax \else
\input \jobname.aux \closein14 \fi
\openout15=\jobname.aux






\newcount\pgn \pgn=1
\def\foglio{\number\numsec:\number\pgn
\global\advance\pgn by 1}
\def\foglioa{A\number\numsec:\number\pgn
\global\advance\pgn by 1}

\footline={\rlap{\hbox{\copy200}}\hss\tenrm\folio\hss}

\def\TIPIO{
\font\setterm=amr7 
\def \settepunti{\def\rm{\fam0\setterm}
\textfont0=\setterm   
\normalbaselineskip=9pt\normalbaselines\rm }\let\nota=\settepunti}

\def\TIPITOT{
\font\twelverm=cmr12
\font\twelvei=cmmi12
\font\twelvesy=cmsy10 scaled\magstep1
\font\twelveex=cmex10 scaled\magstep1
\font\twelveit=cmti12
\font\twelvett=cmtt12
\font\twelvebf=cmbx12
\font\twelvesl=cmsl12
\font\ninerm=cmr9
\font\ninesy=cmsy9
\font\eightrm=cmr8
\font\eighti=cmmi8
\font\eightsy=cmsy8
\font\eightbf=cmbx8
\font\eighttt=cmtt8
\font\eightsl=cmsl8
\font\eightit=cmti8
\font\sixrm=cmr6
\font\sixbf=cmbx6
\font\sixi=cmmi6
\font\sixsy=cmsy6
\font\twelvetruecmr=cmr10 scaled\magstep1
\font\twelvetruecmsy=cmsy10 scaled\magstep1
\font\tentruecmr=cmr10
\font\tentruecmsy=cmsy10
\font\eighttruecmr=cmr8
\font\eighttruecmsy=cmsy8
\font\seventruecmr=cmr7
\font\seventruecmsy=cmsy7
\font\sixtruecmr=cmr6
\font\sixtruecmsy=cmsy6
\font\fivetruecmr=cmr5
\font\fivetruecmsy=cmsy5
\textfont\truecmr=\tentruecmr
\scriptfont\truecmr=\seventruecmr
\scriptscriptfont\truecmr=\fivetruecmr
\textfont\truecmsy=\tentruecmsy
\scriptfont\truecmsy=\seventruecmsy
\scriptscriptfont\truecmr=\fivetruecmr
\scriptscriptfont\truecmsy=\fivetruecmsy
\def \eightpoint{\def\rm{\fam0\eightrm}
\textfont0=\eightrm \scriptfont0=\sixrm \scriptscriptfont0=\fiverm
\textfont1=\eighti \scriptfont1=\sixi   \scriptscriptfont1=\fivei
\textfont2=\eightsy \scriptfont2=\sixsy   \scriptscriptfont2=\fivesy
\textfont3=\tenex \scriptfont3=\tenex   \scriptscriptfont3=\tenex
\textfont\itfam=\eightit  \def\it{\fam\itfam\eightit}%
\textfont\slfam=\eightsl  \def\sl{\fam\slfam\eightsl}%
\textfont\ttfam=\eighttt  \def\tt{\fam\ttfam\eighttt}%
\textfont\bffam=\eightbf  \scriptfont\bffam=\sixbf
\scriptscriptfont\bffam=\fivebf  \def\bf{\fam\bffam\eightbf}%
\tt \ttglue=.5em plus.25em minus.15em
\setbox\strutbox=\hbox{\vrule height7pt depth2pt width0pt}%
\normalbaselineskip=9pt
\let\sc=\sixrm  \let\big=\eightbig  \normalbaselines\rm
\textfont\truecmr=\eighttruecmr
\scriptfont\truecmr=\sixtruecmr
\scriptscriptfont\truecmr=\fivetruecmr
\textfont\truecmsy=\eighttruecmsy
\scriptfont\truecmsy=\sixtruecmsy }\let\nota=\eightpoint}

\newfam\msbfam   
\newfam\truecmr  
\newfam\truecmsy 
\newskip\ttglue
\ifnum\tipi=0\TIPIO \else\ifnum\tipi=1 \TIPI\else \TIPITOT\fi\fi

\def\a{\alpha}
\def\b{\beta}
\def\d{\delta}
\def\e{\epsilon}

\def\vf{\varphi}
\def\g{\gamma}
\def\k{\kappa}
\def\l{\lambda}
\def\r{\rho}
\def\s{\sigma}
\def\t{\tau}
\def\th{\theta}

\def\z{\zeta}
\def\o{\omega}
\def\D{\Delta}
\def\L{\Lambda}
\def\G{\Gamma}
\def\O{\Omega}
\def\S{\Sigma}

\def\del #1{\frac{\partial^{#1}}{\partial\l^{#1}}}

\def\E{{I\kern-.25em{E}}}
\def\N{{I\kern-.25em{N}}}
\def\M{{I\kern-.25em{M}}}
\def\R{{I\kern-.25em{R}}}
\def\Z{{Z\kern-.425em{Z}}}
\def\1{{1\kern-.25em\hbox{\rm I}}}
\def\eu{{1\kern-.25em\hbox{\sm I}}}

\def\C{{I\kern-.64em{C}}}
\def\P{{I\kern-.25em{P}}}
\def\eop{{ \vrule height7pt width7pt depth0pt}\par\bigskip}

\def\del{\partial}


\def\AA{{\cal A}}
\def\BB{{\cal B}}
\def\CC{{\cal C}}
\def\DD{{\cal D}}
\def\EE{{\cal E}}
\def\FF{{\cal F}}
\def\GG{{\cal G}}

\def\II{{\cal I}}

\def\KK{{\cal K}}

\def\OO{{\cal O}}
\def\PP{{\cal P}}
\def\SS{{\cal S}}
\def\TT{{\cal T}}
\def\NN{{\cal N}}
\def\MM{{\cal M}}
\def\WW{{\cal W}}

\def\UU{{\cal U}}

\def\YY{{\cal Y}}

\def\RR{{\cal R}}

\def\chap #1#2{\line{\ch #1\hfill}\numsec=#2\numfor=1}
\def\sign{\,\hbox{sign}\,}

\def\sqr#1#2{{\vcenter{\vbox{\hrule height.#2pt
     \hbox{\vrule width.#2pt height#1pt \kern#1pt
   \vrule width.#2pt}\hrule height.#2pt}}}}


\newcount\foot
\foot=1
\def\note#1{\footnote{${}^{\number\foot}$}{\ftn #1}\advance\foot by 1}
\def\tag #1{\eqno{\hbox{\rm(#1)}}}
\def\frac#1#2{{#1\over #2}}
\def\sfrac#1#2{{\textstyle{#1\over #2}}}
\def\text#1{\quad{\hbox{#1}}\quad}

\def\proof{{\noindent\pr Proof: }}

\def\remark{\noindent{\bf Remark: }}
\def\thanks{\noindent{\bf Aknowledgements: }}
\font\pr=cmbxsl10


\font\ch=cmbx12
\font\ftn=cmr8

\font\it=cmti10
\font\bf=cmbx10
\font\sm=cmr7

%
\catcode`\X=12\catcode`\@=11
\def\n@wcount{\alloc@0\count\countdef\insc@unt}
\def\n@wwrite{\alloc@7\write\chardef\sixt@@n}
\def\n@wread{\alloc@6\read\chardef\sixt@@n}
\def\crossrefs#1{\ifx\alltgs#1\let\tr@ce=\alltgs\else\def\tr@ce{#1,}\fi
   \n@wwrite\cit@tionsout\openout\cit@tionsout=\jobname.cit 
   \write\cit@tionsout{\tr@ce}\expandafter\setfl@gs\tr@ce,}
\def\setfl@gs#1,{\def\@{#1}\ifx\@\empty\let\next=\relax
   \else\let\next=\setfl@gs\expandafter\xdef
   \csname#1tr@cetrue\endcsname{}\fi\next}
\newcount\sectno\sectno=0\newcount\subsectno\subsectno=0\def\r@s@t{\relax}
\def\resetall{\global\advance\sectno by 1\subsectno=0
  \gdef\firstpart{\number\sectno}\r@s@t}
\def\resetsub{\global\advance\subsectno by 1
   \gdef\firstpart{\number\sectno.\number\subsectno}\r@s@t}
\def\v@idline{\par}\def\firstpart{\number\sectno}
\def\l@c@l#1X{\firstpart.#1}\def\gl@b@l#1X{#1}\def\t@d@l#1X{{}}
\def\m@ketag#1#2{\expandafter\n@wcount\csname#2tagno\endcsname
     \csname#2tagno\endcsname=0\let\tail=\alltgs\xdef\alltgs{\tail#2,}%
  \ifx#1\l@c@l\let\tail=\r@s@t\xdef\r@s@t{\csname#2tagno\endcsname=0\tail}\fi
   \expandafter\gdef\csname#2cite\endcsname##1{\expandafter
     \ifx\csname#2tag##1\endcsname\relax?\else{\rm\csname#2tag##1\endcsname}\fi
    \expandafter\ifx\csname#2tr@cetrue\endcsname\relax\else
     \write\cit@tionsout{#2tag ##1 cited on page \folio.}\fi}%
   \expandafter\gdef\csname#2page\endcsname##1{\expandafter
     \ifx\csname#2page##1\endcsname\relax?\else\csname#2page##1\endcsname\fi
     \expandafter\ifx\csname#2tr@cetrue\endcsname\relax\else
     \write\cit@tionsout{#2tag ##1 cited on page \folio.}\fi}%
   \expandafter\gdef\csname#2tag\endcsname##1{\global\advance
     \csname#2tagno\endcsname by 1%
   \expandafter\ifx\csname#2check##1\endcsname\relax\else%
\fi
   \expandafter\xdef\csname#2check##1\endcsname{}%
   \expandafter\xdef\csname#2tag##1\endcsname
     {#1\number\csname#2tagno\endcsnameX}%
   \write\t@gsout{#2tag ##1 assigned number \csname#2tag##1\endcsname\space
      on page \number\count0.}%
   \csname#2tag##1\endcsname}}%
\def\m@kecs #1tag #2 assigned number #3 on page #4.%
   {\expandafter\gdef\csname#1tag#2\endcsname{#3}
   \expandafter\gdef\csname#1page#2\endcsname{#4}}
\def\re@der{\ifeof\t@gsin\let\next=\relax\else
    \read\t@gsin to\t@gline\ifx\t@gline\v@idline\else
    \expandafter\m@kecs \t@gline\fi\let \next=\re@der\fi\next}
\def\t@gs#1{\def\alltgs{}\m@ketag#1e\m@ketag#1s\m@ketag\t@d@l p
    \m@ketag\gl@b@l r \n@wread\t@gsin\openin\t@gsin=\jobname.tgs \re@der
    \closein\t@gsin\n@wwrite\t@gsout\openout\t@gsout=\jobname.tgs }
\outer\def\localtags{\t@gs\l@c@l}
\outer\def\globaltags{\t@gs\gl@b@l}
\outer\def\newlocaltag#1{\m@ketag\l@c@l{#1}}
\outer\def\newglobaltag#1{\m@ketag\gl@b@l{#1}}

\def\t@gsoff#1,{\def\@{#1}\ifx\@\empty\let\next=\relax\else\let\next=\t@gsoff
   \expandafter\gdef\csname#1cite\endcsname{\relax}
   \expandafter\gdef\csname#1page\endcsname##1{?}
   \expandafter\gdef\csname#1tag\endcsname{\relax}\fi\next}
\def\verbatimtags{\let\ift@gs=\iffalse\ifx\alltgs\relax\else
   \expandafter\t@gsoff\alltgs,\fi}
\catcode`\X=11 \catcode`\@=\active
\localtags
%
\global\newcount\numpunt
\hoffset=0.cm
\baselineskip=14pt  
\parindent=12pt
\lineskip=4pt\lineskiplimit=0.1pt
\parskip=0.1pt plus1pt

\hyphenation{small}

\catcode`\@=11

\catcode`\@=11


 \centerline {\bf  One-dimensional  random field Kac's model: }
\vskip.3truecm
\centerline { {\bf weak large deviations principle} \footnote* 
   {\eightrm Supported by: CNR-CNRS-Project 8.005,
INFM-Roma;  MURST/Cofin  05-06; GDRE 224 GREFI-MEFI, CNRS-INdAM.}}
\vskip.5cm
 \centerline{
Enza Orlandi \footnote{$^1$}{\eightrm Dipartimento di Matematica,
Universit\'a di Roma Tre, L.go S.Murialdo 1, 00156 Roma, Italy.\hfill
\break orlandi@matrm3.mat.uniroma3.it} and  Pierre Picco 
\footnote{$^2$}{\eightrm CPT, UMR 6207, CNRS, 
Universit\'e de Provence, Universit\'e de la Mediterran\'ee,
Universit\'e du Sud Toulon Var,
Luminy, Case 907, 13288
Marseille Cedex 9, France. picco@cpt.univ-mrs.fr }  } 
\footnote{}{\eightrm {\eightit AMS 2000 Mathematics Subject Classification}:
Primary 60K35, secondary 82B20,82B43.}
\footnote{}{\eightrm {\eightit Key Words}: phase transition, large deviations
random walk, random environment, Kac potential.
 }

\vskip.5cm
\centerline{ \it Dedicated to A.V. Skorohod for the fiftieth birthday  of
its fundamental paper  [\rcite{Sk}].}

\vskip .5cm
{ \bf Abstract}
We  prove a quenched weak large deviations principle for the  Gibbs measures of
a Random Field Kac Model (RFKM) in one dimension. The  external   random    magnetic
field  is  given by   symmetrically distributed  Bernoulli random
variables.
The results are  valid   for values of the
temperature,  $\b^{-1}$,  and  magnitude, $\theta$,  of the field in the region 
where   the free energy of the corresponding random  Curie Weiss
model  has only  two  absolute minima $m_\b$ and $Tm_\b$.
We give an explicit representation of the rate functional which is a  positive random functional 
determined by two distinct  contributions.  One  is
related to  the free energy cost $\FF^*$ to undergo a  phase change  
(the surface tension). The
$\FF^*$   is  the cost of one single phase change and  depends on the temperature and
magnitude of the field.
The other  is a   bulk contribution due to the presence of  the random magnetic field. We  
characterize the  minimizers of this random functional. 
We show that they  are step functions taking 
values  $m_\b$ and $Tm_\b$.   The points of discontinuity   are described by    a 
stationary renewal process  related to the $h-$extrema for 
a bilateral Brownian motion studied by Neveu and Pitman,
where $h$ in our context is a suitable  constant depending on
 the temperature and on  magnitude of the random field.   As an outcome we   have a   complete  
characterization of  the typical  profiles of   RFKM (the ground states) which   was  initiated in
[\rcite{COP1}]    and extended  in [\rcite{COPV}].    

\bigskip \bigskip
\chap {1 Introduction}1
\numsec= 1 \numfor= 1 \numtheo=1

\medskip

We consider a  one-dimensional spin system interacting via a ferromagnetic two-body Kac
potential  and  external random magnetic field given by symmetrically distributed
Bernoulli random variables. Problems where  a stochastic contribution is added to the energy of the
system arise naturally in condensed matter physics where the presence of the impurities
causes the microscopic structure to vary from point to point. Some of the vast literature on these
topics may be   found consulting  [\rcite {Ah}-\rcite {APZ}], [\rcite {B}], [\rcite {BRZ}], [\rcite {BK}], [\rcite
{FFS}- \rcite {IM}], [\rcite {Ku}], [\rcite {SW}].

 Kac's potentials
is a short way to  denote two-body ferromagnetic interactions with
range $\frac 1\g $, where $\g$ is a dimensionless parameter such
that when $\g \downarrow 0$, i.e. very long range, the strength of the
interaction becomes very weak keeping  the total
interaction between one spin and all the others finite.  They
were  introduced in [\rcite {KUH}], and then generalized in
[\rcite {LP}] and   [\rcite {PL}] to present a rigorous validity of the van der Waals
theory of a liquid-vapor  phase transition.  Performing first the
thermodynamic limit of the spin system interacting via Kac's
potential, and then the limit of infinite range, $\g \downarrow 0$, { Lebowitz and Penrose } 
rigorously derived the Maxwell rule, {\it i.e}    the canonical free
energy  of the system is the convex envelope of the corresponding canonical
free energy for the Curie--Weiss model. The consequence is that, in
any dimension, for values of the temperature at which the free
energy corresponding to the Curie-Weiss model is not convex, the  canonical
free energy  of the Kac's model  is not differentiable in the limit $ \g \downarrow 0$. 
These results show that 
long range models  give satisfactory answer for canonical  free
energies.
At the level of Gibbs measures 
the analysis is more delicate since the behavior of Gibbs measures
depends strongly on the dimension.
 
There are several papers trying 
to understand qualitatively
and quantitatively how a refined analysis of the Gibbs measures of the Kac models allows to see some
 features of systems with long, but finite range interaction,  see for instance
[\rcite {CP}], [\rcite {LMP}],   [\rcite {BZ}].

   For $\g$ fixed and different from zero,  if $d=1$, there exists an unique Gibbs state  for the
Kac model while for the Curie--Weiss model the measure induced by the
empirical magnetization weakly
converges, when
the number of sites goes to infinity,  to a convex contribution of two
different Dirac measures.    
In the one dimensional case, the analysis [\rcite {COP}] for Ising
spin and [\rcite {Bo}]  for more general spin, gives  a
satisfactory description of the typical profiles.  In these papers a large deviations principle for
Gibbs measures was  established.  The ground state of the system in
suitable chosen mesoscopic scales, is concentrated sharply near the two values of the minimizers
of the corresponding Curie-Weiss canonical  free energy.   The typical
magnetization profiles are
constant near
one of the two values over lengths of the order
$e^{\frac \b \g   F} $ where
$F$  was explicitly computed and represents the cost in term of  canonical free energy to go from one phase to the
other, i.e the surface tension.  
Moreover, suitably marking the locations of the   phase changes of  
the typical profiles and scaling the space by $e^{-\frac{\b}{ \g} F}$,
one gets as limiting Gibbs distribution of the marks, the one  of   a
Poisson Point Process. The thermal fluctuations are responsible  for 
the stochastic behavior on this scale.

 The
same type of questions could be asked for   the RFKM which  is one of  the simplest disordered spin system. 
This motivated the [\rcite
{COP1}],    [\rcite{COPV}] as well as the present paper.   The answers we found, as explained below, are dramatically  different from the ones obtained without the presence of the random field.  The
 analysis     done  holds   in dimension  $d=1$ and  for values of the temperature
and magnitude of the field in the whole region of two  absolute minima
for the canonical free energy of the
corresponding Random Field Curie Weiss model. This region is denoted $\EE$, see   \eqv (11.2500) for the precise
definition.  In the first paper  [\rcite
{COP1}] we gave the results   for $(\b,\th)$ in a subset of $\EE$, 
under some   smallness condition, whereas in    [\rcite
{COPV}] as well as in this paper we give the result for $(\b,\th)$ in $\EE$
without further constraints.  We will comment later about this, but one should bear in mind that the
results proven in   [\rcite {COP1}]  hold   for almost all
realizations of the random magnetic fields,
the ones proven in   [\rcite
{COPV}]  hold for a set of realizations of the random magnetic fields
of probability that goes to one when $\g\downarrow 0$, while the
ones in the present paper hold  merely in law.

   Let us recall the previous results:  Here, as well in the previous papers,  the first
step is a coarse graining procedure.  Through a
block-spin transformation, the microscopic system is mapped into a
system on $ \TT  = L^{\infty}(\R, [-1,1]) \times L^{\infty}(\R, [-1,1])$, see \eqv (1.13),  for which the
length of interaction becomes of order one (the macroscopic
system).  The macroscopic state of the system is determined by an
order parameter which specifies the phase of the system.
    It has been proven in [\rcite {COP1}] that for almost all
realizations   of the random magnetic fields, for intervals
whose length in macroscopic scale is of order $ (\gamma
\log\log (1/\gamma))^{-1}$ the typical block spin profile is either rigid,
taking one of the two values ($m_\b$ or $Tm_\b$)  corresponding to the minima of the canonical 
free energy  of  the random field Curie Weiss model, or makes at
most one transition from one of the minima to the other.
In the following,  we will denote these two minima  the  $+$ or $-$ 
phases.  It was also proven in  [\rcite {COP1}],  that if   the system is considered on an interval of
length $  \frac 1 \gamma (\log \frac 1 \gamma)^p$, $p\ge 2$,  the
typical profiles are not rigid over any interval of length  larger  or
equal to 
$L_1(\gamma)=\frac 1 \gamma (\log \frac 1 \gamma)(\log\log
\frac 1\gamma)^{2+\rho}$, for any $\rho>0$.

In [\rcite {COPV}] the following was proved:
On a set of realizations of the random field of overwhelming probability 
(when $\g\to 0$) it
is possible to construct random intervals of length of order
$\frac1 \g$  (macro scale)  and to associate a    random   sign in
such a way that,
typically with respect to the Gibbs measure,  the
magnetization profile is rigid on these intervals and, according to
the sign,
it belongs to the $+$ or $-$ phase.
Hereafter, ``random'' means that it depends on the realizations of the
random fields (and on $\b,\th$).
A description of the transition from one phase to the other was
also discussed in [\rcite{COPV}].  We  recall these results   in Section 2.
The main problem in the proof of the previous results is the ``non
locality'' of the system, due to the presence of the random field.
There is  an interplay between the ferromagnetic two-body
interaction which attracts spins alike  and the  presence of the
random field which would like to have the spins aligned according to its
sign. It is relatively easy to see that the fluctuations of the random field
over intervals in macro scale $\frac 1 \g$ play an important role.  To
determine the beginning and the end of the random interval
where the profiles are rigid and the
sign attributed to it, it is essential to verify  other  local
requirements for the random field. We need a detailed
analysis of suitable functions of the random fields
in all subintervals of the   interval of
order  $\frac 1 \g$. In fact,  it could happen that even though at large
the random fields undergoe  to a positive (for example) fluctuation,
locally there  are negative fluctuations which make not
convenient (in terms of the cost of the total free energy) for the system to have a
magnetization profile close to the + phase in that interval.

Another problem in  the previous  analysis is due to the fact that  the
measure induced by the block-spin transformation contains
multibody interaction of arbitrary order. Estimated roughly as in [\rcite {COP1}], this would
give a contribution proportional to the length of the interval in
which the transformation is done, there the length of intervals was
$(\g\log\log(1/\g))^{-1}$ and the $(\log\log(1/\g))^{-1}$ help us to
get a small contribution.
Here  we are interested in intervals of length $\frac 1 \g$. Luckily enough,
exploiting the randomness of the one body interaction,   it is enough to  estimate  the Lipschitz
norm of the multibody potential. Using cluster expansion tools,
this can be estimated through the representation of the multibody
interaction as an absolute convergent series.

 In this paper we first extend the results of [\rcite{COPV}] by
defining a random profile $u^*_\g$ which belongs to 
$BV([-Q(\g),Q(\g)],\{m_\b,Tm_\b\})$, the set of function from $[-Q(\g),Q(\g)]$ to $
\{m_\b,Tm_\b\}$ having bounded variation. Here $Q(\g)\uparrow
\infty$  when $\g\downarrow 0$ in a convenient  way. On a probability subspace of the
random magnetic field  configurations of overwhelming probability, we identify
a suitable neighborhood of $u^*_\g$ that has  a overwhelming Gibbs measure. 

Then we prove that when $\g\downarrow 0$ the limiting distribution of
the interdistance between the jump points of $u^*_\g$ with respect to the distribution of the
random magnetic fields is the  Neveu-Pitman [\rcite{NP}]
stationary renewal process
of $h$-extrema of a bilateral Brownian motion.  The value of $h$   depends on $\b$ and $\th$. 
Surprisingly  the  residual life distribution  of the renewal process that we   obtained is
the same (setting $h=1$) of  the one determined independently by
Kesten  [\rcite {Ke}]   and Golosov  [\rcite {Go}] 
representing  the limit distribution  of the point of localization of Sinai's
random walk in random environment, see Remark  \eqv (2a), in Section 2.  

This allows us to define the limiting (in Law) typical profile $u^*$ that
belongs to $BV_{{\rm loc}}(\R,\{m_\b,Tm_\b\})$, the set of functions
from $\R$ to $\{m_\b,Tm_\b\}$ that have bounded variations on each
finite interval of $\R$.  The total variation of $u^*$ on $\R$ is infinite. 

Note that here, the Gibbs measure is strongly concentrated on a random
profile that we relate to a renewal process, the randomness being the
one of the random magnetic fields. The phase change   of this   random
profile occurs on such a small scale that we cannot see the thermal
fluctuations that were responsible in the case without magnetic field
of the previously described  Poisson Point Process. At the  same scale
where we find the  renewal process, the system without magnetic
fields is  completely rigid, constantly equal to $m_\b$ or
$Tm_\b$.
Having  exhibited the typical profile  $u^*_{\g}$ and its limit in Law $u^*$
the next natural question concerns  the large deviations with respect
to this typical profile. 
Formally we would like to determine  a  positive functional $\G (u) $
for $ u \in \AA$, where $\AA \subset \TT $,    so that 
$$   \mu^{\o}_\g  [ \AA ] \sim \exp \{ - \frac {\b} \g \inf_{u \in \AA} \G ( u ) \}. \Eq (C.1)$$
When $ \AA \equiv \AA(u) \subset \TT  $ is a convenient, see
\eqv(D.30),  
neighborhood  of $ u \in BV_{\rm {loc} } (\R,  \{ Tm_\b,m_\b\})$ and
$u$ is a suitable local
perturbation of  the typical profile  $u^*_\g$, \eqv (C.1)    should be understood
as 
$$ \lim_{\g \downarrow 0} \left [ - \frac \g \b \log \mu^{\o}_\g  [
\AA (u)] \right ]  
=\G ( u ).  
 \Eq (C.2)$$ 

One has to give a probabilistic sense to the above convergence. It
appears that contrarily to the large deviation functional
associated to the global empirical magnetization (the canonical free
energy of the RFKM, see \eqv(2.14000)) which is not random, $\G(u)$ is random 
and the above convergence holds  in Law. In fact $\G(u)$ can be
expressed in term of $u$, the limiting
$u^*$ and the bilateral Brownian motion.  It represents, in the chosen  limit,  the random cost for the system
to deviate from the equilibrium value $u^*$.   The interplay between the surface free energy
$\FF^*$ (the cost  of   one single phase change)  and the random  bulk contribution appears in a rather clear
way.   Note that in \eqv(C.2) the functional is evaluated at $u$ even if the
considered neighborhood of $u$ does not shrink when $\g\downarrow 0$ to $u$.
This fact allows us to avoid to face difficult measurability problems
when performing infimum over family of sets.
The random functional $\G$ in \eqv(C.2) could be seen as the `` De Giorgi  Gamma-limit in
Law''  
for a sequence of 
intermediate random functionals   
obtained through a coarse graining procedure over $\TT$.
Since a precise definition of such a convergence is beyond the
scope of the paper, see however [\rcite{DMM}], and presents more
complications than simplifications we will not pursue it here. 

\vskip1.cm

The plan of the paper is the following. In Section 2 we give the
description of the model and present the main results. In Section
3 we   recall the  coarse graining procedure.  
In section 4 we prove  the main estimates to derive upper and lower
bound to deduce the  large deviation estimates. 
In Section 5 we prove the above mentioned convergence in Law of the
localization of the jumps of $u^*_\g$ to the stationary renewal
process of Neveu--Pitman. 
In section 6 we give the proof of the main results.

\bigskip

\bigskip
\vskip .5cm
 \noindent{\bf Acknowledgements}
 We are indebted to Errico Presutti who  gave us years ago the
expression  of the random functional, see \eqv(D.5a).  We   thank Jean--Fran\c{c}ois Le Gall
for mentioning to us  the  article by  Neveu and Pitman, Jean Bertoin,
Zhang Shi  and Isaac Meilijson for  illuminating discussions.

\vskip 1truecm
\chap{2 Model, notations  and main results }2
\numsec= 2
\numfor= 1
\numtheo=1 

\medskip
\noindent{\bf 2.1. The model}
\medskip
Let $(\O,\AA,\P)$ be a probability space on which we  define 
$h \equiv \{h_i\}_{i\in \Z}$, a family of
independent, identically distributed Bernoulli random variables with
$ \P[h_i=+1]=\P[h_i=-1]=1/2$. They represent random signs of external
magnetic fields
acting on a spin system on $\Z$, and whose magnitude  is denoted by
$\th>0$.
The configuration space is $\SS\equiv  \{-1,+1\}^\Z$.
If $\s \in \SS$ and $i\in \Z$,
$\s_i$ represents the value of the spin at site $i$. 
The pair interaction among spins is given by a Kac potential of the form
$J_\g(i-j)\equiv \g J(\g (i-j)),\, \g>0$.
We require that for $r \in \R$: (i) $J(r)\geq 0$ (ferromagnetism);
 (ii) $J(r)=J(-r)$ (symmetry); (iii) $J(r)\le ce^{-c'|r|}$ for $c,c'$
positive constants
(exponential decay); (iv) $\int J(r)dr=1$ (normalization).
For  sake of  simplicity we fix $J(r)= \1_{ [|r|\le 1/2]} (r)$,  where
we denote by
$\1_{A} (\cdot)$  the indicator function of the set $A$.

For $\L \subseteq \Z$ we set $\SS_\L=\{-1,+1\}^\L$; its elements
are   denoted by $\s_\L$; also,
if $\s \in \SS$, $\s_\L$ denotes its restriction to $\L$. Given
$\L\subset \Z$  finite and a realization of the magnetic fields,
the   Hamiltonian in the  volume $\L$, with free boundary conditions,
is the  random variable on $(\O,\AA,\P)$
given by
$$
H_\g(\s_\L)[\o]= -\frac 12 \sum_{(i,j) \in \L \times \L}
 J_\g(i-j) \s_i \s_j -\th\sum_{i\in \L} h_i[\o]\s_i. 
\Eq(2.1)
$$
 In the following we drop
the $\o$ from the notation. 

The corresponding {\sl Gibbs measure} on the finite volume $\L$,
at inverse temperature $\b>0$ and free
boundary condition  is then a random variable with values
on the space of probability measures on $\SS_\L$. We denote it by 
$\mu_{\b,\th,\g,\L}$  and it is defined by
$$
\mu_{\b,\th,\g,\L}(\s_\L)
= \frac 1{Z_{\b,\th,\g,\L}} \exp\{-\b H_\g(\s_\L)\} \quad \quad
\s_\L \in \SS_\L,
\Eq(2.3)
$$
where  $Z_{\b,\th,\g,\L}$ is the normalization factor 
called  partition function.
To take into account the interaction between the spins in $\L$ and
those outside $\L$ we set
$$
W_{\g}(\s_\L,\s_{\L^c}) = -\sum_{i\in \L} \sum_{j\in \L^c} J_\g(i-j)
\s_i \s_j.
\Eq(2.2)
$$
If $\tilde \s \in \SS$, the Gibbs
measure on the finite volume $\L$ and boundary condition
$\tilde \s_{\L^c}$ is  the  random probability measure on $\SS_\L$,
denoted by $\mu_{\b,\th,\g,\L}^{\tilde \s_{\L^c}}$  and
defined by
$$
\mu_{\b,\th,\g,\L}^{\tilde \s_{\L^c}}(\s_{\L})=
\frac 1{Z_{\b,\th,\g,\L}^{\tilde \s_{\L^c}}}
\exp\left\{- \b (H_{\g}(\s_{\L})+ W_{\g}(\s_{\L},\tilde
\s_{\L^c}))\right\},
\Eq(2.300)
$$
where again the partition function $Z_{\b,\th,\g,\L}^{\tilde
\s_{\L^c}}$ is  the normalization factor.

Given a  realization of $h$ and   $\g>0$, there is a unique
weak-limit of
$ \mu_{\b,\th,\g,\L}$ along a family of volumes $\L_{L}= [-L,L] \cap \Z$,  $  L
\in \N $; such limit
is called the infinite volume Gibbs measure
$\mu_{\b,\th,\g}$. The limit does not depend on the boundary
conditions, which may be taken
$h$-dependent, but it is a random element, i.e.,
different realizations of $h$ give a priori
different infinite
volume Gibbs measures.

\medskip
\noindent{\bf 2.2. Scales }
\medskip
 When dealing with local long range interaction,  as we did in 
[\rcite {COP}],  [\rcite {COP1}] and  [\rcite {COPV}],   
the  analysis of the configurations that are typical  
for $\mu_{\b,\th,\g}$ in the limit 
 $\g \downarrow  0$,
involves a block spin transformation  which
transforms the  microscopic system on $\Z$
in a 
system on $\R$.   Such changes of scales are standard in Kac type problems. 
Here, notations are particularly troublesome  because we have three
main different scales and according to the
case it is better to work with one or the other.
There will be also
intermediate scales  that we will discuss later.
For historical
reasons the three main scales are called: 
microscopic, macroscopic and  Brownian scale.  More properly they
should be denoted microscopic, mesoscopic and    
macroscopic.
Since  in the previous papers,  [\rcite {COP1}] and  [\rcite
{COPV}], the intermediate scale was 
called {\it macroscopic},  we
continue  to   call  it in such a way to avoid confusion.
Then we will   call  mesoscopic   scales  
all the intermediate scales  between the microscopic and macroscopic scales. 
These   mesoscopic  scales are not  intrinsic to the system but superimposed to study it.

\medskip
 \noindent \item{$\bullet$} { \it The microscopic and macroscopic scales. } 
\medskip
The basic space is the ``microscopic space'', i.e. the lattice $\Z $
whose elements are denoted by $i,j$ and so on.
The microscopic scale 
corresponds to 
the length measured according to the lattice
distance.
The spin $\s_i$ are indexed by $\Z$ 
and the range of interaction in this scale  is of order $\frac 1 \g$.  

The macroscopic regions correspond to intervals of $\R$  that
are of order $\frac 1 {\g} $ in  the  microscopic  scale 
; {\it i.e.}
if $ I \subset \R $,  is an interval in the macroscopic scale then it will
correspond to the interval $\frac {I} {\g} $ in the microscopic scale. 
Since the   range of the interaction  
is  of order $\g^{-1}$ in   the   microscopic  scale, in     the macroscopic
scale it  becomes of order 1.

\medskip
\noindent \item{$\bullet$} { \it The Brownian scale}
\medskip

  The Brownian scale is linked to the random magnetic fields.
The Brownian  regions correspond to intervals of $\R$  that are of
order $\frac 1 {\g^2} $ in the    microscopic 
scale; i.e. if $[-Q,Q] \subset \R $, $Q>0$ is an interval in  Brownian  scale then
it will correspond to   $[-\frac {Q} {\g^2},\frac
{Q} {\g^2}] $ in the   microscopic scale.   
In the Brownian scale the range of
interaction is of order $\g$.

\medskip
  
\medskip
 \noindent \item{$\bullet$} { \it The  partition of $\R$. } 
\medskip
 Given a rational positive number  $\d$, 
 $\DD_\d$ denotes the partition of $\R$ into  intervals $\tilde A_\d(u)=[u\d, (u+1)\d )$
  for   $u\in \Z$.  If $\d =n\d'$  for some  $n\in \N$,  then $\DD_\d$
 is coarser than $\DD_{\d'}$.  
A function $f(\cdot)$  on $\R$ is $\DD_\d$--measurable if it is
constant on each interval of $ \DD_\d$.  A region $ \L$ is
$\DD_\d$--measurable
if its   indicator function is  $\DD_\d$--measurable. 
For $r \in \R$,
we denote by $D_{\d}(r)$ the interval of $\DD_\d$  that contains
$r$. Note that for any  $r \in [u\d, (u+1)\d )$, we have that $D_{\d}(r)= \tilde A_\d(u)$.
  To avoid rounding problems in the following,   we will consider
intervals  that are always   $\DD_\d$--measurable. 
 If $I\subseteq\R$  denotes a macroscopic interval
we set 
$$
{\cal C}_{\d}(I)=\{u \in\Z; \tilde A_\d(u ) \subseteq I\}.
\Eq(PJG.1)$$

\medskip
 \noindent \item{$\bullet$} { \it The mesoscopic   scales  } 
\medskip

The smallest mesoscopic scale involves a parameter $0<\d^*(\g)<1$ satisfying
certain conditions of smallness that will be fixed later.
However we assume that $\d^*\g^{-1}\uparrow \infty$ when $\g \downarrow 0$.
The elements of $\DD_{\d^*}$ will be denoted by
$\tilde A(x)\equiv [x\d^*,(x+1)\d^*)$, with $x \in \Z$.
The partition $\DD_{\d^*}$ 
induce a partition of $\Z$ into blocks $A(x)=\{ i \in \Z ; i \g \in \tilde
A(x)\}\equiv\{a(x),\dots,a(x+1)-1\}$ with length of order
$\d^*\g^{-1}$ in the microscopic scale.

For notational
simplicity, if no confusion arises, we omit  to write the explicit dependence on 
  $\g,\d^*$.
 To avoid rounding problems, we assume  that $\g=2^{-n}$ for some integer $n$, with
$\d^*$ such that $ \d^* \g^{-1}$ is an integer, so that $a(x)= x \d^*
\g^{-1}$, with $x \in \Z$.
When considering another mesoscopic scale, say $\d>\d^*$, we
always assume that $\d^{-1}\in \N$ and $\d=k\d^*$ for some integer
$k\ge 2$.

\medskip
 \noindent   { \bf 2.3  Basic Notations. } 

\medskip
  \noindent \item{$\bullet$} { \it block-spin magnetization  } 
\medskip
Given a realization of $h$ and for each configuration
$\s_\L$, we could   have   defined for each block $ A(x)$ 
a pair of numbers where the first is 
the average magnetization over the sites with positive $h$ and the second to 
those with negative $h$. However it appears, [\rcite{COP1}], to be more 
convenient to use another random partition   of  $ A(x)$
into two sets of the same cardinality. 
This allows to separate on each block the expected contribution
of the random field from its local fluctuations. More precisely we have the following.

Given a realization $h[\o]\equiv (h_i[\o])_{i\in \Z}$, we set
$ A^+(x)= \big\{ i \in A(x); h_i[\o]=+1 \big\}$
and  $A^-(x)= \big\{ i \in A(x); h_i[\o]=-1 \big \}$. 
Let $ \l(x)\equiv {\rm sgn}( |A^+(x)|-(2\g)^{-1}\d^*)$, where ${\rm sgn}$
is the sign function, with the convention that ${\rm  sgn}(0)=0$. For convenience
we assume  $\d^* \g^{-1}$ to be even, in which  case:
$$
\P\left[ \l(x)=0 \right]= 2^{-\d^* \g^{-1}}
{{\d^* \g^{-1}}\choose {\d^* \g^{-1}/2}}.
\Eq(2.7)
$$
We note that  $\l(x)$ is a symmetric random variable.
When $\l(x)=\pm 1$ we set
$$
l(x)\equiv
\inf \{l \ge  a(x)\, : \sum_{j=a(x)}^l \1_{\{ A^{\l(x)}(x)\}}(j)\geq
\d^* \g^{-1}/2\}
\Eq(2.8)
$$
and consider the following decomposition of $A(x)$:
$B^{\l(x)}(x)= \left\{ i\in A^{\l(x)}(x); i \le l(x) \right\}$
and $ B^{-\l(x)}(x)=A(x)\setminus B^{\l(x)}(x)$.
When $\l(x)=0$ we set $B^+(x)=A^+(x)$ and $B^-(x)=A^-(x)$.
We set $ D(x)\equiv A^{\l(x)}(x)\setminus B^{\l(x)}(x)$. In this
way,   the set $B^{\pm}(x)$ depends on  the 
  realizations of the random field, but  the cardinality 
$|B^{\pm}(x)|=\d^* \g^{-1}/2$  is the same for all realizations.  
Set
$$
  m^{\d^*}(\pm,x, \s)=\frac {2\g}{\d^*} \sum_{i\in B^{\pm}(x)} \s_i.
\Eq(2.10)
$$
We   call  block spin  magnetization of the block $A(x)$  the vector  
$$
  m^{\d^*}(x, \s)=(  m^{\d^*}(+,x, \s),  m^{\d^*}(-,x, \s) ). 
\Eq(2.10a)
$$
The  total empirical magnetization of the block $A(x)$  is, of course, given by 
$$
\frac \g{\d^*} \sum_{i\in A(x)}\s_i=
 \frac 12 (m^{\d^*}(+,x,\s)+m^{\d^*}(-,x,\s))
\Eq(2.12)
$$
and   the contribution of the magnetic field  to the Hamiltonian \eqv(2.1) is 
$$
\frac \g{\d^*} \sum_{i\in A(x)} h_i\s_i =\frac 12
(m^{\d^*}(+,x,\s)-m^{\d^*}(-,x,\s)) +
\l(x) \frac {2\g}{\d^*} \sum_{i\in D(x)} \s_i .
\Eq(2.13)
$$

\medskip
  \noindent \item{$\bullet$} { \it spaces of the   magnetization profiles } 
\medskip
Given a  volume $\L\subseteq \Z$ in the original microscopic spin
system, it corresponds to the macroscopic volume $I=\g \L=\{\g i;
i \in \L\}$, assumed to be  $\DD_{\d^*}$--measurable.
 The block spin transformation, as considered in [\rcite {COP1}] and
[\rcite{COPV}],
is the random map which associates to the spin configuration $\s_\L$
the vector $(m^{\d^*}(x, \s))_{x \in {\cal C}_{\d^*}(I)}$,  see \eqv (2.10a),  
with
values in the set
$$
 \MM_{\d^*} (I)\equiv \prod_{x\in {\cal C}_{\d^*} (I)}
\left  \{-1, -1+ \frac {4\g}{\d^*}, -1+ \frac {8\g}{\d^*},\dots,
1- \frac {4 \g}{\d^*},1 \right\}^2.
\Eq (2.14)
$$

We  use the same notation $\mu_{\b,\th,\g,\L}$
to denote
both, the Gibbs measure on
$\SS_\L$, and the probability measure   induced on $ \MM_{\d^*}(I)$, through the
block spin transformation, i.e.,
a coarse grained version of the original measure.
Analogously, the infinite volume limit (as $\L\uparrow\Z$) of the
laws    of the block spin 
$(m^{\d^*}(x))_{x \in {\cal C}_{\d^*}(I)}$
under the Gibbs measure
will also be denoted by  $\mu_{\b,\th,\g}$.

\medskip
\noindent 
We denote   a generic element in 
 $\MM_{\d^*}(I)$ 
  by 
 $$
 m^{\d^*}_I\equiv(m^{\d^*}(x))_{x \in {\cal C}_{\d^*}(I)}
\equiv (m^{\d^*}_1(x),m^{\d^*}_2(x))_{x \in {\cal C}_{\d^*}(I)}.  \Eq (2.14a)
$$ 
Since  $I$ is   assumed to be  $\DD_{\d^*}$--measurable, we can identify $m^{\d^*}_I$ with
the element of 
$${\cal T}= \{ m \equiv (m_1,m_2) \in L^{\infty}(\R) \times L^{\infty}(\R);
\|m_1\|_\infty\vee\|m_2\|_\infty\le 1\} \Eq (1.13) $$ 
 piecewise   constant,   equal   to $m^{\d^*}(x)$ 
on each $\tilde A(x)= [x\d^*, (x+1)\d^* ) $ for $x \in {\cal C}_{\d^*}(I)$, and 
vanishing   outside $I$. Elements of $\TT$ will be   called
magnetization profiles. 
Recalling that $I=\g \L$, the block spin transformation can be
identified with a map  
from the space of spin configurations
$\{-1,+1\}^{\L}$ (with $\L$ a microscopic volume)
into the subset of $\DD_{\d^*}$--measurable functions of 
$L^{\infty}(I) \times L^{\infty}(I) $ (with $I=\g\L$
a macroscopic  volume). 

For $\d\ge \d^*$, recalling that
$\forall r \in [u\d, (u+1)\d )$, we have $D_{\d}(r)= \tilde A_\d(u)$,
we  define for  $ m= (m_1,m_2) \in    {\cal T}$  
and  $i=1,2$ 
$$ 
m ^\d_i (r) =  \frac {1}{\d} \int_{D_{\d}(r)} m_i(s) {\rm}ds.  
\Eq (7.300) 
$$
This defines a map from $\cal T$ into the subset of $\DD_{\d}$--measurable
functions of $\cal T$. 
We define also a map from $\cal T$ into itself by 
 $$
 (Tm)(x)= (-m_2(x),-m_1(x))
 \quad \forall x\in \R. \Eq(4.5P)
 $$
In the following we denote  the total magnetization at  the site $x \in \R$ 
$$\tilde m (x)= \frac {m_1(x)+ m_2(x)} 2. \Eq  (mag) $$ 
\medskip
 \noindent \item{$\bullet$} {\it The Random Field Curie--Weiss model}
\medskip 
The  Lebowitz -Penrose theory, [\rcite{LP}],  is easy to prove for the 
Random Field Kac Model see [\rcite{COP1}], Theorem 2.2.
Namely,  performing first the
thermodynamic limit of the spin system interacting via Kac's
potential and then the limit of infinite range, $\g \to 0$, the canonical
free energy  of the Random Field Kac model  is the convex envelope of the corresponding
canonical free energy for the Random   Field   Curie-Weiss model.

  The canonical free energy for the Random   Field   Curie-Weiss model
derived in [\rcite{COP}]  is 
$$
f_{\b,\th}(m_1,m_2)=
-\frac{(m_1+m_2)^2}{8}-\frac{\th}{2}(m_1-m_2)+\frac{1}{2\b}(\II(m_1)+\II(m_2)),
\Eq (2.14000)
$$
where 
$\II(m)=\frac{(1+m)}{2}\log\big(\frac{1+m}{2}\big)+
\frac{(1-m)}{2}\log\big(\frac{1-m}{2}\big)$.
In Section 9 of  [\rcite{COPV}], it was proved that  
$$
\EE=\cases{
0<\th<\th_{1,c}(\b), &for $1<\b<\frac 32$;\cr
0<\th\le \th_{1,c}(\b)& for$ \b\ge \frac 32$,\cr}
\Eq(11.2500)
$$
where 
$$
\th_{1,c}(\b)=\frac 1 \b \,{\rm arctanh}\,(1-\frac 1\b)^{1/2},
\Eq(11.248)
$$
is the maximal  region of the two
parameters $(\b,\th)$, whose closure contains $(1,0)$  in which
$f_{\b,\th}(\cdot, \cdot)$  has exactly three critical points 
$m_\b,0,Tm_\b$. The two equal minima correspond to  
$m_\b=(m_{\b,1},m_{\b,2})$ and $Tm_{\b}=(-m_{\b,2},-m_{\b,1})$ and $0$ a local maximum. 
Calling  $\tilde m_\b = \frac { m_{\b,1}+m_{\b,2}}2$, on $\EE$ we have 
$$
\frac{\b}{2\cosh^2(\b(\tilde m_\b +\th)} +\frac{\b}{2\cosh^2\b(\tilde
m_\b-\th)}<1. 
\Eq(PP.1)
$$
Moreover, for all $(\b,\th) \in \EE$, the minima are quadratic and
therefore there exists a strictly positive constant
$\k (\b,\th)$ so that for each $m \in [-1,+1]^2$
$$
f_{\b,\th}(m)-f_{\b,\th}(m_{\b}) \geq \k (\b,\th)
\min\{\|m-m_\b\|^2_1,\|m-Tm_\b\|^2_1\}, \Eq(2.19)
$$
where
$ \| \cdot \|_1$ is  the $\ell^1$ norm in $\R^2$.

\medskip
 \noindent \item{$\bullet$} { \it The  spatially homogeneous phases  } 
\medskip 

We   introduce the so called ``excess  free energy
functional" $\FF(m)$,  $m \in {\cal T}$:
 $$
 \eqalign {& \FF(m)= \FF(m_1,m_2) \cr & = \frac 1 4 \int \int
J(r-r') \left [\tilde m(r)- \tilde m(r')\right ]^2 dr dr'+ \int
\left [ f_{\b,\th}(m_1(r),m_2(r)) - f_{\b,\th}(m_{\b,1},m_{\b,2})
\right] dr }\Eq (AP.1) $$   with
 $ f_{\b,\th}(m_1,m_2) $  given by \eqv (2.14000)
  and $\tilde
m(r)=(m_1(r)+m_2(r))/2$.
The functional $\FF$ is  well defined and non-negative,
although it may  take the value $+\infty$.
Clearly, the absolute minimum of $\FF$ is attained at the functions
constantly equal to $m_\b$ (or constantly equal to $Tm_\b$), the minimizers  of
$f_{\b,\th}$.  These two minimizers of $\FF$ are called the spatially homogeneous phases. 
The functional  $\FF$ represents the continuum approximation
of the deterministic contribution to the free 
energy of the system (cf. \equ(3.11P)) normalized by subtracting    
$f_{\b,\th}(m_{\b})$, the free  energy of the homogeneous phases. 
Notice that  $\FF$  is invariant  under
the  $T$-transformation, defined in  \eqv (4.5P).

 \medskip
 \noindent \item{$\bullet$} { \it The  surface tension } 
\medskip
   In  analogy to systems in higher  dimensions,
we denote by surface tension the  free energy cost 
needed by    the system to undergo to  a  phase change.
It has been proven in [\rcite {COP2}] that under 
the condition $m_1(0)+m_2(0)=0$, and for $(\b,\th) \in \cal E$, 
there exists a unique  minimizer 
$\bar m=(\bar m_1,\bar m_2)$, 
of $\FF$ over the set
$$
\MM_{\infty}=\{(m_1,m_2)
\in {\cal T}; \limsup_{r\to -\infty} m_i(r) < 0<\liminf_{r\to +\infty} m_i(r),
i=1,2 \}.
 \Eq(AP.5)
$$
Without the condition $m_1(0)+m_2(0)=0$, there is a continuum of minimizers 
obtained  translating $\bar m$. 
The minimizer 
$\bar m(\cdot )$ is infinitely differentiable 
and converges exponential fast, as  $ r \uparrow +\infty$ (resp. $-\infty$)
to the limit value $m_\b$, (resp.$Tm_\b$). 
Since   $\FF$ is invariant by  the $T$-transformation, see \eqv (4.5P), 
interchanging $ r \uparrow +\infty$ and $ r \downarrow -\infty$  in \eqv (AP.5),
there exists  one other family of minimizers
obtained translating $T \bar m$. 
We denote  by  $$\FF^*\equiv  \FF^*(\b, \th)=  \FF(\bar m)= \FF(T\bar m)>0,  \Eq (min1) $$ 
the surface tension. 
\smallskip

\medskip
 \noindent \item{$\bullet$} { \it how to detect local equilibrium  } 
\medskip

As in [\rcite {COP1}], the description of the profiles is based on
the behavior of local averages of
$m^{\d^*}(x)$ over $k$ successive blocks in the block spin
representation, where
$k\geq 2$ is a positive integer. Let $\d=k\d^*$   be  such that  $1/\d
\in \N$. Let $\ell \in \Z$,   $[\ell,\ell+1)$   be   a macroscopic block of length 1, 
$  {\cal C}_\d([\ell,\ell+1))$,    as     in \eqv (PJG.1), and $\z>0$. We define the   block spin  
variable
$$
\eta^{\d,\z}(\ell)=\cases{
1,& if $\,\forall_{u \in {\cal C}_{\d}([\ell,\ell+1))}\,
\frac {\d^*}{\d}\sum_{x\in {\cal C}_{\d^*}([u\d,(u+1)\d))}
\|m^{\d^*}(x,\s)- m_\b\|_1 \le\z$;\cr
 -1, &if $\,\forall_{u \in {\cal C}_{\d}([\ell,\ell+1))}\,
\frac {\d^*}{\d}\sum_{x\in {\cal C}_{\d^*}([u\d,(u+1)\d))}
\|m^{\d^*}(x,\s)-Tm_\b\|_1 \le\z$;\cr
0, & otherwise.\cr }
\Eq (2.190)
$$
where  for a  vector $v=(v_1,v_2)$, $\| v \|_1 = |v_1|+|v_2|$.
 When $
\eta^{\d,\z}(\ell)=1 $, (resp. $-1$), we 
say that a spin configuration $ \s 
\in
\{-1,1\}^{\frac 1 \g [\ell, \ell+1 ) } $ has magnetization   close 
to $m_\b$,
(resp. $Tm_\b$), with accuracy $(\d, \z)$  in $
[\ell, \ell+1 ) $.   Note that $
\eta^{\d,\z}(\ell)= 1 $ (resp $-1$) is equivalent to 
$$ \forall  y \in [\ell,  \ell+1 ) \qquad   \frac 1 \d 
\int_{D_\d (y)} dx \|  m^{\d^*} (x,\s) -v  
\|_1 \le \z \Eq (E.11) $$
for $v= m_\b$ (resp. $Tm_\b$), 
since for any $u \in {\cal C}_\d( [\ell, \ell+1 ))$, 
for all $y\in [u\d,(u+1)\d)\subset [\ell,\ell+1)$, 
$D^{\d}(y)=[u\d,(u+1)\d)$. 
We say that   a magnetization profile
$ m^{\d^*}(\cdot)$, in a macroscopic interval  $I \subseteq \R$, 
is close to  the equilibrium phase $\t$, for $\t\in \{-1,+1\}$,
with   accuracy $(\d,\z)$
when 
$$ \{ \eta^{\d,\z}(\ell) = \t,\, \forall  \ell \in I \cap \Z \}  \Eq (equi) $$
or equivalently
 if
$$  
\forall  y \in I  \qquad  \frac 1 \d \int_{D^\d(y)} dx \| m^{\d^*} (x,\s)-v\|_1 \le \z 
 \Eq (D.1a)
$$ 
 where $v=m_\b$ if $\t=+1$ and $v=Tm_\b$ if $\t=-1$. 
In the following  
the letter $\ell$ will   always   indicate an element of $\Z$.  This will allow to 
 write
\eqv (equi)  as 
  $ \{ \eta^{\d,\z}(\ell) = \t,\, \forall    \ell \in I  \}   $. 
  In   \eqv (D.1a)   the interval $I$ is always given in the macro--scale.
The definition   \eqv(D.1a) can be used   for function   $v$   
 more general
 than the constant ones. 
 In particular,  given 
$v=(v_1,v_2)
\in {\cal T} $,  $ \d=n\d^*$ for some positive integer $n$, $\z>0$,  
 and
 $[a,b)$ an interval  in   Brownian scale,  
we say that a spin configuration $ \s \in
\{-1,1\}^{[\frac {a} {\g^2}, \frac {b} {\g^2})} $ has magnetization
profile  
  close   to $v$ with 
accuracy $(\d, \z)$  in the interval $[a,b)$ if $\s $ belongs to the set  
$$ \left \{ \s \in 
\{-1,1\}^{[\frac {a} {\g^2}, \frac {b} {\g^2})}:
\forall  y\in [\frac {a} \g ,  \frac {b} \g )
\qquad   \frac 1 \d 
\int_{D^\d (y)} dx \|  m^{\d^*} (x,\s)   
-v^{\d^*}(x)
\|_1 \le \z \right  \}.  
\Eq (D.3)
$$

\vskip.5truecm

\noindent In view of  the results on the typical configurations
 obtained in [\rcite{COPV}] the above notion  is too strong. 
In  fact the typical   profiles  form
long runs of length of order $\g^{-1}$ (in the macroscopic scale) 
of $\eta^{\d,\z}(\cdot)=1$ that are followed
by short runs of $\eta^{\d,\z}(\cdot)=0$ that are in turn followed by
long runs of $\eta^{\d,\z}(\cdot)=-1$.
The typical profiles undergo to a   phase change within the runs of 
$\eta^{\d,\z}(\cdot)=0$ .
The length of these runs, see  Theorem  2.4  in 
[\rcite{COPV}], is smaller than  $2R_2=2R_2(\g)\uparrow
\infty$   in the macroscopic scale, see \eqv (erre2).
In the  Brownian scale, this length becomes $2\g R_2$ and one obtains that 
$\g R_2\downarrow 0$.
So in    Brownian scale, when $\g\downarrow 0$,  the
localization of  the phase   change    shrinks  to a point 
: the point of a jump.
For small $\g>0$,  the results in  [\rcite{COPV}]  allow  to   localize  these   points
within an interval of length $2\r>>2\g R_2$ centered around well
defined points  depending  on the realizations of the random field.
We call $\r$  the  fuzziness and $\r=\r(\g)\downarrow 0$ in the Brownian  scale.

With this in mind, a candidate
for  the limiting support of $\mu_{\b,\th,\g}$ when $\g\downarrow 0$
is an appropriate neighborhood 
of  functions on $\R$, (considered in the Brownian scale), 
taking two values $m_\b = (m_{\b,1}, m_{\b,2})$ or
$Tm_\b= (-m_{\b,2}, -m_{\b,1}) $
that have finite   variation. To fix the notations, we recall  the standard
definitions. 
Let us define, for any bounded interval $[a,b)\subset \R$
(in the Brownian scale)   $BV([a,b),\{m_\b,Tm_\b\})$  
as the set of   right continuous  bounded variation functions on $[a,b)$ with value in
$\{m_\b,Tm_\b\}$. Since we consider mainly  bounded variation functions
with value in $\{m_\b,Tm_\b\}$, we write $BV([a,b))\equiv BV([a,b),\{m_\b,Tm_\b\})$.
Since any bounded variation function $u$ is the difference of two
increasing functions, it has a left limit. 
We call the jump at $r$ the quantity $Du(r)=u(r)-u(r_{-})$ where
$u(r_{-})=\lim_{s\uparrow r}u(s)$.  If $r$ is such that  $Du(r)\neq 0$ we call $r$ a point of
jump of $u$,  and in such a case  $\|D u (r)\|_1=4 \tilde m_\b $.  We 
denote  by $ N_{[a,b)} (u)$   the number of jumps of $u$ on $[a,b)$
and  by
$V_a^b(u)$ the variation of $u$ on $[a,b)$, i.e. 
$$
V_{a}^{b} (u)\equiv  \sum_{a\le r<b } \|D u (r)\|_1 =  N_{[a,b)} (u)  2 [m_{\b,1}+ m_{\b,2}] =
4 \tilde m_\b N_{[a,b)} (u)   
<\infty . 
\Eq (D.4)
$$  
Note that $Du(r)\neq 0$ only on points of jump of $u$ and therefore the sum in
\eqv(D.4) is well defined. 
We denote by $BV_{{\rm loc}}
\equiv BV_{{\rm loc}} (\R, 
\{ m_\b, Tm_\b\})$   the set  of  functions  from
$\R$ with values in $\{m_\b,Tm_\b\}$
which  restricted to any  bounded
interval have bounded  variation
but not necessarily
having bounded variation on $\R$.
If  $u \in BV_{{\rm loc}} (\R,  \{ m_\b, Tm_\b\})$,
 then, see \eqv (mag),  $ \tilde u \in BV_{{\rm loc}}(\R,  \{\tilde m_\b, -\tilde m_\b\}
)$ where $\tilde m_\b$ is defined before \eqv(PP.1).

\vskip0.5cm

  \noindent 
  Since   a phase  change     can be better detected  in   macro units
we state the following definition which  corresponds to     Definition 2.3 of   [\rcite{COPV}].
\vskip0.5cm
  \noindent {\bf \Definition  (1)}   { \bf The macro interfaces } {\it 
Given an interval $[\ell_1,\ell_2]$ (in  macro-scale)  and a  positive integer
 $2R_2 \le  | \ell_2-\ell_1|$,
 we   say that a  single phase change   occurs within $[\ell_1,\ell_2]$
on a length $R_2$
if there exists $ \ell_0\in (\ell_1+R_2,\ell_2-R_2) $ so that $\eta^{\d,\z}
(\ell)=\eta^{\d,\z}(\ell_1)\in \{-1,+1\}, \forall \ell 
 \in[\ell_1, \ell_0 -R_2]$;
$\eta^{\d,\z}(\ell)=\eta^{\d,\z}(\ell_2)=-\eta(\ell_1), \forall  \ell \in
[\ell_0+R_2,\ell_2]$, and 
$\{\ell \in [\ell_0 -R_2,\ell_0+R_2]: \eta^{\d,\z}(\ell)=0\}$ is a set of
consecutive integers.
We denote by $ \WW_1([\ell_1,\ell_2], R_2,\z)$  the  set 
of   configurations $\eta^{\d,\z}$ with 
these properties.}
\vskip0.5cm
  \noindent In words, on $ \WW_1([\ell_1,\ell_2], R_2,\z)$, 
there is an unique
run of $\eta^{\d,\z}=0$, with no more than $2R_2$ elements,  inside
the interval $[\ell_1,\ell_2]$.  
To take into account that  for the typical profiles  
the point of jumps  are  determined with  fuzziness
$\r$, it is convenient to    associate to   $u  \in  BV ([a,b) ) $  a
partition of the interval $[a,b)$  (in Brownian scale) 
as follows : 
\vskip0.5cm
\noindent {\bf \Definition  (part.a)} $\, \,$ 
 {\bf Partition associated to BV functions}  {\it  
   Given $u \in BV([a,b))$, $\rho> \d=n\d^* $,     with $8\r+8\d$  
smaller than the minimal distance between two points of jumps of $u$, let    $C_i (u)$,
$i=1,..,N_{[a,b)}(u)$, (see \eqv(D.4)),  be  the smallest   $\DD_\d$ measurable
interval that contains an interval
of diameter $ 2 \r$, 
centered  at  the $i-$th   jump  of
$u$ in $[a,b)$. 
We have $C_i (u)\cap C_j (u) =\emptyset$ for $i \neq j$.

Let $C(u)=
\cup_{i=1}^{N_{[a,b)}(u)}C_i (u)$.
We set $B(u)=  [a,b) \setminus C(u) $ and   write $[a,b)  = C(u) \cup B(u) $. 
We denote by $C_{i,\g}(u)=\g^{-1}C_i(u)$,   $ C_\g(u)= \g^{-1}C(u)$
and  $ B_\g(u)= \g^{-1}B(u)$ the
elements of the  induced partition on the macroscopic scale. }
\vskip0.5cm

\noindent Whenever we deal with functions in
${\cal T }$ we   will always  assume that their  argument varies  on the
macroscopic  scale. So $m  \in {\cal T }$  means that $ m(x), x \in I $  where $ I \subset \R $ is an
interval in the macroscopic  scale. 
Whenever we deal with bounded variation functions, if not
further specified, we will always  assume that their  argument varies  on the Brownian
scale.  Therefore $u \in BV([a,b))$ means that $ u (r), r \in [a,b)$ and $[a,b)$
is considered in the Brownian scale.
This means that  in the macroscopic  scale we need to
write $ u (\g x )$  for $x \in [\frac a \g ,  \frac b \g)  $.
For $u \in BV([a,b))$,  we define for $x \in
 [\frac a \g ,  \frac b \g) $ {\it i.e} in the macroscopic scale, 
$$
u^{\g, \d^*}(x)=\frac{1}{\d^*} \int_{D_{\d^*}(x)} u(\g s)\,ds. 
\Eq(BVnew)
$$
 Given $[a,b)$ (in the Brownian scale),  $u$ in
$BV([a,b))$,   $ \r> \d=n\d^* >0$, with  $8\r+8\d$   
satisfying  the condition of  Definition \eqv(part.a),
$\z>0$, 
we say that a spin configuration $ \s \in
\{-1,1\}^{[ \frac a {\g^2} , \frac b{\g^2} )} $ has  magnetization
profile close   to $u$ with accuracy $(\d, \z)$ and {\it fuzziness} $\rho$ if
$\s \in
\PP^{\rho}_{\d,\g,\z,[a,b)} (u) $ where
$$ \eqalign {
&\PP^{\rho}_{\d,\g,\z,[a,b)} (u)= \cr
& \left \{ \s \in 
\{-1,1\}^{[ \frac a {\g^2}, \frac b {\g^2})}:   \forall  y \in B_\g (u), \,   \frac 1 \d 
\int_{D^\d (y)}  \|  
m^{\d^*} (x,\s)-   u^{\g,\d^*}( x)
\|_1\,dx \le \z   \right \}   
\bigcap_{i=1}^{N_{[a,b)}(u)}  \WW_1([C_{i,\g}(u)], R_2,\z).  
 } \Eq (D.30)
$$
In  \eqv (D.30) we consider the spin configurations
close  with accuracy $(\d, \z)$ to $m_\b$ or $Tm_\b$ in $ B_\g(u)$  according to the
value of $u^{\g, \d^*} (\cdot)$.  In  $C_\g(u)$ 
we require that the spin configurations have only one jump in  each
interval $C_{i,\g}(u)$,
$i=1,..N$,    and are close with accuracy  $(\d, \z)$ to the
right and to the left of this interval to the value of $u$ in those
intervals of $B_\g(u)$ that are adjacent to  $C_{i,\g}(u)$.  
 With all these definitions in hand we can slightly improve 
the main results of [\rcite{COPV}].
 \vskip0.5cm \noindent  
\noindent{\bf \Theorem (copv) [COPV]} {\it Given $(\b,\th) \in \EE$,
see \eqv(11.2500), 
there exists
$\g_0(\b,\th)$ so that for  $0<\g \le \g_0(\b,\th)
$,   for   
$Q = \exp [(\log  g(1/\g))/\log\log g(1/\g)]$,  with $g(1/\g)$ a
suitable positive, increasing 
function such that $\lim_{x\uparrow \infty} g(x)=+\infty$,
$\lim_{x\uparrow \infty} g(x)/x=0$  for suitable values of
$\d>\d^*>0$, $\r>0$, $\z>0$, $a'>0$, $R_2$   there exists    $\O_1  \subset \O $
with 
$$
\P[\O_1] \ge 
1-  K(Q)
 \left(\frac 1{ g({1}/{\g})}\right)^{a'}
\Eq(5.02)$$
where
$$
K(Q)= 2+ 5(V(\b,\th)/(\FF^*)^2)Q\log[Q^2g(1/\g)],
\Eq(M.1)
$$
 $\FF^*=\FF^*(\b,\th)$ is  defined in \eqv (min1)
and
$$
V(\b,\th) = \log \frac { 1+ m_{\b,2} \tanh (2 \b \th) } {  1- m_{\b,1}
\tanh (2 \b \th)}.
\Eq (C.1a)
$$
For  $\o \in \O_1$  we explicitly construct    $ u^*_\g (\o) \in BV ([-Q, Q])$
so that the minimal distance between jumps of $u^*_\g$  within $[-Q,+Q]$
is bounded from below by $8\r+8\d$,
$$
\mu_{\b,\th,\g}\Big(  \PP^{\rho}_{\d,\g,\z,[-Q,Q]} ( u^*_\g (\o)) \Big)
\ge 1 -2  K(Q)  e^{-\frac \b \g \frac 1 { g(1/\g) }},\Eq (copvprob)
$$
and
$$
V_{-Q}^{Q}(u^*_\g)\le 4\tilde m_\b K(Q).
\Eq(totalrecall)
$$

}

\vskip0.5cm
\noindent The previous Theorem is a  direct consequence of Theorem
2.1, Theorem 2.2
and Theorem 2.4 proven in  [\rcite {COPV}], together with 
Lemma \eqv (36)  that gives the   value   \eqv(M.1). The control of
the minimal distance  between jumps of $u^*_\g$ is done at the end of Section 5.

To facilitate the reading
we did not write explicitly  in the statement of Theorem \eqv(copv)
the  choice  done of the parameters $\d,\d^*,\z,g, R_2$ nor the explicit
construction of $u^*_\g$. We  dedicate the  entire Subsection 2.5  to recall
and motivate the choice of the parameters done in  [\rcite {COPV}] as
well as in this paper. The  $u^*_\g(\o)$  in   Theorem  \eqv (copv)  is a function  in
$BV ( [-Q,Q])$ associated to the sequence of maximal elongations
and their sign as determined  in  [\rcite
{COPV}] Section 5.  
For the moment it is
enough to know that it is possible to  determine  random points 
$ \a^*_i=\a^*_i(\g, \o)$ and a  random  sign $\pm 1$ associated to intervals
$[\e\a^*_i,\e\a^*_{i+1})$ in the Brownian scale, where $\e=\e(\g)$ has
to be suitably chosen. 
These random intervals  are the so called  maximal elongations. 
We   denote 
$$
u^*_\g (\o) (r) \equiv
\left \{
\eqalign { &  m_\b, \,\,\,\,\quad r \in
[\e\a^*_i, \e\a^*_{i+1})  \,\,\quad
\hbox {if the sign of
elongation $[\e\a^*_i, \e\a^*_{i+1})$ is } =+1
\cr & Tm_\b, \quad r \in [\e\a^*_i, \e \a^*_{i+1})   \quad  \hbox {
if the sign of elongation $[\e\a^*_i, \e \a^*_{i+1})$ is} =-1.
 } \right.     \Eq (L.1)
$$
for $ i \in  \{\k^*(-Q)+1,\dots,-1,0,1, \dots, \k^*(Q)-1\}$, 
where 
$$ \k^*(Q)=\inf(i\ge 0: \e\a^*_i>Q),  \qquad 
 \k^*(-Q)=\sup(i\le 0: \e \a^*_i<-Q) \Eq (num)  $$
with the convention  
 that $\inf(\emptyset)=+\infty,
 \sup(\emptyset)=-\infty$ and  $\e\a^*_0<0$ and $\e\a^*_1>0$, that is just a
relabeling of the points determined in [\rcite {COPV}],  Section 5.
The $\k^*(-Q)$ and $\k^*(Q)$ are random numbers and 
Lemma \eqv(36) gives that, with a $\P$--probability absorbed  in
\eqv(5.02), we have 
 $ |\k^*(-Q)|\vee \k^*(Q) \le K(Q)$, with $K(Q)$  given in  \eqv(M.1). 
This implies \eqv(totalrecall)

\medskip
\noindent{\bf 2.4. The main results}
\medskip
Let   $u$ and $u^* \in  BV_{{\rm loc}} $.   Denote by  $(W(r), r\in \R)$  the    Bilateral Brownian  motion
(BBM) , {\it i.e.}  the
Gaussian process   with independent
increments that satisfies  $\E(W(r))=0$, $\E(W^2(r))=|r|$ for $r \in \R$  (and
therefore $W(0)=0$) and  by $ \PP$ its    Wiener   measure
on $  C \left (\R , \BB(C(R))\right ) $.  
Let   $W$ be
a real valued continuous function  from
$\R$ to $\R$,   that is a realization of  a   BBM.
  Let   $[a,b)\subset \R$ be a finite interval    and 
denote  by $ \NN_{[a,b)}(u, u^*) $ the  points of jump of $u$ or $u^*$ in $[a,b)$: 
$$ \NN_{[a,b)}(u, u^*) = \{ r \in [a,b):  \|D u
(r)\|_1 \neq 0\,\,\hbox {or}  \,\,  \|D u^*
(r)\|_1 \neq 0 \}.  \Eq (G.1)$$
Note that by right continuity if $\|D u(a)\|_1 \neq 0 $ then   $ a \in  \NN_{[a,b)}(u, u^*)$,
while if  $\|D u (b)\|_1 \neq 0 $ then $b \notin   \NN_{[a,b)}(u, u^*)$. 
Since  $u$ and $u^*$ are   $  BV_{{\rm loc}}$ functions $
\NN_{[a,b)}(u, u^*)$ is a finite set of points.
We    index  in increasing order the points
in $ \NN_{[a,b)} (u, u^*)$ and  by an abuse of notation we  denote 
$\{ i \in  \NN_{[a,b)} (u, u^*)\}$ instead of  $\{ i : r_i \in  \NN_{[a,b)} (u, u^*)\}$. 
     Define  for  $u \in BV_{{\rm loc}} $,  the following  finite volume functional 
$$ \eqalign { &  \G_{[a,b)}  (u|u^*, W)\cr &=  \frac 1 {  2 \tilde m_{\b}}
   \sum_{ i\in  \NN_{[a,b)}(u, u^*) }    \left  \{ \frac {\FF^*} 2  \big[   \|D u
(r_i)\|_1
-     \|D u^*(r_i)\|_1  
\big ]    -
V (\b,\th)
( \tilde u(r_i)-  \tilde u^*(r_i))[ W (r_{i+1} )- W (r_i)]  \right \}  .}\Eq
(D.5)
$$
 The functional in \eqv (D.5) is always well defined since it is sum of finite terms. 
 In the following   $u^* \equiv  u^* (W) $
is a $ BV_{{\rm loc}} $  function   determined by the realization of  the
BBM.
We construct it  through   the $h-$extrema  of  BBM where  $ h = \frac
{2\FF^*} {V(\b,\th)} $.
In  Section 5, we recall  the construction done
by Neveu and Pitman,  [{\rcite{NP}],    together with all its 
relevant properties.
Here we only recall what is needed to state the main theorems.   Denote, as in [{\rcite{NP}] , by   
$\{ S_i \equiv S_i^{(h)};  \in \Z \}$   the points of $h-$ extrema with the labeling convention 
that $\dots S_{-1}<S_0\le 0< S_1<S_2\dots$.  
They  proved that $\{ S_i \equiv S_i^{(h)};  \in \Z \}$ is   a stationary renewal process, and
gave  the  Laplace
transform of the inter-arrival times. 
The  $u^*=u^*(W)$ is the  following random $BV_{\rm loc}$ function:
 $$
u^* (r) =
\left \{ \eqalign  { & m_\b, \,\,\,\, \quad \hbox {for} \quad
r \in  [ S_{i}, S_{i+1}),   \quad \hbox {if }\quad 
S_{i} 
 \quad \hbox {is    a point of $ h$--minimum for
$W$ };   \cr & 
 Tm_\b, \quad \hbox {for} \quad  r  \in [ S_{i+1}, S_{i+2}).   }  \right. \Eq (SS.1a)$$
 $$  u^* (r) = \left \{ \eqalign  { & Tm_\b, \quad \hbox {for} \quad  r
\in  [ S_{i}, S_{i+1}),   \quad \hbox {if }\quad
 S_{i} 
 \quad \hbox {is    a point of $h$--maximum for
$W$ };   \cr & 
 m_\b,\,\,\,\, \quad  \hbox {for} \quad  r  \in [ S_{i+1}, S_{i+2}).   } \right. \Eq (SS.1)$$
 For   $W $  and $u^*(W)$ chosen as described,  we denote  for $u \in BV_{{\rm loc}} $
 $$  \G (u|u^*, W)= \sum_{j \in \Z}   \G_{[S_j, S_{j+1})}  (u|u^*(W), W). \Eq
(D.5b)
$$ 
Since    $ |[S_j, S_{j+1})|$ is $\PP$ a.s a finite interval, 
$\G_{[S_j, S_{j+1})}  (u|u^*(W), W)$ is well defined. Actually it  can
be proven,
see Theorem \eqv (min) stated below,  that   the sum  is positive and therefore the functional
in \eqv (D.5b)  is well defined   although it may  be infinite.
The $  \G (\cdot |u^*, W)$  provides an   extension  of   the
functional   \eqv (D.5)
in  $\R$.   One can  formally  write the functional \eqv (D.5b) as 
$$  \G (u|u^*, W)=  \frac 1 {  2 \tilde m_{\b}} \left  \{    \frac { \FF^*} 2
  \int_\R dr     \big[   \|D u
(r)\|_1
-     \|D u^*(r)\|_1  
\big ]     -   V (\b,\th) \int_\R 
    ( \tilde u(r)-  \tilde u^*(r))
d W (r) \right \}  .
  \Eq (D.5a)$$
  The stochastic integral
in \eqv(D.5a) should   be defined but since we use merely \eqv(D.5b) that
is well defined, we leave  \eqv(D.5a) as a formal definition.  
    We have the following result:
\medskip
\noindent{\bf \Theorem (min)} {\it    $\PP$ a.s. one  can 
construct    an unique $u^*= 
u^*(W)
\in  BV_{{\rm loc}} $  such that for any  $u \in BV_{{\rm loc}}   $, 
$  \G (u|u^*, W) \ge 0 $.  } \vskip0.5cm
\medskip
  \noindent  Theorem \eqv (2)  and Corollary \eqv (2bis) stated next    relate  $ u^*_\g
(\o)$  defined in \eqv(L.1)  with  
$u^* (W)$ the minimizer of $  \G (u|u^*, W)$.

\medskip
\noindent {\bf \Theorem (2)} {\it  Given $(\b,\th) \in \EE$, see \eqv
(11.2500),  choosing
the parameters as  in Subsection 2.5, setting $h= \frac {2\FF^*} {V(\b,\th)} $,   we have  that 
$$
\lim_{\g \to 0} \e(\g)\a^*_i(\g) \mathrel{\mathop=^{\rm Law}}  S_i^{(h) } \equiv S_i  \Eq(law)
$$
for   $i \in \Z  $.
 The $\{S_i, i \in \Z \} $ is a stationary renewal process on $\R$.  
The $S_{i+1} -  S_{i}$, (and $S_{-i}-S_{-i-1}$) for   $i>1$
are independent, equidistributed, with Laplace transform
$$ \E[e^{-\l (S_{i+1} -  S_{i})}]= [\cosh (h   \sqrt {2 \l})]^{-1}
\, {\rm for}\,\, \l \ge 0  
\Eq(Laplace1)
$$
(mean $ h^2$)  and distribution given by 
$$ \frac d {dx} \left ( \P[S_2-S_1 \le x] \right)  =  \frac{\pi}{2}
\sum_{k=0}^{\infty} {(-1)^k} \frac {(2k+1)} {h^4}
\exp\left[-(2k+1)^2\frac{\pi^2}{8} \frac{x}{h^2}\right]  
\quad \rm {for}\,\,  x> 0.
\Eq(D.1)
$$
 Moreover $S_1$ and $-S_0$ are equidistributed, have Laplace transform
$$
\E[e^{-\l S_1}]=
\frac{1}{h^2 \l}
\left(1-\frac{1}{\cosh (h\sqrt{2\l})}\right)\,{\rm for}\,\, \l\ge 0 
\Eq(M.6)
$$
and  distribution  given by 
$$\frac d {dx} \left ( \P[S_1 \le x] \right)= \frac{4}{\pi}
\sum_{k=0}^{\infty} \frac{(-1)^k}{(2k+1) {h^2}}  
\exp\left[-(2k+1)^2\frac{\pi^2}{8} \frac{x}{h^2}\right]   
\quad \rm {for}\,\,  x> 0.
\Eq(6)
$$
}
The  formula \eqv (Laplace1) was given by Neveu and Pitman in [\rcite{NP}] and
is reported  here for completeness.    
 
 \medskip
\noindent {\bf \Corollary  (2bis)} {\it  
Under the same hypothesis of Theorem \eqv (2),  for the topology induced by the Skorohod metric that makes
$BV_{\rm loc}$ a complete separable space, see \eqv(Sk4) we have   }
$$
\lim_{\g \downarrow 0} u^*_{\g} \mathrel{\mathop=^{\rm Law}}u^*.
\Eq(law2)
$$
The proof of Theorem \eqv (2) and Corollary \eqv (2bis)  are
given in Section 5.

\vskip0.5cm 
\noindent{ \Remark (2a)}
Note that the Laplace transform \eqv(M.6) is the one of 
the limiting distribution of the age or the residual life of
a renewal process whose Laplace transform of inter-arrival
times is given is \eqv(Laplace1).
The explicit expression given in \eqv(6) 
is    the same   found by
H. Kesten, [\rcite{Ke}], and independently by Golosov,   [\rcite{Go}],
for
the limiting distribution of the
point of localization of the Sinai random walk in one dimension
given that this point is positive, [\rcite{S}].
The formula \eqv(D.1) can be easily obtained from \eqv(6) knowing that
\eqv(6) is the limiting distribution of the age of the above  renewal process.

\vskip0.5cm 
\noindent{\Remark (2b)}
An immediate consequence of Theorems \eqv(copv) and \eqv(2) is 
that to construct the limiting (in Law when $\g\downarrow 0$) typical
profile of the Gibbs measure one can proceed in the following way:
Starting on the right  of the origin take a sample of a random variable
with distribution \eqv(6) and put a mark there and call it $S_1$, make the
same on the left of the origin and call the mark $S_0$.
Then the limiting typical profile will be constant on $[S_0,S_1)$.
To determine if it is equal to $m_\b$ or $Tm_\b$, just 
take a sample of a symmetric Bernoulli random variable $\t$ with value
in $\{m_\b,Tm_\b\}$ and take the profile accordingly. 
To continue the construction, use the above renewal
process to determine  marks $S_2$ and $S_{-1}$ then take for the
typical limiting profile in $[S_1, S_2)$ and $[S_{-1},S_0)$ 
the one with $T\t$ defined in \eqv(4.5P) with $T^2$ the identity, then continue. 

\vskip.5truecm


Recall that the  results of   Theorem \eqv (copv)  are obtained  in a
probability subset $\O_1$ depending on $\g$ and $u^*_\g$ is defined   
only on  the  interval $[-Q,+Q]\equiv [-Q(\g), Q(\g)] $
which is  finite for any fixed $\g$, see
\eqv (queue). 
To our aim it is    convenient
to consider  $u^*_\g\in BV([-Q,+Q])$ as an element of $BV_{\rm loc}$
by replacing $u^*_\g$ by $u^{*Q}_\g$ where for  $ u \in BV([-Q,+Q])$, 
$$
u^{Q}(r)=\cases { u(r\wedge Q), &if $r\ge 0$;\cr
u(r\vee (-Q),&if $r<0$.\cr}
\Eq(exten)
$$
In Theorem \eqv (2) we related the  profile $u^*_\g (\o) \in BV
([-Q(\g), Q(\g)])$, (or what is the same $ u^{*Q}_\g$ in $BV_{\rm
loc}$) to $u^*(W) \in BV_{{\rm loc}} $.
Next result which is a weak large deviation principle,
connects the  random functional $  \G ( \cdot |u^*, W)$, \eqv  (D.5b),
with  a functional  obtained when estimating the (random) Gibbs measures of
the neighboorhood $\PP^\r_{\d,\g,\z,[-Q,+Q]}(u)$ for $u$ 
belonging  to a class of perturbations of  $ u^*_\g(\o) $ that has to
be specified.   
Let us denote for $Q\equiv Q(\g)$ and   $ f(Q)$ a positive increasing real function,  
$$
{\cal U }_{Q }(u^*_\g ) = \left \{  u \in     BV_{\rm loc};
\,\,\,u^Q(r)=u^{*Q}_\g(r), \forall\, |r|\ge Q-1,\,\,\,
V_{-Q}^Q (u)\le V_{-Q}^Q (u^*_\g) f(Q)   \right \}. 
\Eq (D.40bisbis)
$$
The last  requirement   in \eqv  (D.40bisbis) imposes that the number of
jumps of $u$ does not grow too fast with respect to the ones of $u^*_\g$. 
Note that  $u$ in  ${\cal U }_{Q(\g)
}( u^*_\g(\o) )$ is   a random function   depending  on $\g$
that is  $u\equiv u(\g)$ and $u\equiv u(\g,\o)$ if one needs to
emphasize the $\o$ dependence.

 \vskip0.5cm 
\noindent{\bf \Theorem (main1)} {\it Given $(\b,\th) \in \EE$,
let $ u^* \in BV_{\rm loc}$  be the  $\PP$ a.s.   minimizer for $\G (\cdot|u^*, W)$
given in \eqv(SS.1a), \eqv(SS.1).
 Setting the parameters as in Subsection 2.5,  taking 
$$
f(Q)=e^{(\frac{1}{8+4a}-b)(\log Q)(\log\log Q)},
\Eq(effebase)
$$
with $0<b< 1/(8+4a)$, 
for $u(\g,\o) \in {\cal U}_Q(u^*_\g(\o)) $ such that
$$
\lim_{\g\downarrow 0} (u(\g),u^*_\g){\mathop=^{\rm Law}}(u,u^*)
\Eq(limitinlaw)
$$
we have 
$$
\lim_{\g \downarrow 0}\left [
- \g \log \mu_{\b,\th,\g}\Big ( \PP^\rho_{\d,\g,\z, [-Q,Q]} (u   (\g)) \Big )
\right ] \mathrel{\mathop=^{\rm Law}} \G (u|u^*,W).  
\Eq (main)$$
 } 
\vskip0.5cm

\noindent  
Let us consider  some examples: 
Suppose  $u_1 \in {\cal U}_Q(u^*_\g(\o))  $ is such that 
for some $L>0$ 
 $$
 u_1   (\o) (r) = v (r) \1_{[-L,L]}(r) +u^*_\g (\o)\1_{[-Q,Q]\setminus
[-L,L]}(r), \Eq (DD.1)
$$ 
where  $v \in  BV_{{\rm loc}} $ is non random function. 
When $L$ is  a fixed  number independent on $\g$
then $(u_1 , u^*_\g)$ converges in Law when $\g\downarrow 0$ to 
$$( v (r) \1_{[-L,L]}(r) +u^* (r,W)\1_{\R \setminus [-L,L]}(r),u^*(r,W))$$
and the  functional in
the r.h.s.  of \eqv(main) is computed on  
$ u(r) = v (r) \1_{[-L,L]}(r) +u^* (r,W)\1_{\R \setminus [-L,L]}(r)$.
When $ L\equiv L (\g)$ in \eqv (DD.1)  goes  to infinity when $ \g \downarrow 0$
then $(u_1,u^*_\g)$ converges in Law to $(v,u^*)$ and 
the functional in the r. h. s. of \eqv(main) is computed 
on the function  $u(r)\equiv v (r)$.
Theorem \eqv(main1) is a consequence of    accurate  estimates, see
Proposition \eqv(31),
where  approximate terms and errors are explicitly computed.
\medskip
 \noindent   { \bf 2.5  Choice of the parameters } 
We regroup here the choice of the parameters that will be used all
along this work. This choice is similar to the one done in
[\rcite{COPV}]. 
First one chooses a function $g$ on $[1,\infty)$ such that 
$g(x)>1, g(x)/x \le 1, \forall x>1$ and 
$\lim_{x\uparrow \infty} x^{-1}g^{38}(x)=0$.
Any increasing function slowly varying at infinity can be modified to satisfy such
constraints.
A possible choice is  $g(x)=1\vee \log x$ or any iterated of it.
For $\d^*$,   which represents the smallest coarse graining
scale, we have two constraints: 
$$
\frac{(\d^*)^2}{\g}
g^{3/2}(\frac{\d^*}{\g}) \le 
 \frac 1{\b \k (\b,\th) e^3 2^{13}},
\Eq(TE.2)
$$
where $\k(\b,\th)$ is the constant in \eqv 
(2.19)
and 
$$
(\frac {2\g}{\d^*})^{1/2} \left(\log \frac{1}
{\g\d^*} +\frac{\log g(\d^*/\g)}{\log\log g(\d^*/\g)}\right)
\le \frac{1}{32}.
\Eq(MTE.2)
$$
A possible choice of $\d^*$ is 
$$\d^*=\g^{\frac{1}{2} +d^*} \quad {\rm for\, some}\, 0<d^*<1/2.
\Eq(deltastar)
$$ 
 The first constraint, \eqv(TE.2),   is needed    to 
represent in a manageable  form  the multibody interaction 
that comes from the block spin transformation, see  Lemma \eqv (62);
the second one, \eqv(MTE.2) is needed to estimate  the Lipschitz norm when
applying a concentration inequality
to  some function of the random potential.}

Taking $g$ slowly varying at infinity, the conditions  \eqv(deltastar)
\eqv(TE.2) and \eqv(MTE.2) are satisfied by taking $\g$ small enough.
For $(\z, \d)$,   the accuracy chosen to  determine how   close is  the
local magnetization  to the equilibrium values,
there exists a $\z_0=\z_0(\b,\th)$ such that 
$$
\frac {1}{[\k(\b,\th)]^{1/3}g^{1/6}(\frac{\d^*}{\g})}
<\z \le \z_0,
\Eq(TE.1)$$
and $\d$ is  taken   as 
$$
\d= \frac 1{5(g(\frac {\d^*}{\g}))^{1/2}}.
\Eq(5.97037)
$$
The fuzziness $\r$ is chosen as
$$
\r=\left(\frac{5}{g(\d^*/\g)}\right)^{1/(2+a)},
\Eq(rho)
$$
where $a$ is an arbitrary positive number. Note that $\d\le \r$ and
$\r/\d\uparrow \infty$, so in Definition \eqv(part.a) we have just a
constraint  of the form $\g\le \g_0(u)$.
Furthermore
$\e$ that appears in
\eqv(L.1) is chosen as 
$$\e=(5/g(\d^*/\g))^4, \Eq (epsilon) $$
$R_2$ that appears in Definition \eqv(part.a) is chosen as 
$$
R_2=c(\b,\th) (g(\d^*/\g))^{7/2}
\Eq(erre2)
$$
for some positive $c(\b,\th)$, and 
$$Q = \exp [(\log  g(\d^*/\g))/\log\log g(\d^*/\g)].
\Eq(queue)
$$
Note that choosing $\d^*$ as in \eqv(deltastar), in Theorem
\eqv(copv) we have called
$g(1/\g)$ what is really $g(\g^{-(1/2)+d^*})$ however since $g$ is
rather arbitrary this is the same. 
Note also that since $g$ is slowly varying at infinity,
as we  already mentioned $\g R_2\downarrow 0$ when $\g\downarrow0$. 
At last, note  that the only constraint on $\z$ is \eqv(TE.1). 
In particular  one can also choose $\z$ in Theorem \eqv(copv) as 
$$
\z=\z(\g)\equiv 2 \frac {1}{[\k(\b,\th)]^{1/3}g^{1/6}(\frac{\d^*}{\g})}
\Eq(zeta)
$$
that goes to zero with $\g$.

\vskip 1truecm

\chap{3 The block spin representation and Basic  Estimates }3
\numsec= 3
\numfor= 1
\numtheo=1
\medskip
As explained in the introduction the first step  is a coarse graining 
procedure. The microscopic spin system is mapped  into a block spin 
system  (macro scale), for which the length of the interaction 
becomes 1.  In this section we state the  results of this 
procedure, see  Lemmas \eqv 
(40) and \eqv (MM).    The actual  computations    are   
straightforward, but  tedious.   Once  this procedure has been accomplished one 
needs to estimate and  represent  in a        form, convenient for 
further computations,  the
multibody interaction, which is a byproduct of the coarse graining 
procedure,  and the main stochastic contribution to the coarse 
grained energy.  The multibody interaction
is represented as a convergent series  applying a well known 
Statistical Mechanics technique, the Cluster expansion, see Lemma 
\eqv (62).  The  main stochastic  term
is represented with the help of Central Limit Theory, in Proposition 
\eqv (PP).   We then 
give some basic estimates  which we will apply in Section 4.

\noindent  With $C_{\d^*}(V)$ as in \eqv(PJG.1),
let $\S_V^{\d^*}$ denote the
sigma--algebra of $\SS$ generated  by $m^{\d^*}_{V}(\s)\equiv
(m^{\d^*}(x,\s)$, $x\in \CC_{\d^*}(V))$, where $ m^{\d^*}(x,\s)=
(m^{\d^*}(+,x,\s),m^{\d^*}(-,x,\s))$, cf. \equ(2.10).
We take  $I= [i^-,i^+)\subseteq\R$ with $i^{\pm} \in \Z$. The
interval $I$ is assumed to be  $\DD_{\d^*}$--measurable and we set
$\partial^+ I\equiv \{x\in \R\colon i^+ \le x < i^+ +1\}$,
$\partial^- I\equiv \{x\in \R\colon i^- -1 \le x < i^-\}$, and
$\partial I=\partial^+I\cup\partial^-I$.
For  $(m^{\d^*}_I,m^{\d^*}_{\partial I})$  in $\MM_{\d^*}(I \cup
\partial  I)$, cf. \eqv(2.14), we set
 $\widetilde
m^{\d^{*}}(x)=({m_{1}^{\d^{*}}(x)+m_{2}^{\d^{*}}(x)})/{2}$,
$$ E(m^{\d^*}_I)\equiv -\frac {\d^*}2 \sum_{(x,y)\in 
\CC_{\d^*}(I)\times
\CC_{\d^*}(I)} J_{\d^*}(x-y) \tilde m^{\d^*}(x)\tilde m^{\d^*}(y),
\Eq(3.4)
$$
 
$$ E (m^{\d^*}_I, m^{\d^*}_{\partial^{\pm}I}) \equiv
-\d^*\sum_{x\in \CC_{\d^*}(I)}\sum_{y\in
\CC_{\d^*}(\partial^{\pm}I)} J_{\d^*}(x-y) \tilde
m^{\d^*}(x)\tilde m^{\d^*}(y), \Eq(3.5P)
$$
where $J_{\d^*}(x)= \d^*J (\d^*x)$. Further denote 
$$
\eqalign{ {\widehat \FF}(m^{\d^*}_I|m^{\d^*}_{\partial^{ }I} )=&
 E(m^{\d^*}_I)+ E(m^{\d^*}_I,m^{\d^*}_{\partial^{ }I})  -\frac
{\th \d^*}2 \sum_{x\in \CC_{\d^*}(I)} (m^{\d^*}_1(x)-m^{\d^*}_2(
x))\cr  &- \d^* \sum_{x\in \CC_{\d^*}(I)} \frac {\g}{\b\d^*} \log
{{\d^*\g^{-1}/2} \choose
{{{\frac{1+m^{\d^*}_1(x)}{2}\d^*\g^{-1}/2}}}}
  {{\d^*\g^{-1}/2} \choose 
{{{\frac{1+m^{\d^*}_2(x)}{2}\d^*\g^{-1}/2}}}},\cr }\Eq(3.11P)
$$

$$
\GG(m^{\d^*}_{I})
\equiv
 \sum_{x \in \CC_{\d^*}(I)}
 \GG_{x,m^{\d^*}(x)}(\l(x))
 \Eq(3.14P)
$$
where  for each $x\in\CC_{\d^*}(I)$, $\GG_{x,m^{\d^*}(x)}(\l(x))$ is
the cumulant generating function:
$$
\GG_{x,m^{\d^*}(x)}\left(\l(x)\right)
\equiv -\frac1\b\log  \E_{x,m^{\d^*}(x)}^{\d^*} (e^{ 2\b\th 
\l(x)\sum_{i\in
D(x)}\s_i}),
\Eq(3.13P)
$$
of the ``canonical" measure on $\{-1,+1\}^{A(x)}$,  defined through
$$
\E_{x,m^{\d^*}(x)  }^{\d^*}(\vf)=
\frac  {\sum_{\s} \vf(\s)\1_{
\{m^{\d^*}(x,\s)=m^{\d^*}(x)\}}} {\sum_{\s}
\1_{ \{m^{\d^*}(x,\s)=m^{\d^*}(x)\}}},
\Eq(3.9P)
$$
the sum being over $\s \in \{-1,+1\}^{A(x)}$.
Finally denote 
$$
V (m^{\d^*}_I) \equiv V_I (m^{\d^*}_I,h)  = -\frac{1}{\b}\log \E_{
m^{\d^*}( I)}[\prod_{\scriptstyle{x\neq y} \atop \scriptstyle {x,y
\in \CC_{\d^*}(I)\times \CC_{\d^*}(I)} } e^{-\b U
(\s_{A(x)},\s_{A(y)})}].
  \Eq(3.140)
$$
where
$$
U (\s_{A(x)},\s_{A(y)})= -  \sum_{i \in A(x),j  \in A(y)}   \g
  \big[ J(\g|i-j|)-J(\d^*|x-y|) \big ]\s_i\s_j.  \Eq (C.1P)
  $$
and 
$$
\E_{ m^{\d^*}_I}[f]\equiv
 \frac { \sum_{\s_{\g^{-1}I}  } \prod_{x_1
\in \CC_{\d^*}(I) } \1_{ \{m^{\d^*}(x_1,\s)=m^{\d^*}(x_1)\}}e^{ 2
\b \th \l(x_1) \sum_{i\in D(x_1)} \s_i } f(\s) }
 { \sum_{\s_{\g^{-1}I}  }  \prod_{x_1 \in
\CC_{\d^*}(I)} \1_{ \{m^{\d^*}(x_1,\s)=m^{\d^*}(x_1)\}}e^{ 2 \b\th
\l(x_1) \sum_{i\in D(x_1)} \s_i }  }. \Eq (3C.1)
$$
 Let   $F^{\d^*}$   be  a
$\S_I^{\d^*}$-measurable bounded function, $m^{\d^*}_{\partial
I} \in \MM_{\d^*}(\partial I)$ and  
$\mu_{\b,\th,\g}\left(F^{\d^*}|\S^{\d^*}_{\partial I}\right)$  the
conditional expectation of $F^{\d^*}$ given the $\s$--algebra
$\S^{\d^*}_{\partial I}$.  We obtain:  

\medskip
\noindent{\bf \Lemma(40)}{\it  $$
\mu_{\b,\th,\g}\left(F^{\d^*}\bigm|\S^{\d^*}_{\partial^{ }
I}\right) (m^{\d^*}_{\partial^{ }I}) = \frac {e^{\pm \frac \b \g 2
\d^*}}  {Z_{\b,\th,\g,I}(m^{\d^*}_{\partial^{ }I})} 
\sum_{m^{\d^*}_{I} \in {\cal M}_{\d^*}(I)}
F^{\d^*}(m^{\d^*}) e^{-\frac \b \g \left\{{ \widehat
\FF}(m^{\d^*}_I|m^{\d^*}_{\partial^{ }I} )+\g \GG(m^{\d^*}_I)+
\g  V (m^{\d^*}_I)\right\}},
\Eq(3.15P)
$$
 where equality    has to be interpreted as an upper bound for
$\pm=+1$ and a lower bound for $\pm=-1$, and 
$$
 Z_{\b, \g,\th,I}(m^{\d^*}_{\partial I}) =
\sum_{m^{\d^*}_{I}\in {\cal M}_{\d^*}(I) } e^{-\frac \b \g \left\{{ 
\widehat
\FF}(m^{\d^*}_I|m^{\d^*}_{\partial^{ }I} )+\g \GG(m^{\d^*}_I)+
\g  V (m^{\d^*}_I)\right\}}.  \Eq(3.8P1)
$$}

\medskip \noindent   
 That is, up to the error  terms ${e^{\pm \frac \b \g 2 \d^*}}$, we
are able to describe the  system in terms of the  block spin
variables giving  a rather explicit form to the deterministic and
  stochastic part. The  explicit derivation of Lemma \eqv  (40)  
is done  in Section 3 of    [\rcite{COPV}].  
Here we only point out that  since 
$$|\1_{\{\g|i-j|\le 1/2\}}-\1_{\{\d^*|x-y|\le 1/2\}}|\le
\1_{\{-\d^* +1/2\le \d^*|x-y|\le \d^* +1/2\}} \Eq (S.9)
$$
one can estimate  
$$ |U (\s_{A(x)},\s_{A(y)})| \le   \g   (\frac {\d^*}
\g)^2 \1_{\{1/2-\d^*\le \d^*|x-y|\le 1/2+\d^*\}}.
 \Eq (H.1) $$
Therefore,  given $m^{\d^*}_I\in \MM_{\d^*}(I)$, we easily obtain from
\eqv(H.1) 
$$
   \left | H(\s_{\g^{-1}I}) + \th \sum_{i\in \g^{-1}I} h_i\s_i -\frac
1{\g} E(m^{\d^*}_I)  \right| = \frac 1 \b
 \left |  \log \Big [ \prod_{x \in \CC_{\d^*}(I)}\prod_{y \in
\CC_{\d^*}(I)}
e^{-\b  U (\s_{A(x)},\s_{A(y)})  } \Big] \right|   
\le |I| \d^*  \g^{- 1},
\Eq (S.2)
$$
for 
$
 \s \in  \{\s\in \g^{-1} I:
m^{\d^*}(x,\s)=m^{\d^*}(x),\,\forall {x \in \CC_{\d^*}(I)}\}$.
The following lemma 
gives an explicit integral form of the deterministic
part of the block spins system. For $ m\in \cal{T}$, $ f_{\b,\th}(m)$ 
defined in \eqv (2.14000),   let us call
$$
\eqalign{ \tilde \FF(m_I|m_{\del I})&= \int_I f_{\b,\th}(m(x))\,dx
+ \frac  14 \int_I\int_I J(x-y)[\tilde m(x)-\tilde m(y)]^2\,dxdy\cr 
&+\frac 12
\int_I\,dx \int_{I^c} J(x-y) [\tilde m(x)-\tilde m(y)]^2\, dy.\cr 
}\Eq(5.Rio2)
$$
 \medskip
\noindent{\bf \Lemma(MM)}{\it  
Set  $m^{\d^*}_{I\cup \del I}\in
\MM_{\d^*}(I\cup \del I)$,  $m(r)=m^{\d^*}(x)$ for $r \in
 [x\d^*,(x+1) \d^*) $ and $x\in \CC_{\d^*}(I\cup \del I)$, then  one has
$$
|\widehat \FF(m^{\d^*}_I|m^{\d^*}_{\partial I} )- \tilde
\FF(m_I|m_{\del I})+ 
\frac{\d^*}{2} \sum_{y \in \CC_{\d^*}(\del
I)} \big[\tilde m^{\d^*}(y)\big]^2 \sum_{x \in \CC_{\d^*}(I)}
J_{\d^*}(x-y)
|\le |I|\frac{\gamma}{\d^*} \log\frac{\d^*}{\gamma}.
\Eq(5.Rio3)
$$
} 
\proof
 
Using Stirling formula, see [\rcite{Ro}],  we get
$$\eqalign{
\Big| \d^* \sum_{x\in \CC_{\d^*}(I)} \frac
1{2\b}\left(\II(m_1^{\d^*})+\II(m_2^{\d^*})\right) - \d^*
\sum_{x\in \CC_{\d^*}(I)} &\frac {\g}{\b\d^*} \log
{{\d^*\g^{-1}/2} \choose
{{{\frac{1+m^{\d^*}_1(x)}{2}\d^*\g^{-1}/2}}}}
  {{\d^*\g^{-1}/2} \choose 
{{{\frac{1+m^{\d^*}_2(x)}{2}\d^*\g^{-1}/2}}}}
  \Big |\cr & \le \frac
{1}{\b}|I| \frac{\g}{\d^*}\log \frac{\d^*}{\g},} \Eq(3.11P1)
$$
where $\II(\cdot)$ is defined after  \eqv (2.14000).
Recalling the definition of $f_{\beta,\th}(m)$,
 cf. \equ(2.14000),  the lemma is proven. \eop

There are two random terms in \eqv(3.15P):
$\GG(m^{\d^*}_I)$,  the main random contribution, and $V(m^{\d^*}_I)$, the random
expectation of the deterministic term \eqv (C.1P).   To treat them we will use the
following classical deviation  inequality for Lipschitz function of Bernoulli random
variables. See [\rcite{LT}] or  [\rcite{COP1}] for a short proof.

\medskip
\noindent{\bf \Lemma(41)}{\it Let $N$ be a positive integer and
$F$ be a real function on $\SS_N=\{-1,+1\}^{N}$ and for all $i\in
\{1,\dots,N\}$ let
$$
\|\partial_i F\|_{\infty}\equiv \sup_{(h,\tilde h): h_j=\tilde
h_j,\forall j\neq i} \frac{\left|F(h)-F(\tilde
h)\right|}{|h_i-\tilde h_i|}. \Eq(4.13PP)
$$
 If $\P$ is the symmetric Bernoulli measure and
$\|\partial(F)\|^2_{\infty} =\sum_{i=1}^N
\|\partial_i(F)\|^2_{\infty}$ then, for all $t>0$
$$
\P\left[ F-\E(F)\geq t\right]\le  e^{-\frac{t^2}{4
\|\partial(F)\|^2_{\infty}}} \Eq(4.14)
$$ and also}
$$
\P\left[ F-\E(F)\le -t\right]\le  e^{-\frac{t^2}{4
\|\partial(F)\|^2_{\infty}}}. \Eq(4.140)
$$
\medskip \noindent 
  When  considering volumes $I$ that are not too large, we use the
following simple fact that follows from \eqv(3.14P) and \eqv(3.13P) 
$$
|\GG(m^{\d^*}_I)|\le 2\th \sup_{\s_I\in\{-1,+1\}^{I/\g}} \sum_{x\in 
\CC_{\d^*}(I)}
\Big|\sum_{i\in D(x)}\s_i\,\Big| \le  2\th \sum_{x\in
\CC_{\d^*}(I)}|D(x)|.
\Eq(NewPP.1)
$$
Lemma \eqv(41) implies the following rough estimate:

\noindent{\bf \Lemma(60002)}({\bf The rough estimate}) {\it For all
$\d^*>\g >0$ and   for all positive integer $p$,
that satisfy
$$ 12(1+p){\d^*}\log\frac 1\g\le 1
\Eq(4.00510)
$$
 there exists $\O_{RE}=\O_{RE}(\g,\d^*, p)\subseteq \O$  with 
$\P[\O_{RE}]\ge 1-\g^2$
such that on $\O_{RE}$ we have:
$$
\sup_{I\subseteq [-\g^{-p},\g^{-p}]}
 \frac{\sum_{x\in \CC_{\d^*}(I)} \left(|D(x)|-\E[|D(x)|]\right)}
{\sqrt{|I|}} \le \frac {\sqrt{3(1+p)}}{\g} \sqrt{\g \log
\frac{1}{\g}} \Eq(4.0051)
$$
and, uniformly with respect to all intervals
 $I\subseteq [-\g^{-p},\g^{-p}]$,
$$
\sup_{m^{\d^*}_I\in \MM_{\d^*}(I)}\g |\GG(m^{\d^*}_I)||\le 2\th
\left( \frac{|I|}{2} \sqrt{\frac{\g}{\d^*}} + \sqrt {3(1+p)} 
\sqrt{|I|\g
\log\frac 1\g} \right)
\le  2\th |I| \sqrt
{\frac{\g}{\d^*}}.
\Eq(4.0052)
$$
} \medskip
\noindent This lemma  is a direct consequence of Lemma
\eqv(41), since  $
 |D(x)|=
 (|D(x)|-\E[|D(x)|]) + \E[|D(x)|]
$
, $|D(x)|=|\sum_{i\in A(x)}h_i|/2$, and $\E[|D(x)|]\le 
\frac{1}{2}\sqrt{{\d^*}/{\g}}$ by
Schwarz inequality.
 When we use the estimate \eqv(4.0052), $V(m^{\d^*}_I)$ is
estimated using \eqv(S.2) and one has 
$$
\sup_{m^{\d^*}_I \in \MM_{\d^*}(I)} \g |V(m^{\d^*}_I)|\le \d^*|I|.
\Eq(NewPP.2)
$$
\vskip .5truecm

\noindent  However when \eqv(4.0052) and \eqv(NewPP.2) give useless 
results,
one can use Lemma \eqv(41) to estimate $V(m^{\d^*}_I)$
and at some point $\|\del_i V(m^{\d^*})\|_\infty$ will be needed. 
In  Theorem
 8.1  in   [\rcite{COPV}],  with the help of the cluster expansion, 
we prove the
following.

\medskip
\goodbreak
\noindent{\bf \Lemma(62)}{\it For any finite interval $I$, let
$$
\|\del_i V_I\|_\infty \equiv
 \sup_{(h,\tilde h): h_j=\tilde h_j,\forall j\neq i}
\frac{\left |  V_I (m^{\d^*}_I, h) -  V_I (m^{\d^*}_I,\tilde h)
\right|} {|h-\tilde h|}. \Eq (H.220)
 $$
Then,  for all  $\b>0$, for all $\d^*>\g> 0$,
such that
$$
\frac {(\d^*)^2 }{\g}\le  \frac{1}{6e^3\b} \Eq(H.2201)
$$
we have
$$
\sup_{I\subseteq\Z }\sup_{i\in I} \|\del_i V_I\|_\infty \le \frac 1
\b \frac {S }   {1-  S   } , \Eq (H.22)
 $$
where  
  $  0<S \le 6 e^3 \b  \frac {(\d^*)^2} {\g} $. }
\vskip0.5cm

\noindent Together with the above estimates for $V_I$, we   need an
explicit expression for $\GG(m_I^{\d^*})$.
Since  $D(x)\subseteq B^{-\l(x)}(x)$, 
$ \GG_{x,m^{\d^*}(x)}\left(\l(x)\right)$, see \eqv(3.13P),  depends  
only on one component of $ m^{\d^*}(x)$, precisely on 
$m^{\d^*}_{\frac{3+\l(x)}{2}}$.
In fact, we have
$$
\GG_{x,m^{\d^*}(x)}\left(\l(x)\right)
= -\frac{1}{\b}\log
\frac
{\sum_{\s\in\{-1,+1\}^{B^{-\l(x)}(x)}}
\1_{\{m^{\d^*}_{\frac{3+\l(x)}{2}}(x,\s)=m^{\d^*}_{\frac{3+\l(x)}{2}}\}}
e^{2\b\th\l(x)\sum_{i\in D(x)}\s_i}}
{\sum_{\s\in\{-1,+1\}^{B^{-\l(x)}(x)}}
\1_{\{m^{\d^*}_{\frac{3+\l(x)}{2}}(x,\s)=m^{\d^*}_{\frac{3+\l(x)}{2}}\}}},
\Eq(3.13PM)
$$
since the sums over the spin configurations in 
$\{-1,+1\}^{B^{\l(x)}(x)}$
-- the ones that depend on $m^{\d^*}_{\frac{3-\l(x)}{2}}$-- 
cancel out between the numerator and denominator in \eqv(3.9P).
The formula  \eqv (3.13PM) is  almost  useless. One can think about 
making an expansion in $\b\th$ as we basically did in [\rcite{COP1}],
Proposition 3.1  
where $\b\th$ was assumed to be as small as needed. 
 Since here we assume $(\b,\th ) \in \cal E$, one has to find 
another small quantity. Looking at the term $\sum_{i\in D(x)}\s_i$
in \eqv(3.13P) and  setting   
$$
p(x)\equiv p(x,\o)= |D(x)|/|B^{\l(x)}(x)| = 2\g
|D(x)|/ \d^*,
\Eq(2.1690)
$$
it is easy to see that  
for  $I\subseteq \R$,
if 
$$
\left(\frac{2\g}{\d^*}\right)^{1/2} \log \frac{|I|}{\d^*}\le \frac 
{1}{32}, 
\Eq(Newone)
$$
we have
$$
\P\left [ \sup_{x\in \CC_{\d^*}(I)} p(x) > (2\g/\d^*)^{\frac{1}{4}}
\right]
\le e^{-\frac{1}{32}\left(\frac{\d^*}{2\g}\right)^{\frac 12}}.
\Eq(2.17003)
$$
Depending on the  values of $m^{\d^*}_{\frac{3+\l(x)}{2}}$,
$\GG_{x,m^{\d^*}(x)}\left(\l(x)\right)$ has a behavior that
corresponds  to the classical Gaussian, Poissonian, or  
Binomial regimes, as explained in [\rcite{COP1}]. 
  It turns out, see     Remark  4.11 of  [\rcite{COPV}], that 
we need accurate  estimates    only for those values of   
$m^{\d^*}_{\frac{3+\l(x)}{2}}$ for which  
$\GG_{x,m^{\d^*}(x)}\left(\l(x)\right)$ is in  the  Gaussian regime.  
In this regime,  applying the Central Limit Theorem,  we obtain  a 
more convenient representation of  
$\GG_{x,m^{\d^*}(x)}\left(\l(x)\right)$ which is the content of next 
proposition.    Let
$g_0(n)$ be a positive  increasing real function with
$\lim_{n\uparrow \infty}g_0(n)=\infty$ such that $g_0(n)/n$ 
is decreasing to $0$ when $n\uparrow \infty$.

\noindent {\bf \Proposition(PP)}{\it
For all $(\b,\th)\in \cal E$, 
there exist $\g_0=\g_0(\b,\th)$ and
$d_0(\b)>0$ such that for $0<\g\le \g_0$,
$\g/\d^*\le d_0(\b)$,
on the set
 $\{\sup_{x\in\CC_{\d^*}(I)} p(x) \le (2\g/{\d^*})^{1/4}\}$, if
$$
|m^{\d^*}_{\frac{3+\l(x)}{2}}(x)|\le
1-
\left(\frac{g_0(\d^*\g^{-1}/2)}{\d^*\g^{-1}/2}\vee \frac{16p(x)\b\th}
{1-\tanh(2\b\th)}\right),
\Eq(4.PM0)
$$
then
$$
\eqalign{
&\GG_{x,m^{\d^*}(x)}\left(\l(x)\right)=
-\frac{1}{\b}
\log
\frac {\Psi_{\l(x)2\b\th,p(x),m^{\d^*}_{\frac{3+\l(x)}{2}}(x)}}
 {\Psi_{0,0,m^{\d^*}_{\frac{3+\l(x)}{2}}(x)}}
\cr
&
-\frac{1}{\b}|D(x)|\left[ \log\cosh(2\b\th)+
\log\left(1+\l(x)m^{\d^*}_{\frac{3+\l(x)}{2}}(x)\tanh(2\b\th)\right)
+
\hat\varphi(
m^{\d^*}_{\frac{3+\l(x)}{2}}(x),2\l(x)\b\th, p(x) )\right],\cr
}\Eq(4.PM)
$$
where
$$
\left|\hat\varphi(
m^{\d^*}_{\frac{3+\l(x)}{2}}(x),2\l(x)\b\th, p(x) )\right|
\le
\left(\frac{2\g}{\d^*}\right)^{1/4}
 \frac{32\b\th(1+\b\th)}
{(1-|m^{\d^*}_{\frac{3+\l(x)}{2}}(x)|)^2
(1-\tanh(2\b\th))}
\Eq(5.PM)
$$
and
$$
\left|\log
\frac {\Psi_{\l(x)2\b\th,p(x),m^{\d^*}_{\frac{3+\l(x)}{2}}(x)}}
 {\Psi_{0,0,m^{\d^*}_{\frac{3+\l(x)}{2}}(x)}}\right|\le
\frac{18}{g_0(\d^*\g^{-1}/2)} + 
\left(\frac{2\g}{\d^*}\right)^{1/4}c(\b\th),
\Eq(6.PM)
$$
with  
$$
c(\b\th)=\frac{\tanh^2(2\b\th)(1+\tanh^2(2\b\th))^2}
{[1-\tanh^2(2\b\th)]^2[1-\tanh(2\b\th)]^6}.
\Eq(3.25PM)
$$
 }
\vskip0.5cm \noindent
The proof of Proposition \eqv (PP) is given in   Proposition  3.5  
of   [\rcite{COPV}]. 
In the following we deal with quotients  of quantities 
(partition functions)
of the type \eqv (3.8P1) with boundary conditions that might be 
different between numerator and denominator. For this reason 
it  is convenient to introduce the following notations. 
 Let $I$ any finite interval.  
We set    $m^{\d^*}_{\del
I}= (m^{\d^*}_{\del^- I},m^{\d^*}_{\del^+ I})$ and,  see   
\eqv(3.8P1), we denote
 $$
Z_{\b,\th,\g,I}\left(
m^{\d^*}_{\del^- I}=m_{s_1}
,m^{\d^*}_{\del^+ I}=m_{s_2}\right)\equiv
Z_I^{m_{s_1},m_{s_2}}
\Eq(5.52P1)
$$
where $(m_{s_1},m_{s_2}) \in \{m_-,0,m_+\}^2$ 
and for  $m_{s_1}=0$, we set 
in \eqv(3.8P1)
$E(m^{\d^*}_I,m^{\d^*}_{\del^- I})=0$ while 
for $m_{s_2}=0$ we set $E(m^{\d^*}_I,m^{\d^*}_{\del^+ I})=0$.
In a similar way, recalling \eqv(3.15P), if $F^{\d^*}$ is
$\S_I^{\d^*}$--measurable we set 
$$
\frac{Z_I^{m_{s_1},m_{s_2}}(F^{\d^*})}{Z_I^{m_{s_1},m_{s_2}}} \equiv
\frac{ \sum_{m^{\d^*}_{I} \in {\cal M}_{\d^*}(I)} F(m^{\d^*}_I)
e^{-\frac \b \g \left\{ {\widehat
\FF}(m^{\d^*}_I|m^{\d^*}_{\partial^{-}I}=m_{s_1}
m^{\d^*}_{\partial^{ +}I}=m_{s_2}) +\g \GG(m^{\d^*}_I)+ \g  V
(m^{\d^*}_I)\right\}}} {Z_I^{m_{s_1},m_{s_2}}}. \Eq(5.52P4)
$$
Further,  let   $m^{\d^*}_\b$ be one of  the points in
$ \left  \{-1, -1+ \frac {4\g}{\d^*},\dots, 1-
\frac {4 \g}{\d^*},1 \right\}^2$ which is closest to $m_\b$.
 Let $m^{\d^*}_{\b,I}$ be the function
which coincides with $m^{\d^*}_\b$ on $I$ and vanishes outside
$I$ and  $\RR^{\d,\z}(\eta)$ for $\eta = \pm 1$ the set of 
configurations which are  close  with accuracy  $(\d, \z)$,  see \eqv 
(equi), to $m_{\b} $  when $\eta=1$ and to   
$Tm_{\b} $ when $\eta=-1$.  By definition, $|m^{\d^*}_\b-m_\b| \le 8 
\g/\d^*$ and
choosing suitable the parameters we  obtain that 
$ m^{\d^*}_\b$ (resp. $T m^{\d^*}_\b$) is in $\RR^{\d,\z}(+1)$, (resp
$\RR^{\d,\z} (-1)$).
According to the results presented in   Section 2, 
 the typical configurations profiles  are  long runs close to  one 
equilibrium value 
 followed by a jump, then again long runs  close to  the other  
equilibrium value an so on. 
 It is therefore comprehensible that the following quantities will 
play an important role.
 $$
 \frac{Z_I^{0, 0 }(\1_{T \RR^{\d,\z}(\eta) } ) }{Z_I^{0, 0} ( \1_{  
\RR^{\d,\z}(\eta) } )} \equiv \frac{Z_I^{0, 0 
}(T\RR^{\d,\z}(\eta))}{Z_I^{0, 0} (   \RR^{\d,\z}(\eta))}
. \Eq(5.PP4)
$$
Since  the two minima of  $f_{\b,\th}$, see \eqv (2.14000), are $m_\b$ 
and $Tm_\b$,  
we have $T \RR^{\d,\z}(\eta)= \RR^{\d,\z}(-\eta)$, and  we write \eqv 
(5.PP4) as 
$$
\frac
{Z^{0,0}_{I}(\RR^{\d,\z}(-\eta))}
{Z^{0,0}_{I}(\RR^{\d,\z}(\eta))} 
\equiv e^{ \b  \D^\eta\GG(m^{\d^*}_{\b,I })}
\frac {Z_{I,0}^{0,0} \left (\RR(-\eta)\right)} 
{ Z_{I,0}^{0,0} \left  ( 
 \RR(\eta) \right)  } 
\Eq(5.61)
$$
where
$$
\D^{\eta} \GG(m^{\d^*}_{\b, I })
\equiv \eta  
\left[\GG(m^{\d^*}_{\b,
  I })-\GG(Tm^{\d^*}_{\b,  I })\right ]  =  - \eta
\sum_{x\in \CC_{\d^*}( I )} X(x) ,
\Eq(4.800)
$$

$$
X(x) =
\GG_{x, m^{\d^*}_\b}  (\l(x))-
 \GG_{x, T m^{\d^*}_{\b} }(\l(x)),   
\Eq(4.10a)
$$ and 
$$   \eqalign { &  \frac {Z_{I,0}^{0,0} \left (\RR(-\eta) \right)} 
{ Z_{I,0}^{0,0} \left  ( 
 \RR(\eta) \right)  } \cr &  
\equiv
\frac {\sum_{m^{\d^*}_{I} \in {\cal M}^{\d^*}(I)}
\1_{\{\RR^{\d,\z }(\eta)\}} e^{-\frac \b \g\left\{{\widehat
\FF}(m^{\d^*}_{I},0)+ \g\D_0^{-\eta}\GG(m^{\d^*}_{I}) + \g  V 
(Tm^{\d^*}_{I})\right\}}
 }
{\sum_{m^{\d^*}_{I} \in {\cal M}_{\d^*}(I)}
\1_{\{\RR^{\d,\z }(\eta)\}} e^{-\frac \b\g\left\{{\widehat
\FF}(m^{\d^*}_{I},0)+ \g\D_0^\eta \GG(m^{\d^*}_{I}) + \g  V 
(m^{\d^*}_{I})\right\}} },  } 
\Eq(4.9PP)
$$
and $$
\D^\eta_0 \GG(m^{\d^*}_{  I })\equiv
\sum_{x\in \CC_{\d^*}(\tilde I_{12})} \D^\eta_0 \GG^h_{x,m^{\d^*}(x)}
\Eq(4.90P)
$$
where 
 recalling \eqv(3.13P),
$$
\D^\eta_0 \GG^h_{x,m^{\d^*}(x)}=
\GG_{x,T^{\frac{1-\eta}2}m^{\d^*}(x)}(\l(x))-
 \GG_{x, T^{\frac{1-\eta}2}m^{\d^*}_{\b}(x)}(\l(x))
\Eq(4.10)
$$ with $T^0$ equal to the identity.  
When  flipping  $h_i$ to  $-h_i$, for all $i$,  then
$\l(x)\rightarrow -\l(x)$, $B^+(x)\rightarrow B^-(x)$ while
$D(x)$ does not change. Therefore,
$$
\frac {Z^{0,0}_{I}(\RR(-\eta))}
{Z^{0,0}_{I}(\RR(\eta))}(h) = \frac
{Z^{0,0}_{I}(\RR(\eta))}
{Z^{0,0}_{I}(\RR(-\eta))}(-h),
\Eq(4.7)
$$ which implies  that $\log
\frac {Z^{0,0}_{I}(\RR(-\eta))}
{Z^{0,0}_{I}(\RR(\eta))}(h)$ is a 
symmetric  random variable, in particular has mean zero. 
Further the  $X(x)$  in \eqv (4.10a) is a 
symmetric  random variable as it can be directly  checked  inspecting \eqv (3.13P). 
Therefore $\log
\frac {Z^{0,0}_{I,0}(\RR(-\eta))}
{Z^{0,0}_{I,0}(\RR(\eta))}(h)$
  is a symmetric random variable 
having  mean zero and it  has been estimated in [\rcite
{COPV}] applying Lemma \eqv(41). In  [\rcite
{COPV}]    this term
was denoted  
$\frac{Z_{-\eta,0,\d,\z}(I)} {Z_{\eta,0,\d,\z}(I)}$.  The estimate is reported in 
  the next Lemma. \vskip0.5cm
\noindent{\bf \Lemma(401) }
 {\it Given $(\b, \th) \in \cal E $,
 there exist
 $\g_0=\g_0(\b,\th)>0$, $d_0=d_0(\b,\th)>0$, and
 $\z_0=\z_0(\b,\th)$ such that for all $0<\g\le \g_0$, for all 
 $\d^*>\g$ with
 $\g/\d^* \le d_0$, for all $0<\z <\z_0$
 that satisfy the following condition
 $$
 \z  \ge \left(5184 (1+c(\b\th))^2 (\frac{\g}{\d^*})^{1/2}\right)
 \vee \left(12\frac{e^3\b}{c(\b,\th)} \frac{(\d^*)^2}{\g}\right)^2
 \Eq(4.12009)
 $$
 where $c(\b\th)$ is given in \eqv(3.25PM) and $c(\b,\th)$ is another
$\b,\th$ dependent constant, then  for all $a> 0$,
 $$ \P\left[{\max_{I\subseteq  \D_Q}}^*
 \max_{\tilde I_{12}\subseteq I}\left|
  \log  \frac {Z_{\tilde I_{12},0}^{0,0} \left (\RR(-\eta) \right)} 
 { Z_{\tilde I_{12},0}^{0,0} \left  ( 
  \RR(\eta) \right)  } \right|
 \ge \b\frac{4a+12\z }{\g}
 \right]
 \le
 \frac{2  Q }{\e }
 \frac{  e^{-\sfrac{u}{\e} }}{1-e^{-\sfrac{u}{\e} }}
  \Eq(5.78)
 $$
 where ${\max_{I\subseteq  \D_Q}}^*$ denote the maximum over the 
 intervals
  $I\subseteq \D_Q$ such that $|I|=\e\g^{-1}$ and
 $
  u\equiv
 \frac{a^2 \b^2}{8 \z c^2(\b,\th)}.
 $
 }
 \vskip0.5cm \noindent 
To apply Lemma \eqv(41) one needs a control of  the Lipschitz  norm of the object to
be estimated.  The  Lipschitz  norm  of  $\log   \frac {Z_{I,0}^{0,0} \left
(\RR(-\eta)
\right)}  { Z_{I,0}^{0,0} \left  ( 
 \RR(\eta) \right)  }$ is  given in Lemma 4.9 of [\rcite
{COPV}].  
The   estimate  in Lemma \eqv (401) holds for interval $I$, $|I|=\e\g^{-1}$.      
To treat intervals longer than $ \frac \e \g$,   see Lemma \eqv (33),   Section 4,  one
needs
  a  non trivial extension of Lemma \eqv (401)      and a convenient choice of the
parameters involved in the estimate.    This was done  in the proof of  Lemma 6.3 of 
[\rcite {COPV}].    The estimate  \eqv  (5.78)   is useful when $ \e$ is small ($\e \to
0$).  When  dealing with intervals of order $ \frac 1 \g$, ($\e=1$) to get an useful
estimate we need to have 
$ u \to \infty$.  The only way to  obtain this, see the choice of $u$ in Lemma \eqv
(401), is to let $ \z (\g) \to 0$ as $ \g \to 0$.  But $ \z$ is the accuracy we choose and
we would like to have $ \z$  satisfying \eqv  (4.12009), small but not
going to zero when 
$
\g
\to 0$.    So the main effort is to show, that eventhough the accuracy  to define the
vicinity of the profiles to $ m_\b$ or $Tm_\b$ is kept  finite,   it is possible to find
with overwhelming Gibbs probability and $\P$ a.s., blocks in which the typical
magnetization profiles are indeed
 at distance less than $ \z_5$   to $ m_\b$ or $Tm_\b$,  with   $ \z_5(\g) \to 0$ as $ \g
\to 0$. 
This allows to replace the $ \z$ in the definition of $u$ with  $
\z_5$. This   is   done in  Theorem 7.4 of  [\rcite{COPV}] and it will
be applied
when proving  Proposition \eqv(31) in  Section 4. 

To treat the  term in \eqv (4.800) we 
apply  Proposition \eqv(PP) 
on the set   $ \{ p(x) \le (2\g/\d^*)^{1/4}\} $ and we obtain a  very 
convenient representation  for 
$X(x)$ 
$$
X(x)=
-\l(x)|D(x)|\left[ \log\frac{1+m^{\d^*}_{\b,2}\tanh(2\b\th)}
{1-m^{\d^*}_{\b,1}\tanh(2\b\th)}+ \Xi_1(x,\b\th,p(x))\right]
-\l(x) \Xi_2(x,\b\th,p(x))
\Eq(2.170)
$$
where $\Xi_1$ and  $\Xi_2$ are easily obtained from \eqv(4.PM).
Furthermore, choosing $g_0(n)=n^{1/4}$ in  Proposition \eqv(PP), it 
follows
 that
$$
\left|\Xi_1(x,\b\th,p(x))\right|\le
 64  \frac {\b\th (1+ \b \th) }  {(1-m_{\b,1})^2 (1-\tanh(2\b\th))}   
(2 \frac \g {\d^*})^{1/4} \Eq(2.1692)
$$
and
$$
\left|\Xi_2(x,\b\th,p(x)\right| 
\le  (2\frac \g {\d^*})^{1/4} \left [ 36+2 c(\b \th) \right ]
\Eq(2.1700)
$$
where $c(\b\th)$ is given in \eqv (3.25PM).
The $X(x)$ are in fact symmetric  random variables
as it follows from \eqv (2.170). 
We  have  that
$$
\eqalign{
\E[ X(x)\1_{\{p(x)\le (2\g/\d^*)^{1/4}\}}]&=0,\cr
\E[X^2(x)\1_{\{p(x)\le (2\g/\d^*)^{1/4}\}} ]&= \frac {\d^*} {\g}
c(\b,\th,  \g /\d^*)}
\Eq (B.20)
$$
where  $ c(\b,\th, \g /\d^*)$  satisfies
$$
\eqalign{
 c(\b,\th, \g /\d^*) \le& \left( V(\b,\th)\right)^2
\left[ 1+ ( \g/\d^*)^{\frac 1 5 }\right]^2 \cr
 c(\b,\th, \g /\d^*) \ge& \left( V(\b,\th)\right)^2
\left[ 1- ( \g/\d^*)^{\frac 1 5 }\right]^2 \cr
}\Eq(2.21)
$$
and 
 $ V(\b,\th) $ is defined in \eqv (C.1a). 
By  the  results in [\rcite{COPV}],   the runs of configurations close to 
$m_\b$ or to $Tm_\b$ are of order  $1$ in Brownian scale  (  $\frac 1 
{\g^2}$  in micro units), so    it is convenient to 
partition $\R$ into blocks of  
length $\e $, in the Brownian scale;
i.e. each block in micro units is of length $\frac \e {\g^2}$  and 
the basic
assumption is that $ \e \equiv \e (\g)$, $\lim_{\g \to 0}  \e (\g) =
0$,
$\frac \e {\g^2} > \frac { \d^*} \g $,  so that
each block of length  $\frac
\e {\g^2}$   contains at least one block
$A(x)$; to avoid rounding problems we assume $\e/\g\d^*\in \N$, and
that the basic initial partition $A(x)\colon x \in
C_{\d^*}(\R)$ is a refinement of the present one.  
We  define for $\a \in \Z$: 
 $$
 \chi^{(\e)}(\a)\equiv \g \sum_{x:\d^*x \in \tilde A_{\e/\g}(\a)}
X(x)\1_{\{p(x)\le (2\g/\d^*)^{1/4}\}}, \Eq(3.1)
  $$
where   $\tilde A_{\e/\g}(\a)= [\a \frac \e
\g,  (\a+1) \frac \e \g)$
and for the
sake of simplicity the $\g$, $\d^*$ dependence is not explicit.
To simplify further, and if no confusion arises, we
shall write simply $\chi (\a)$.  Note that  $\chi(\a)$ is a symmetric 
random variable and  from
\eqv (B.20) 
$$
\eqalign {  \E[\chi(\a)]& =0\cr
\E[\chi^2(\a)]&=\e c(\b,\th, \g /\d^*).}  
\Eq(2.20)
$$
 It was proved in [\rcite{COPV}], Lemma 5.4,  that there exists
$d_0(\b,\th)>0$
such that if $\g/\d^*\le d_0(\b,\th)$ then for all $\l\in \R$ we have
$$
\E\left[e^{\l\chi(\a)}\right] \le e^{\frac {3\l^2}{4} \e V^2(\b,\th)}
\Eq(PP.10)
$$
where $V(\b,\th)$ is  defined in \eqv (C.1a).

\vskip 1truecm
\chap{4 Finite volume estimates}4
\numsec= 4
\numfor= 1
\numtheo=1

\vskip0.5cm 

In this section, we give upper and lower bounds of
the  infinite volume random Gibbs  probability \break
$\mu_{\b,\th,\g} \big ( \PP^\r_{\d,\g,\z,[-Q,Q]}
(u  )\big)$ in term of finite volume quantities, see  Proposition
\eqv (31). This is the fundamental ingredient in the proof  of  Theorem
\eqv(main1).
By assumption  $u  \in
\UU_Q(u^*_\g(\o))$, see  \eqv (D.40bisbis),   where  for $ \o \in \O_1$, the
probability subset in
Theorem \eqv (copv), $u^*_\g (\o) \in  BV ([-Q,+Q]) $ is   the  profile  defined in
\eqv (L.1).
  There is no lost of generality to assume that there exists a positive  integer 
$L$,  $L<Q$ such that
$$
u_\g (r)=u^*_\g(r),\,\,\, \forall\, |r|\ge L.
\Eq(LF)
$$  
   To avoid  the case that  a jump of $u^*_\g(\o)$ occurs at 
$L$ or $-L$,   we require that  
$$ \{-L\}\cup \{ L \}   \notin   \cup_{i= \k^*(-Q)}^{\k^*(Q)}  [\e\a^*_i-2\r, \e\a^*_i+2\r], \Eq
(GG.1) $$ 
where $\k^*(\pm Q)$  are    defined in   \eqv(num) and $\r$ is chosen as in \eqv(rho). 
To see that requirement \eqv(GG.1) is harmless, let
$$
\O_3 \equiv \O_3 (Q)=  \bigcup_{L \in [1,Q] \cap \Z}\left \{ \o: \{-L\}\cup \{ L \}  \in
\cup_{i= k(-Q)}^{k(Q)}  [\e\a^*_i-2\r,
\e\a^*_i+2\r]\right \}.
\Eq (n.1)
$$ 
We have  the following result. 

\noindent{\bf \Lemma (310)} {\it 
There exist   $\g_0(\b,\th)>0$ and  $a>0$  such that 
for $\g\le \g_0=\g_0(\b,\th)$  we have  
$$  \P [\O_3]  \le  \frac Q{(g(\frac  {\d^*}  \g ))^{\sfrac {1\wedge a} {8(2+a)}}}\le 
\frac 1{(g(\frac {\d^*} \g ))^{\sfrac {1\wedge a} {10(2+a)}}}
\,. \Eq (n.3) $$ }

\proof Note that 
$$ \O_3 \subset \bigcup_{L \in [1,Q] \cap \Z } \left \{ \exists i \in \{\k^*(-Q),\dots,\k^*(Q)\}, \quad \e\a_i^*\in
[L-2\r,L+2\r] \cup [-L-2\r,-L+2\r] \right \}. \Eq (n.a1) $$ 
To estimate the probability of the event  \eqv (n.a1),  we use Lemma \eqv(36) where it is proven
that   uniformly with respect to $Q$  and  with 
$\P$--probability  larger  than  $1-(5/g( {\d^*} /\g))^{\frac{a}{8(2+a)}}$,
$\k^*(Q)$ and $\k^*(-Q)$ are bounded by $ K(Q)$  given  in \eqv (M.1).  The other ingredient is the
estimate   of the 
  probability that $\e\a^*_0\,\, {\rm or}\,\, \e\a_1^* \in [-2\r,+2\r]$. This is done in  Theorem 5.1 of
[\rcite{COPV}] (see  formula 5.29, 5.30 and 6.66 of  [\rcite{COPV}]).  Then  
  for some  $c(\b,\th)$, $a>0$, when  $\g\le \g_0(\b,\th)$ we have  the following: 
$$
\eqalign{
&\P\left [\exists i \in \{\k(^*-Q),\dots,\k^*(Q) \}:\e\a_i^*\in [L-2\r,L+2\r]\right]\cr
&\quad \le 
2 c(\b,\th) K(Q)   [g(\d^*/\g)]^{-1/(4(2+a))} +\left(\frac{5}{g(\frac  {\d^*}  \g)}\right)^{\frac{a}{8(2+a)}}
\le  \frac{1}{(g( \frac {\d^*} \g))^{\sfrac{1\wedge a}{8(2+a)}}}
. }\Eq(PP.01)
$$
By subadditivity one gets \eqv(n.3), recalling that $Q=\exp\left[\log
(g(\frac  {\d^*}  \g))/\log\log (g(\frac  {\d^*} \g))\right]$.   \eop

>From now on, we will always  consider  $\o\in \O_1\setminus \O_3$ and 
 since the union  is  over $L \in [1,Q] \cap \Z$
in \eqv(n.1), this probability  set is the same for all $u\in
\UU_{Q}(u^*_\g(\o) )$.
For  $u \in \UU_Q(u^*_\g(\o))$   and  $L$  so  that \eqv (LF) holds, 
denote  
$$r_1=\inf(r: r>-Q, \, \|D u(r)-D u^*_\g(r)\|_1 >0); \qquad
r_{\hbox {last}} =\sup(r: r < Q, \, \|D u(r)-D u^*_\g(r)\|_1 >0) \Eq (jumps1) $$
where $Du$ is defined before \eqv(D.4). The $r_1$  is the first point  
starting from $-Q$ where $u$ differs from $u^*_\g(\o)$  and $
r_{\hbox {last}}$
is the last point smaller than $Q$ where $u$ differs from  $u^*_\g(\o)$.
We denote by   
$$r_i\qquad  i=1,.. N_1, \qquad r_{N_1}\equiv r_{\hbox {last}} \Eq (jumps)$$ 
 the points of jumps of $u$ or $ u^*_\g(\o)$, between $r_1> -Q $
and   $r_{\hbox {last}}<Q$, in  increasing order.
Note that    $r_i$ could be  a point of jump  for both  $u$ and  $
u^*_\g(\o)$ and
$$
N_1\le N_{[-L,+L]}(u)+N_{[-L,+L]}(u^*_\g(\o)).
\Eq(enne1)
$$
We have the following result.
 \vskip0.5cm   

\noindent{\bf \Proposition  (31)} {\it Take  the parameters as in
Subsection 2.5.  Let $\O_1 $ be  the probability subspace of Theorem
\eqv (copv)  and $\O_3$ defined
in \eqv(n.1),  $ \O_4$ defined in  Corollary \eqv (diff5)   and $\O_5$ defined in \eqv (diff3a). 
 On $\O_1\setminus (\O_3 \cup \O_4\cup  \O_5)$, with $ \P [\O_1\setminus (\O_3 \cup \O_4\cup  \O_5)) ] \ge 1 -
3 (g (\frac {\d^*} \g))^{-\sfrac{1\wedge a}{10(2+a)}} $, for some
$a>0$, for all 
$u \in \UU_Q(u^*_\g(\o)) $ such that
$$
N_{[-Q,+Q]}(u)\le N_{[-Q,+Q]}(u^*_\g)  e^{(\frac{1}{8+4a}-b)(\log Q)(\log\log Q)}
\Eq(condu)
$$
for $0<b<1/(8+4a)$,  there exists a $\g_0=\g_0(\b,\th,u)$ such that
for all $0< \g\le \g_0(\b,\th,u)$, we have that }
$$ 
\eqalign { & 
\frac{\g}{\b}  \log \left
[ \mu_{\b,\th,\g} \Big ( \PP^\r_{\d,\g,\z,[-Q,Q]} (u)  \Big)  \right]=  
\cr & 
- \FF^*
\sum_{i=1}^{N_1}  \left[\frac{
\|D  u(r_i)\|_1-\|D  u^*_\g(r_i)\|_1}
{4\tilde m_\b}
\right]
+\sum_{i=1}^{  N_1 }
\frac{\tilde u( r_{i})-\tilde u^*_\g(r_{i})}{2\tilde m_\b}
\left[\sum_{\a:\,\e\a \in [r_i, r_{i+1})} \chi(\a)\right]  \pm  g(\d^*/\g)^{-b}\cr
}\Eq(L.66)
$$
{\it we have an upper bound for $\pm=+1$ and a lower bound 
for $\pm=-1$.}

\vskip0.5cm \noindent 
 Since the proof of Proposition \eqv (31)  is rather long,  we divide
it in several intermediate steps.
It is convenient to state the following  definitions.
\vskip0.5cm \noindent 
{\bf \Definition  (31.a)} $\,$ {\bf Partition associated to a couple $
(u,v)$  of $BV ([a,b]) $. }
{\it  Let $u$ and $v$ be  in $BV([a,b])$.
   We associate  to  $(u, v)$ the partition of $[a,b]$ obtained 
by taking $C(u, v) =C(u) \cup C(v)$ and $B(u,v)= [a, b]\setminus C(u,
v)$. The  $C(u)$ and $C(v)$ are the elements of the partitions  in
Definition \eqv (part.a).
We  set $C(u, v) =\cup_{i=1}^{\bar N_{[a, b]}}
C_i(u, v)$, where     $\bar N_{[a, b]} \equiv  \bar N (u,v, [a, b])$ is
the number of disjoint intervals  in  $C(u, v)$,  
$\max
\{ N_{[a, b]}(u) ,   N_{[a, b]}(v) \} \le  
\bar N_{[a, b]}
\le N_{[a, b]}(u)+ N_{[a, b]}(v)   $.  
 }
 \vskip0.5cm \noindent 
By  definition, for $i\neq j$,   $C_i(u) \cap C_j(u)=\emptyset $ and    $C_i(v) \cap
C_j(v)=\emptyset $, however when $u$ and $v$ have jumps at distance
less than $\r$,   $C_i(u) \cap C_j(v) \neq
\emptyset$    for some $i \neq j$ and in this case one element of
$C(u,v)$ is $C_i(u)\cup C_j(v)$. 
  \vskip0.5cm \noindent
\Remark (ottorho)    The condition that $\r$ and $\d$ are small
enough in such a way 
that the distance between two  successive jumps of $u$ or $v$ is larger
than $8\r+8\d$, see Definition  \eqv(part.a),   implies that the distance
between any two distinct $C_i(u,v)$ is
at least $2\r+2\d$.   This means that in a given $C_i(u, v)$ there is 
at most two jumps, one of $u$ and the other of $v$.   
 \vskip0.5cm \noindent
The   partition in Definition \eqv (31.a) induces a  partition on the
rescaled (macro) interval 
$\frac 1 \g [a,b]=C_\g(u,v)\cup B_\g(u,v)$ 
where $C_\g(u,v)=\cup_{i=1}^{\bar N_{[a,b]}} C_{i,\g}(u,v)$ and 
$ C_{i,\g}(u,v)= \g^{-1} C_i(u,v)$. 
\smallskip
We will use   Definition \eqv (31.a) 
for  the couple $ (u , u^*_\g(\o)) $ for $u \in \UU_Q(u^*_\g(\o))$, 
$[a,b] = [r_1, r_{\hbox {last}}] $, see \eqv (jumps1).
For simplicity we denote  
$$ \bar N (u,
u^*_\g(\o),[r_1, r_{\hbox {last}}]) \equiv \bar N.   $$
Of course,   $ N_1 \ge  \bar N $, see \eqv(enne1). 
 We write  in macroscale
$$ C_\g(u, u_\g^*) = \cup_{i=1}^{\bar N} [a_i,b_i), \quad
[a_i,b_i) \cap [a_j,b_j)= \emptyset  \qquad 1\le i \neq
j \le \bar N.
\Eq (EP.11)$$
  
 \vskip0.5cm \noindent {\Remark (gammau)}  In Proposition
\eqv(31),  $\g_0=\g_0(\b,\th,u)$ depends on $u$ since 
$8\r(\g)+8\d(\g)$ has to be smaller  than the distance
between two successive jumps of $u$.  

 \vskip0.5cm \noindent
Since   the estimates  to prove  Proposition \eqv (31) are done in intervals written in macroscale   
we  make the following convention: 
 $$
\eqalign {
& m(x)= u(\g x), \quad  m^* (x)= u^*_\g(\g x) \quad \hbox{for}  \quad x\in \frac 1\g
[-Q,Q] \cr
& \quad q_1=-\frac Q \g,
\quad q_2= \frac Q \g; \quad v_1=-\frac L \g,
\quad v_2= \frac L \g;
\quad x_i = \frac {r_i} \g, \quad i=1,\dots,    N_1      
\cr
&\quad \PP^\r_{[q_1,q_2]}
(m)\equiv
\PP^\r_{\d,\g,\z,[-Q,Q ]} (u_\g),
} \Eq (Ma.2) $$
where we recall that  $r_i$  are the points where $u$ or $ u^*_\g$ jumps, see \eqv
(jumps).
Furthermore, let us define 
$$ \eta (\ell,  v )= 
\left \{  
\eqalign { &  0  \qquad \hbox { if}\quad  \ell \in  C_\g(v);\cr 
& 1 \qquad \hbox {when $ m(x)$ equal to $m_\b$  for } x \in   B_\g(v);\cr
& -1 \quad \hbox {when $m(x)$ equal to $Tm_\b$ for }  x \in   B_\g(v).\cr
} \right . $$
Note that  $\eta (\ell,  v )$ is associated to the  function
$v$ not to a block-spin configuration. 
 \vskip0.5cm \noindent  
{\bf \Definition  (35.a) }   {\it 
For  $ \d$ and $\z$ positive,  for  two integers $p_1<p_2$
define 
$$
\OO^{\d,\z}_0([p_1,p_2])
\equiv\big\{\eta^{\d,\z}(\ell)=0,\,\forall
\ell \in [p_1,p_2]\big\}
\Eq(zeri)
$$
and for $\bar \eta \in \{-1,+1\}$,
$$
\RR^{\d,\z}(\bar \eta,[p_1,p_2])
\equiv \big\{ \eta^{\d,\z}(\ell)=\bar \eta, \,\forall \ell \in [p_1,p_2]\big\}.
\Eq(bareta)
$$}
 The first  set $\OO^{\d,\z}_0([p_1,p_2])$ contains configurations for which the 
 block spin variable, see \eqv (2.190), is  $ \z$ far from the equilibrium values.
  The second set  $\RR^{\d,\z}(\bar \eta,[p_1,p_2])$  
 contains configurations for which the 
 block spin variable  is  $ \z$ close  to  the equilibrium value $m_\b$ when  $\bar \eta=1$
 or  $Tm_\b$ when  $\bar \eta=-1$.
\vskip0.5cm \noindent 
Using a simple modification of the rather involved proof of Theorem 7.4 in
[\rcite{COPV}] 
one gets the following.

{\noindent \bf \Proposition(P1)}
{\it Take the parameters as in Subsection 2.5. Let $\O_1 \setminus
\O_3$ be
the probability subset with $\O_1$   in Theorem \eqv(copv) and $\O_3$ defined
in \eqv(n.1).   There exist 
$\g_0=\g_0(\b,\th)$ and  $\z_0$ such that for all $\o\in \O_1$, for
all $\bar \eta\in\{-1,+1\}$, for all $\ell_0\in \N$,  for all $\d,\z,\z_5$ with
$1>\d>\d^*>0$, and any $\z_0>\z>\z_5\ge 8\g/\d^*$, for all
$[\bar p_1,\bar p_2]\subset [p_1,p_2]\subset [q_1,q_2]$ with
$\bar p_1-p_1\ge \ell_0,\, p_2-\bar p_2 \ge \ell_0$, we have
$$
\mu_{\b,\th,\g}\left(
\RR^{\d,\z}(\bar \eta,[p_1,p_2])
\cap \OO^{\d,\z_5}_0([\bar p_1, \bar p_2])
\right)\le
e^{-\frac{\b}{\g} \left\{
(\bar p_2-\bar p_1)\big(\frac{\k(\b,\th)}{4}\d\z_5^3-48(1+\th)
\sqrt{\frac{\g}{\d^*}}\big) -2\z e^{-\a(\b,\th,\z_0) 2\ell_0 }
-4\ell_0 \sqrt{\frac{\g}{\d^*}}\right\}
}. 
\Eq(shrink)
$$
Here $\a(\b,\th,\z_0)$ is a strictly positive constant for all
$(\b,\th) \in \EE$, $\k(\b,\th)$ is the same as in \eqv(2.19).
Moreover
$$
\sup_{ [p_1,p_2] \subseteq [-\g^{-p},\g^{-p}]}
\frac{ Z_{[p_1,p_2]}^{0,0}\left(
\RR^{\d,\z}(\bar \eta,[p_1,p_2])
\cap \OO^{\d,\z_5}_0([\bar p_1, \bar p_2])
\right)}
{ Z_{[p_1,p_2]}^{0,0}\left(
\RR^{\d,\z}(\bar \eta,[p_1,p_2])
\right)}
\Eq(laila.10)
$$
satisfies the same estimates as \eqv(shrink). }
 \vskip0.5cm \noindent
\Remark(errorsr)   The terms with $\sqrt{\g/\d^*}$ in the right hand
side of \eqv(shrink) comes from the rough estimates see Lemma
\eqv(60002).  The fact that  $\a(\b,\th,\z_0)>0$ is a consequence of
\eqv(PP.1). The term $\k(\b,\th)\d\z_5^3$  comes from  estimating the contribution  of \eqv(2.19)  for spin
configuration in  $\OO^{\d,\z_5}_0([\bar p_1, \bar p_2])$. 
 \vskip0.5cm \noindent 
We have the following.  
 \smallskip
{\noindent \bf \Lemma(32) (reduction to finite volume) }
{\it Under the same hypothesis of   Proposition \eqv(31) and on the   probability
space $\O_1\setminus  \O_3$,
for  $\z_5$ that satisfies 
$$
\d\z_5^3\ge 384
(1+\z\frac{\g}{\d^*} +\th)\frac{1}{\k(\b,\th)\a(\b,\th,\z_0)}\sqrt{\frac{\g}{\d^*}}\log
\frac{\d^*}{\g},
\Eq(zeta51)
$$
we have 
$$
\eqalign{
&\mu^{\o}_{\b,\th,\g} \Big 
( \PP^\r_{[q_1,q_2]}(m)  \Big)\ge 
e^{-\frac{\b}{\g}(4\z_5+8\d^*)}
\left({1-2K(Q) e^{-\frac{\b}{\g} \frac{1}{g(\d^*/\g)}}}-
2e^{-\frac{\b}{\g}\frac{\k(\b,\th)}{8}\d\z_5^3}\right)\times\cr
&\quad 
\frac{
Z^{0,0}_{[v_1-1,v_2+1]}  
\Big(\PP^\r_{[v_1,v_2]} (m), 
\eta^{\d,\z_5}(v_1-1)=\eta(v_1-1,m^*), 
\eta^{\d,\z_5}(v_2+1)=\eta(v_2+1,m^*)\Big)}
{Z^{0,0}_{[v_1-1,v_2+1]}
\Big(
\PP^\r_{[v_1,v_2]} (m^*), 
\eta^{\d,\z_5}(v_1-1)=\eta(v_1-1,m^*), 
\eta^{\d,\z_5}(v_2+1)=\eta(v_2+1,m^*)\Big)} \cr
}\Eq(PPPP.1)
$$
and
$$
\eqalign{
&\mu^{\o}_{\b,\th,\g} \Big 
( \PP^\r_{[q_1,q_2]}
(m) 
\Big)\le  2
e^{-\frac{\b}{\g} \left\{ L_1 \frac{\k(\b,\th)}{8}\d\z_5^3 \right \} }
+
e^{\frac{\b}{\g}(4\z_5+8\d^*)}  \times   \cr &\quad 
\sum_{\scriptstyle v_1-L_1-1\le n'_0\le  v_1\atop \scriptstyle
v_2\le n'_{\bar N+1} \le v_2+L_1+1 }
\frac{ 
Z^{0,0}_{[n'_0,n'_{\bar N+1}]}
\Big(\PP^\r_{[n'_0,n'_{\bar N+1}]} (m), 
\eta^{\d,\z_5}(n'_0)=\eta(n'_0,m^*), 
\eta^{\d,\z_5}(n'_{\bar N+1})=\eta(n'_{\bar N+1},m^*)\Big)}
{Z^{0,0}_{[n'_0,n'_{\bar N+1}]}
\Big(\PP^\r_{[n'_0,n'_{\bar N+1}]} (m^*), 
\eta^{\d,\z_5}(n'_0)=\eta(n'_0,m^*), 
\eta^{\d,\z_5}(n'_{\bar N+1})=\eta(n'_{\bar N+1},m^*)\Big)},\cr
}\Eq(PP.1a)
$$
where  $ L_1$ satisfies $L_1+\ell_0 \le \r/\g$ and 
$$
\ell_0=\frac{ \log (\d^*/\g)}{\a(\b,\th,\z_0)}.
\Eq(ellezero)
$$ 
}

\noindent \proof
Recalling \eqv(Ma.2) and \eqv(LF), for $\o\in \O_1\setminus \O_3$ see \eqv(n.1),
one has 
 $$\eta(\ell,m)= \eta(\ell,m^*) \neq 0 \quad \hbox {for} \quad 
 \ell \in [v_1-\frac {\r}{\g} ,v_1 +\frac {\r}{\g}] \cup
 [v_2-\frac {\r}{\g} ,v_2 +\frac {\r}{\g}]. \Eq (prop1)$$ 
Therefore the  spin configurations in  $ \PP^\r_{[q_1,q_2]} (m) $
satisfy  
$$
\eqalign{
&\eta^{\d,\z}(\ell)(\s)=\eta(\ell,m)=\eta(\ell,m^*)
=\eta(v_1,m^*)\neq 0
;\,\forall \ell \in [v_1-\frac {\r}{\g} ,v_1 +\frac {\r}{\g}]\cr
&\quad {\rm and}\quad 
\eta^{\d,\z}(\ell)(\s)=\eta(\ell,m)=\eta(\ell,m^*)=\eta(v_2,m^*) \neq 0;
\,\forall 
\ell \in [v_2-\frac {\r}{\g} ,v_2 +\frac {\r}{\g}].
}\Eq(prop2)
$$
We start proving the lower bound \eqv (PPPP.1).
Within  the proof  we present  a  fundamental   procedure, the {\it  cutting} which 
allows us to estimate  the infinite volume Gibbs measure with 
finite volume quantities.
This procedure will be constantly used in the following. We explain 
here in  details, referring to it  when needed. 

\noindent {\bf The lower  bound}

\noindent We just  impose extra constraints
at $v_1-1$ and $v_2+1$, that is
 $$
\mu_{\b,\th,\g} \Big ( \PP^\r_{[q_1,q_2]} (m)  \Big)
\ge 
\mu_{\b,\th,\g} \Big ( \PP^\r_{[q_1,q_2]}
(m),\eta^{\d,\z_5}(v_1-1)=\eta(v_1-1,m^*),
\eta^{\d,\z_5}(v_2+1)=\eta(v_2+1,m^*)
\Big).
\Eq(PP.3)
$$ 
As mentioned in Section 2.1 there is an unique infinite volume
Gibbs measure that can be obtained as the weak limit of finite volume
Gibbs measure with 0 boundary conditions. So
to  estimate  the infinite volume  Gibbs measure in \eqv(PP.3)
we  start considering the Gibbs measure  in a volume
$[-a,a] $ with $a>0$ big enough so that $[q_1,q_2] \subset [-a,a] $. 
We  write
$$
\frac { Z^{0,0}_{[-a,a]} \Big ( \PP^\r_{[q_1,q_2]}
(m),\eta^{\d,\z_5}(v_1-1)=\eta(v_1-1,m^*),
\eta^{\d,\z_5}(v_2+1)=\eta(v_2+1,m^*)
\Big) } { Z^{0,0}_{[-a,a]}}. 
\Eq(whole1)
$$
The goal is 
to    estimate \eqv  (whole1)   uniformly  with respect to $a$.
This will be achieved  by  {\it cutting } at $v_1-1$ and $v_2+1$, which we explain now.  
Divide the interval  $[-a,a]$ in   three pieces
$[-a,v_1-2]$,$[v_1-1,v_2+1]$, and
$[v_2+2,a]$.  Then, associate the interaction between the first and the second
interval to the first interval, and the one between the second and the third to the third 
interval.
Use \eqv(S.2) with $I=[v_1-2,v_1]$ to make the  block spin
transformation there, this will give an error term $\b 2\d^*/\g$. 
Use  $\eta^{\d,\z_5}(v_1-1)\neq 0$ to get that  for all configurations $\s$
$$
\big|E(m^{\d^*}_{[v_1-2,v_1-1)}(\s),
m^{\d^*}_{ [v_1-1,v_1)}(\s'))-
E(m^{\d^*}_{ [v_1-2,v_1-1)}(\s),
T^{\frac{1-\eta(v_1-1,m)}{2}}
m^{\d^*}_{\b, [v_1-1,v_1)} )\big |\le \z_5
\Eq(Yasmina10)
$$
for $\s'$ such that
$\eta^{\d,\z_5}(v_1-1)(\s'_{[v_1-1,v_1)})=\eta^{\d,\z_5}(v_1-1)$
where  $m_\b^{\d^*}$  is  defined after \eqv (5.52P4). 
Therefore one sees  that up to an error $e^{\pm \frac{\b}{\g}( 2\d^*+\z_5)}$
we can replace in \eqv (whole1)
the $\s,\s'$  interaction  between $[-a,v_1-2]$ and $ [v_1-1, a]$
by an interaction between $\s$ and a constant profile
$  T^{\frac{1-\eta(v_1-1,m)}{2}}m^{\d^*}_{\b,[v_1-1,v_1)}$.  Making
similar computations
in the intervals  $[v_2,v_2+1)$, $ [v_2+2, a)$ and recalling
\eqv(5.52P4),   one gets 
$$
\eqalign{
&
\frac{Z^{0,0}_{[-a,a]}\Big( \PP^\r_{[q_1,q_2]} (m),\, 
\eta^{\d,\z_5}(v_1-1)=\eta(v_1-1,m^*),
\eta^{\d,\z_5}(v_2+1)=\eta(v_2+1,m^*)
\Big)}
{Z^{0,0}_{[-a,a]}}  \cr
&=  e^{\pm \frac{\b}{\g} (2\z_5+4\d^*)}
\frac{1}{Z^{0,0}_{[-a,a]}}
Z^{0,m^*}_{[-a,v_1-2]}\left (\PP^\r_{[q_1,v_1-2]}(m^*)\right )
Z^{m^*,0}_{[v_2+2,a]}\left (\PP^\r_{[v_2+2,q_2]}(m^*)\right )
\times \cr   &\quad \times
Z^{0,0}_{[v_1-1,v_2+1]}
\Big(\PP^\r_{[v_1,v_2]} (m),\, 
\eta^{\d,\z_5}(v_1-1)=\eta(v_1-1,m^*), \,
\eta^{\d,\z_5}(v_2+1)=\eta(v_2+1,m^*)\Big)
.\cr
}\Eq (cutting1)
$$
In the first term on the right hand side of \eqv(cutting1), using
\eqv(LF), one has  $\PP^\r_{[q_1,v_1-2]}(m^*)=\PP^\r_{[q_1,v_1-2]}(m)$ and
$\PP^\r_{[ v_2+2, q_2]}(m^*)=\PP^\r_{[ v_2+2, q_2]}(m)$. 
 Furthermore the boundary condition  $Z^{0,m^*}_{[-a,v_1-2]}(\cdot)$ is written
in term of $m^*$, but since on $v_1-1$ we have $\eta(v_1-1,m)=\eta(v_1-1,m^*)$
we could have also written $Z^{0,m}_{[-a,v_1-2]}(\cdot)$.
Similar considerations hold  for the partition function in $[ v_2+2, a]$.
The above procedure  which allows  to  factorize the partition
function
up to some minor error, see  \eqv(cutting1), will be denoted 
{\it cutting at $v_1-1$  and  $v_2+1$}.
\vskip0.5cm \noindent 
\Remark(cut)
To perform a cutting at some point $\ell$
and  to get an error term  $e^{\pm \frac{\b}{\g}( 2\d^*+\z_5)}$, one needs  to
have $\eta^{\d,\z_5}(\ell)\neq 0$ at this point. Trying to cut at a point $\ell$
where $\eta^{\d,\z}(\ell)=0$ gives an error  term $e^{\frac{\b}{\g}( 2\d^*+1)}$
that will definitively ruin all future estimates. 
Trying to cut at a point $\ell$
where $\eta^{\d,\z}(\ell)\neq 0$ gives an error  term
$e^{\frac{\b}{\g}( 2\d^*+\z)}$.  Since we are not imposing that $\z$
goes to zero, this will also ruin all the future estimates.

 \vskip0.5cm \noindent 
Multiplying and dividing   \eqv(cutting1)  by 
$$
Z^{0,0}_{[v_1-1,v_2+1]}
\Big(\PP^\r_{[v_1,v_2]}(m^*),\,
\eta^{\d,\z_5}(v_1-1)=\eta(v_1-1,m^*), \,
\eta^{\d,\z_5}(v_2+1)=\eta(v_2+1,m^*)\Big)
\Eq(P.220)
$$
and then regrouping, one gets
$$
\eqalign{
&
\frac{Z^{0,0}_{[-a,a]}\Big(\PP^\r_{[q_1,q_2]} (m), 
\eta^{\d,\z_5}(v_1-1)=\eta(v_1-1,m^*), 
\eta^{\d,\z_5}(v_2+1)=\eta(v_2+1,m^*)\Big)}
{Z^{0,0}_{[-a,a]}} \ge \cr
&\quad
\ge e^{-\frac{\b}{\g} (2\z_5+4\d^*)} \frac{1}{Z^{0,0}_{[-a,a]}}
  Z^{0,m^*}_{[-a,v_1-2]}\left (\PP^\r_{[q_1,v_1-2]}(m^*) \right )
Z^{m^*,0}_{[v_2+2,a]}\left (\PP^\r_{[v_2+2,q_2]}(m^*) \right )
\times\cr&\quad\quad\times
Z^{0,0}_{[v_1-1,v_2+1]}
\Big(\PP^\r_{[v_1,v_2]}(m^*),
\eta^{\d,\z_5}(v_1-1)=\eta(v_1-1,m^*), 
\eta^{\d,\z_5}(v_2+1)=\eta(v_2+1,m^*)\Big)\times
\cr
&\quad\quad\times 
\frac{
Z^{0,0}_{[v_1-1,v_2+1]}
\left (\PP^\r_{[v_1,v_2]} (m),\, 
\eta^{\d,\z_5}(v_1-1)=\eta(v_1-1,m^*), \,
\eta^{\d,\z_5}(v_2+1)=\eta(v_2+1,m^*)\right )}
{Z^{0,0}_{[v_1-1,v_2+1]}
\Big(\PP^\r_{[v_1,v_2]} (m^*),\,
\eta^{\d,\z_5}(v_1-1)=\eta(v_1-1,m^*), \,
\eta^{\d,\z_5}(v_2+1)=\eta(v_2+1,m^*)\Big)}.\cr
}\Eq(P.211)
$$
 The last term in \eqv (P.211) is the main contribution to  the lower bound  stated in   \eqv (PPPP.1).  
 So to complete the proof we need  to estimate from below the remaining terms in \eqv  (P.211). 
To achieve this we do the conceptual opposite procedure of cutting.
Namely we  {\it glue} the first three
partitions function in the right hand side of \eqv(P.211) at 
$v_1-1$ and $v_2+1$, applying again  \eqv (Yasmina10).
That is
$$
\eqalign{
&\frac{1}{Z^{0,0}_{[-a,a]}}
Z^{0,m^*}_{[-a,v_1-2]} \left (\PP^\r_{[q_1,v_1-2]}(m^*)\right )
Z^{m^*,0}_{[v_2+2,a]}\left (\PP^\r_{[v_2+2,q_2]}(m^*)\right )
\times\cr&\quad\times
Z^{0,0}_{[v_1-1,v_2+1]}
\Big(\PP^\r_{[v_1,v_2]}(m^*),\,
\eta^{\d,\z_5}(v_1-1)=\eta(v_1-1,m^*), \,
\eta^{\d,\z_5}(v_2+1)=\eta(v_2+1,m^*)\Big)\cr
&\ge  e^{-\frac{\b}{\g}(2\z_5+4\d^*)}\times \cr
&\quad \frac{1}{Z^{0,0}_{[-a,a]}} Z^{0,0}_{[-a,a]}\left (
\PP^\r_{[q_1,q_2]}(m^*),\,\eta^{\d,\z_5}(v_1-1)=\eta(v_1-1,m^*),\, 
\eta^{\d,\z_5}(v_2+1)=\eta(v_2+1,m^*)\right ).
}\Eq(P.212)
$$
By definition  \eqv (D.30) and assumption, see \eqv (prop1), when $ \s \in \PP^\r_{[q_1,q_2]}(m^*)$ 
$\eta^{\d,\z}(v_1-1, \s)\neq 0$ and $
\eta^{\d,\z}(v_2+1, \s)\neq 0$,    then 
$$
\eqalign{
&\PP^\r_{[q_1,q_2]}(m^*)=
\left(\PP^\r_{[q_1,q_2]}(m^*),\,\eta^{\d,\z_5}(v_1-1)=0,\eta^{\d,\z_5}(v_2+1)=0\right)\cup
\cr &
\left(\PP^\r_{[q_1,q_2]}(m^*),\,\eta^{\d,\z_5}(v_1-1)=\eta(v_1-1,m^*), 
\eta^{\d,\z_5}(v_2+1)=\eta(v_2+1,m^*)\right).  \cr
}\Eq(orrore2)
$$
Therefore taking the limit when $a\uparrow \infty$ in the right hand
side of \eqv(P.212), using Theorem \eqv(copv) and Proposition \eqv(P1)
with $ \ell_0=\frac{ \log (\d^*/\g)}{\a(\b,\th,\z_0)}$, {\bf $\bar p_2-  \bar p_1 =1$, $  p_2-    p_1 > 2 \ell_0$ }
one gets   
$$
\eqalign{
&\lim_{a\uparrow \infty}
\frac{1}{Z^{0,0}_{[-a,a]}} Z^{0,0}_{[-a,a]}\big(
\PP^\r_{[q_1,q_2]}(m^*),\,\eta^{\d,\z_5}(v_1-1)=\eta(v_1-1,m^*), 
\eta^{\d,\z_5}(v_2+1)=\eta(v_2+1,m^*)\big) \cr
&\quad \geq
1-K(Q)e^{-\frac{\b}{\g}\frac{1}{g(\d^*/\g)}}-
2  e^{-\frac{\b}{\g}\left(\frac{\k(\b,\th)}{4}\d\z_5^3
-48(1+\th)\sqrt{\frac{\g}{\d^*}}
\frac{\log\frac{\d^*}{\g}}{\a(\b,\th,\z_0)}\right)}, 
}\Eq(limite)
$$
where   $ K(Q)$ is  given  in \eqv (M.1).
 Collecting   \eqv (cutting1), \eqv(P.211), and \eqv(P.212), using
\eqv(zeta51)  one gets
$$  \eqalign { & 
\mu^{\o}_{\b,\th,\g} \Big ( \PP^\r_{[q_1,q_2]} (m) \Big)
\ge e^{-\frac{\b}{\g}(4\z_5+8\d^*)} 
\left( 1-2K(Q)e^{-\frac{\b}{\g} \frac{1}{g(\d^*/\g)}}-
2  e^{-\frac{\b}{\g} \frac{\k(\b,\th)}{8}\d\z_5^3}
\right) \times \cr & 
\frac{
Z^{0,0}_{[v_1-1,v_2+1]}  
\Big(\PP^\r_{[v_1,v_2]} (m),\, 
\eta^{\d,\z_5}(v_1-1)=\eta(v_1-1,m^*), \,
\eta^{\d,\z_5}(v_2+1)=\eta(v_2+1,m^*)\Big)}
{Z^{0,0}_{[v_1-1,v_2+1]}
\Big(
\PP^\r_{[v_1,v_2]} (m^*),\, 
\eta^{\d,\z_5}(v_1-1)=\eta(v_1-1,m^*), \,
\eta^{\d,\z_5}(v_2+1)=\eta(v_2+1,m^*)\Big)} } 
\Eq(P.214)
$$
which is  \eqv(PPPP.1).

\noindent {\bf The upper bound}  
In the proof of the lower bound we could  {\it cut}
making an error proportional to $\z_5$ by simply restricting to those
configurations  having magnetization   close to 
 the  equilibrium values  with accuracy $(\d, \z_5)$
in the chosen  $[\ell, \ell+1)$ block .  In the upper bound
obviously this procedure cannot be applied.  We need to
find a  block     where the   spin configurations  have
magnetization   close to the  equilibrium values  with accuracy $(\d, \z_5)$.  
This makes notations more cumbersome. To facilitate the reading,  we
use  indexes with a $'$
to denote  the  points  $\ell$  where $\eta^{\d,\z_5} (\ell) \neq 0$.
We search these points  within the intervals $[v_1 -L_1-1, v_1-1]$ and  $[ v_2-1, v_2 +L_1+1 ]$
 where $ L_1$ is an integer which will be suitable chosen.
>From  \eqv(prop2)  we have that  
 $ \eta^{\d,\z}(\ell)(\s)=\eta(\ell,m^*)
=\eta(v_1,m^*)\neq 0 $  for $ \ell \in [v_1 -L_1-1, v_1-1]$ and
$ \eta^{\d,\z}(\ell)(\s)=\eta(\ell,m^*)
=\eta(v_2,m^*)\neq 0 $  $ \ell \in  [ v_2-1, v_2 +L_1+1 ]$,  provided $ L_1< \frac \r \g $. 
 Then we apply     Proposition \eqv(P1) in  both the intervals, taking
$\ell_0$ as in \eqv(ellezero) and setting   
$$
\eqalign{
&
p_1=v_1-L_1-\ell_0, \,\, p_2=v_1 +\ell_0, \cr
&\bar p_1=v_1 -L_1-1,\qquad \bar p_2 = v_1
}\Eq(shrink2)
$$
for some $L_1$ such that $0<L_1 +\ell_0 \le \r/\g$ to be chosen later. 
We have $\bar p_2-\bar p_1= L_1+1$
and
$$
\eqalign{
\mu_{\b,\th,\g}\left(
\RR^{\d,\z}(\bar \eta,[p_1,p_2])
\cap \OO^{\d,\z_5}_0([\bar p_1, \bar p_2])
\right)
&\le
e^{-\frac{\b}{\g} \left\{
(L_1+1) \big(\frac{\k(\b,\th)}{4}\d\z_5^3-48(1+\z\frac{\g}{\d^*} +\th)
\sqrt{\frac{\g}{\d^*}}\big) -
4 \frac{ \log (\d^*/\g)}{\a(\b,\th,\z_0)}\sqrt{\frac{\g}{\d^*}} \right\}}\cr
&\le 
e^{-\frac{\b}{\g} \left\{
L_1 \frac{\k(\b,\th)}{8}\d\z_5^3  \right\}}.
}\Eq(shrink3)
$$
At the last step, we have used  \eqv(zeta51). We do  
similarly  in the interval   $ [v_2 , v_2 +L_1+1] $ . 
We apply  Proposition \eqv(P1)  setting 
$$
\eqalign{
&
p_3=v_2-\ell_0, \,\, p_4=v_2 +L_1+\ell_0 \cr
&\bar p_3=v_2,\,\,
\bar p_4= v_2+L_1+1,
}\Eq(shrink4)
$$
then one gets that
$ \mu_{\b,\th,\g}\big(
\RR^{\d,\z}(\bar \eta,[p_3,p_4])
\cap \OO^{\d,\z_5}_0([\bar p_3, \bar p_4])
\big)
$
satisfies the same estimate as in \eqv(shrink3).
Therefore one has the basic estimate
$$
\mu_{\b,\th,\g}\big(\PP^\r_{[q_1,q_2]}(m)\big)
\le
\mu_{\b,\th,\g}\big(\PP^\r_{[q_1,q_2]}(m)
\cap \big(\OO^{\d,\z_5}_0([\bar p_1, \bar p_2])\big)^c
\cap \big(\OO^{\d,\z_5}_0([\bar p_3,\bar p_4])\big)^c\big)
+ 2
e^{-\frac{\b}{\g} \left\{ L_1 \frac{\k(\b,\th)}{8}\d\z_5^3\right\}}. 
\Eq(funda)
$$
In the set  $\big(\PP^\r_{[q_1,q_2]}(m)
\cap \big(\OO^{\d,\z_5}_0([\bar p_1,\bar p_2])\big)^c
\cap \big(\OO^{\d,\z_5}_0([\bar p_3, \bar p_4])\big)^c\big)$ there exists 
at least one block variable indexed by  $n'_0$ with    $\bar p_1\le n'_0\le \bar p_2$ such that
$\eta^{\d,\z_5}(n'_0)=\eta(n'_0,m^*)$ and   one block variable 
indexed by  $n'_{\bar N+1}$,     $\bar p_3\le n'_{\bar N+1} \le \bar p_4$
where $\eta^{\d,\z_5}(n'_{\bar
N+1})=\eta(n'_{\bar N+1},m^*)$. 
These   are the blocks where we  will cut. 
 Consider the Gibbs measure in a volume
$[-a,a] $ with $a>0$ large  enough so that $[q_1,q_2] \subset [-a,a] $. 
 We have the simple estimate 
$$
\eqalign{
&\frac{Z^{0,0}_{[-a,a]}\Big(\PP^\r_{[q_1,q_2]} 
\cap \big(\OO^{\d,\z_5}_0([\bar p_1, \bar p_2])\big)^c
\cap \big(\OO^{\d,\z_5}_0([\bar p_3,\bar p_4])\big)^c \Big)}
{Z^{0,0}_{[-a,a]}}\le\cr
&\quad
\sum_{\scriptstyle  \bar p_1\le n'_0\le  \bar p_2\atop \scriptstyle \bar p_3\le
n'_{\bar N+1}\le
\bar p_4}
\frac{Z^{0,0}_{[-a,a]}\Big(\PP^\r_{[q_1,q_2]}
(m),\,\eta^{\d,\z_5}(n'_0)=\eta(n'_0,m),\,
\eta^{\d,\z_5}(n'_{\bar N+1})=\eta(n'_{\bar N+1},m)\Big)}
{Z^{0,0}_{[-a,a]}}.
}
\Eq(fund2)
$$
Consider now a generic  term in the sum in the right hand side of
\eqv(fund2).
   Recalling \eqv(cutting1) and 
 cutting at $n'_0$ and $n'_{\bar N+1}$, in the numerator we get 
$$
\eqalign{
& \frac{Z^{0,0}_{[-a,a]}\Big(\PP^\r_{[q_1,q_2]} (m),\,\eta^{\d,\z_5}(n'_0)=\eta(n'_0,m),\,
\eta^{\d,\z_5}(n'_{\bar N+1})=\eta(n'_{\bar N+1},m)\Big)}{Z^{0,0}_{[-a,+a]}}  \le e^{\frac{\b}{\g}
(2\z_5+4\d^*)}
\frac{1}{Z^{0,0}_{[-a,a]}}\cr
&\times Z^{0,m^*}_{[-a,n'_0-1]}(\PP^\r_{[q_1,n'_0-1]}(m^*))
\times
\cr &\times Z^{0,0}_{[n'_0,n'_{\bar N+1}]}
\Big(\PP^\r_{[n'_0,n'_{\bar N+1}]} (m),\, 
\eta^{\d,\z_5}(n'_0)=\eta(n'_0,m^*), 
\eta^{\d,\z_5}(n'_{\bar N+1})=\eta(n'_{\bar N+1},m^*)\Big)
Z^{m^*,0}_{[n'_{\bar N+1}+1,a]}(\PP^\r_{[n'_{\bar N+1}+1,q_2]}),\cr
}\Eq(P.21)
$$
see \eqv (5.52P4)  to recall notations. Multiplying and dividing by
$$
Z^{0,0}_{[n'_0,n'_{\bar N+1}]}
\Big(\PP^\r_{[n'_0,n'_{\bar N+1}]}   (m^*)\, , 
\eta^{\d,\z_5}(n'_0)=\eta(n'_0,m^*), 
\eta^{\d,\z_5}(n'_{\bar N+1})=\eta(n'_{\bar N+1},m^*)\Big)
$$ 
one gets, after regrouping the terms
$$
\eqalign{
&
\frac{Z^{0,0}_{[-a,a]}\Big(\PP^\r_{[q_1,q_2]} (m),\,  
\eta^{\d,\z_5}(n'_0)=\eta(n'_0,m), 
\eta^{\d,\z_5}(n'_{\bar N+1})=\eta(n'_{\bar N+1},m)\Big)}
{Z^{0,0}_{[-a,a]}}\le  e^{\frac{\b}{\g} (2\z_5+4\d^*)}\times 
\cr
& 
\frac{ 
Z^{0,0}_{[n'_0,n'_{\bar N+1}]}\Big(
\PP^\r_{[n'_0,n'_{\bar N+1}]} (m),\,  
\eta^{\d,\z_5}(n'_0)=\eta(n'_0,m^*), 
\eta^{\d,\z_5}(n'_{\bar N+1})=\eta(n'_{\bar N+1},m^*)\Big)}
{
Z^{0,0}_{[n'_0,n'_{\bar N+1}]}
\Big(\PP^\r_{[n'_0,n'_{\bar N+1}]} (m^*),\,  
\eta^{\d,\z_5}(n'_0)=\eta(n'_0,m^*), 
\eta^{\d,\z_5}(n'_{\bar N+1})=\eta(n'_{\bar N+1},m^*)\Big)}\times \cr
&\frac{
Z^{0,m^*}_{[-a,n'_0-1]}(\PP^\r_{[q_1,n'_0-1]}(m^*))}
{Z^{0,0}_{[-a,a]}}\times \cr
&\quad  Z^{0,0}_{[n'_0,n'_{\bar N+1}]}
\Big(\PP^\r_{[n'_0,n'_{\bar N+1}]} (m^*),\,  
\eta^{\d,\z_5}(n'_0)=\eta(n'_0,m^*), 
\eta^{\d,\z_5}(n'_{\bar N+1})=\eta(n'_{\bar N+1},m^*)\Big)
Z^{m^*,0}_{[n'_{\bar N+1}+1,a]} (\PP^\r_{[n'_{\bar N+1}+1,q_2]})
}.\Eq(P.22)
$$
Now, glueing 
at $n'_0$ and $n'_{\bar N+1}$, as in \eqv(P.212), 
uniformly with respect  to $a$, we have 
$$
\eqalign{
&\frac{ 
Z^{0,m^*}_{[-a,n'_0-1]}(\PP^\r_{[q_1,n'_0-1]}(m^*))}
{Z^{0,0}_{[-a,a]}}\times\cr
&Z^{0,0}_{[n'_0,n'_{\bar N+1}]}
\Big(\PP^\r_{[n'_0,n'_{\bar N+1}]} (m^*),\,  
\eta^{\d,\z_5}(n'_0)=\eta(n'_0,m^*), 
\eta^{\d,\z_5}(n'_{\bar N+1})=\eta(n'_{\bar N+1},m^*)\Big)
\times \cr
&\quad Z^{m^*,0}_{[n'_{\bar N+1}+1,a]}(\PP^\r_{[n'_{\bar N+1}+1,q_2]})
\le e^{\frac{\b}{\g}(2\z_5+4\d^*)} 
\frac{Z^{0,0}_{[-a,+a]}\Big(\PP^\r_{[q_1,q_2]} (m^*),\,  \Big)}
{Z^{0,0}_{[-a,a]}}\le  e^{\frac{\b}{\g}(2\z_5+4\d^*)}. }
\Eq(P.23)
$$
 From \eqv (P.22) and \eqv (P.23) we get
$$
\eqalign{
&\frac{Z^{0,0}_{[-a,a]}\Big(\PP^\r_{[q_1,q_2]} (m),\,  
\eta^{\d,\z_5}(n'_0)=\eta(n'_0,m), \,
\eta^{\d,\z_5}(n'_{\bar N+1})=\eta(n'_{\bar N+1},m)\Big)}
{Z^{0,0}_{[-a,a]}}\le  e^{\frac{\b}{\g} (4\z_5+8\d^*)}\times \cr
&\times 
\frac{ 
Z^{0,0}_{[n'_0,n'_{\bar N+1}]}
\Big(\PP^\r_{[n'_0,n'_{\bar N+1}]} (m),\,   
\eta^{\d,\z_5}(n'_0)=\eta(n'_0,m^*), \,
\eta^{\d,\z_5}(n'_{\bar N+1})=\eta(n'_{\bar N+1},m^*)\Big)}
{Z^{0,0}_{[n'_0,n'_{\bar N+1}]}
\Big(\PP^\r_{[n'_0,n'_{\bar N+1}]} (m^*),\,   
\eta^{\d,\z_5}(n'_0)=\eta(n'_0,m^*), \,
\eta^{\d,\z_5}(n'_{\bar N+1})=\eta(n'_{\bar N+1},m^*)\Big)}.\cr
}\Eq(P.24)
$$
Collecting \eqv(fund2) and \eqv(P.24), one get  \eqv(PP.1a).
This ends the proof of the lemma. \eop

\smallskip 
\noindent

The   configurations in  $\PP^\r_{[v_1, v_2]} (m)$ and  $\PP^\r_{[v_1, v_2]} (m^*)$
are long runs of  $\eta^{\d,\z} (\ell) \neq 0$ followed by phase  changes  
in the intervals $[a_i, b_i)$, for $i=1,\dots \bar N$, see \eqv (EP.11).  So
to estimate the ratio of the partition function
in \eqv (PPPP.1) and \eqv(PP.1a),  it is convenient to   separate the contribution
given by those intervals in which the spin configurations  undergo  to
a phase change, {\it i.e}  in which  the block spin variables  are  $\eta^{\d,\z}
(\ell)=0$,  from those
intervals in which the   block spin variables  are  $\eta^{\d,\z}
(\ell)\neq 0$.
This can be achieved {\it cutting}  at suitable points
the partition function. We require these points to be such that    $\eta^{\d,\z_5} (\ell) \neq 0$
to obtain  error terms which are  negligible. 
 We  start   proving an   upper bound for \eqv(PP.1a).  To facilitate  the
reading, as before,
we use  indexes with $'$  and $''$  to denote  the  points  $\ell$  where $\eta^{\d,\z_5} (\ell) \neq 0$.   
 Denote  a generic term  in    \eqv(PP.1a) by 
 $$
{\cal Z} (n'_0,n'_{\bar N+1})
\equiv \frac
{Z^{0,0}_{[n'_0,n'_{\bar N+1}]}
\Big(\PP^\r_{[n'_0,n'_{\bar N+1}]} (m),\,   
\eta^{\d,\z_5}(n'_0)=\eta(n'_0,m^*), \,
\eta^{\d,\z_5}(n'_{\bar N+1})=\eta(n'_{\bar N+1},m^*)\Big)}
{Z^{0,0}_{[n'_0,n'_{\bar N+1}]}
\Big(\PP^\r_{[n'_0,n'_{\bar N+1}]} (m^*),\,   
\eta^{\d,\z_5}(n'_0)=\eta(n'_0,m^*), \,
\eta^{\d,\z_5}(n'_{\bar N+1})=\eta(n'_{\bar N+1},m^*)\Big)}.
\Eq(up)
$$
 We have the following:
 \smallskip
{\noindent \bf \Lemma(32D)  }
{\it Under the same hypothesis of   Proposition \eqv(31), on the   probability
space $\O_1\setminus   ( \O_3 \cup \O_4)$, and 
for  $\z_5$  as in \eqv (zeta51), 
we have 
 $$
\eqalign { & {\cal Z} (n'_0,n'_{\bar N+1})  \le 
e^{-\frac{\b}{\g} \frac{\FF^*}{2\tilde m_\b}
\sum_{-L\le r \le L} \left[|D\tilde u(r)|-|D\tilde u^*_\g(r)|\right]}
e^{ \frac \b \g 
\sum_{i=1}^{\bar N} \frac{\tilde u(r_i)-\tilde u^*_\g(r_i)}{2\tilde m_\b}
\left[\sum_{\a:\,\e\a\in [r_i, r_{i+1})}\chi(\a)\right]}\cr
&\quad \times  e^{ \frac{\b}{\g}\bar N\left [ 4\z_5+8\d^*+ \g \log\frac{\r}{\g}+\g\log L_1+
\frac{20 V(\b,\th)}{
(g(\d^*/\g))^{1/4(2+a)}}+32\th(R_2+\ell_0+L_1)\sqrt{\frac{\g}{\d^*}}\right]  }\cr &
+ \bar N ^2 e^{\bar N \log \frac{\r}{\g}}   e^{\frac{\b}{\g}(8\d^*+4\z)}
e^{-\frac{\b}{\g}L_1\frac{\k(\b,\th)}{8}\d\z^3_5}. }   
\Eq(up1)
$$ }
\proof 
Recalling\eqv(D.30), and \eqv(EP.11), one sees that in each interval $[a_i,b_i]$,
there is a single phase change   on a length $R_2$ for $m$ or
$m^*$. 
There are three possible cases:

\noindent {\bf Case 1}  $[a_i,b_i] \in C_\g(u)$ and $[a_i,b_i] \in B_\g(u^*_\g)$.
Therefore 
$$
\eqalign{
\eta (a_i, m)&=-\eta(b_i,m)\neq 0\cr
\eta (a_i, m^*)&=\eta(b_i,m^*)\neq 0.\cr
}\Eq(PP.26)
$$
 {\bf Case 2} $[a_i,b_i] \in B_\g(u)$ and $[a_i,b_i] \in C_\g(u^*_\g)$.
Therefore
$$
\eqalign{
\eta (a_i, m)&=\eta(b_i,m)\neq 0\cr
\eta (a_i, m^*)&=-\eta(b_i,m^*)\neq 0. \cr
}\Eq(PP.27)
$$
 {\bf Case 3} $[a_i,b_i] \in C_\g(u)$ and $[a_i,b_i] \in C_\g(u^*_\g)$.
Therefore
$$
\eqalign{
\eta (a_i, m)&=-\eta(b_i,m)\neq 0\cr
\eta (a_i, m^*)&=-\eta(b_i,m^*)\neq 0. \cr
}\Eq(PP.27e)
$$
 In the first two cases   there exists an unique $x_i\in   [a_i,b_i]$, see \eqv (Ma.2), 
so that,  in the  the case {1}, $|D\tilde m(x_i))| >0$ and 
in the case 2,   $|D\tilde m^*(x_i)|>0$.
In the  case 3,  both $m$ and $m^*$ have  one jump in $[a_i,b_i]$.
Recalling \eqv(bareta)   and Definition \eqv(1) we denote 
$$
\WW_1(\ell_i,m)\equiv
\WW_1([\ell_i-R_2,\ell_i+R_2],R_2,\z) \cap
\{\eta^{\d,\z}(\ell_i-R_2)=\eta(a_i,m),\eta^{\d,\z}(\ell_i+R_2)=\eta(b_i,m)\},
\Eq(wave)
$$
the set of configurations undergoing to a phase change  induced by $m$ in $ [\ell_i-R_2,\ell_i+R_2]$.  We 
denote  in the cases   {1}  and {3} 
$$\PP^\r_{   [a_i,b_i]} (m,\ell_i,i)\equiv  \RR^{\d,\z}(\eta(a_i,m),[a_i,\ell_i-R_2-1]) \cap  
\WW_1(\ell_1,m) \cap
\RR^{\d,\z}(\eta(b_i,m),[\ell_i+R_2+1,b_i])      \Eq (OPP.2)  $$ 
and in the case {2}
 $$  \eqalign { & \PP^\r_{   [a_i,b_i]} (m,\ell_i,i)\equiv
\cr &   \RR^{\d,\z}(\eta(a_i,m),[a_i,\ell_i-R_2-1])\cap
\RR^{\d,\z}(\eta(a_i,m),[\ell_i-R_2,\ell_i+R_2])
 \cap
\RR^{\d,\z}(\eta(b_i,m),[\ell_i+R_2+1,b_i]).  }       \Eq (OPP.3)  $$ 
 The set   $\PP^\r_{   [a_i,b_i]} (m,\ell_i,i)$   denotes the  spin
configurations  which,
in the case {1} and {3}, have  a  jump in the
 interval  $ [a_i,b_i] $, starting after  the point  $\ell_i-R_2$ and
ending
before $\ell_i+R_2$ and  close  to  different  equilibrium values 
in $[a_i,b_i] \setminus [ \ell_i-R_2, \ell_i+R_2]$.
In the case {2}, it  denotes the  spin configurations which are in
all
$ [a_i,b_i] $ close to one equilibrium value, namely they do not have jumps.
The $ \ell_i$ in this last case is written  for future use. 
We    use  for both $m$ and $m^*$ the   notation \eqv  (OPP.2)  and \eqv  (OPP.3). 
In the case  {3} both $m$ and $ m^*$ have a jump in $ [a_i,b_i]$.  Obviously  we have 
 $$
\PP^\r_{   [a_i,b_i]} (m)\subset \bigcup_{\ell_i \in [a_i+R_2+1,
b_i-R_2-1]} \PP^\r_{   [a_i,b_i]} (m,\ell_i,i).  \Eq (OPP.1)  $$ 
To get an upper bound for \eqv(up), we  use the subadditivity of  the numerator
in \eqv(up) to treat the
$\cup$ in \eqv(OPP.1) obtaining  a sum over   $\ell_i\in
[a_i+R_2+1,b_i-R_2-1]$.  For the  denominator  we  obtain an upper bound  simply restricting to  the  subset of 
 configurations which is suitable for us, namely 
$$
\PP^\r_{   [a_i,b_i]} (m^*)\supset 
\PP^\r_{   [a_i,b_i]} (m^*,\ell_i,i). 
\Eq(OPP.121)
$$
To short notation,  let $\underline \ell\subset [\underline
a,\underline b] \equiv \{\ell_i\in
[a_i+R_2+1,b_i-R_2-1],\, \forall i,  1\le i\le \bar N\}$ and 
set  
$$ \AA(m, \underline \ell) \equiv  \Big(\PP^\r_{[n'_0,n'_{\bar N+1}]}
(m)
\cap_{i=1}^{\bar N} \PP^\r_{[a_i,b_i]} (m,\ell_i,i),\,  
\eta^{\d,\z_5}(n'_0)=\eta(n'_0,m^*), \,
\eta^{\d,\z_5}(n'_{\bar N+1})=\eta(n'_{\bar N+1},m^*)\Big), \Eq  (short) $$
   $$
 {\cal Z} (n'_0,n'_{\bar N+1},\underline \ell)
\equiv   \frac
{Z^{0,0}_{[n'_0,n'_{\bar N+1}]}
\Big(\AA(m, \underline \ell)  \Big)}
{Z^{0,0}_{[n'_0,n'_{\bar N+1}]}
\Big(\AA(m^*, \underline \ell)\Big)}.
 \Eq(up202)
$$
Therefore, recalling \eqv(up),  we can write
$$
{\cal Z} (n'_0,n'_{\bar N+1})\le \sum_{\underline \ell\subset [\underline
a,\underline b]} {\cal Z} (n'_0,n'_{\bar N+1},\underline \ell).
\Eq(up201)
$$
 The 
number of terms in the sum in \eqv(up201) does not exceed $\prod_{i=i}^{\bar N}
(b_i-a_i)\le \exp(\bar N\log (\r/\g))$.
For future use, when  $B$ is an event let us define
$$
 {\cal Z} (n'_0,n'_{\bar N+1},\underline \ell;B)
\equiv  
  \frac
{Z^{0,0}_{[n'_0,n'_{\bar N+1}]}
\Big(\AA(m, \underline \ell) \cap B  \Big)}
{Z^{0,0}_{[n'_0,n'_{\bar N+1}]}
\Big(\AA(m^*, \underline \ell)\Big)}.
 \Eq(up2021)
$$
For $\ell_0$ defined
in \eqv(ellezero), for the very same $L_1$ to be chosen later and
$\z_5$ that satisfies \eqv(zeta51) , let us denote $\bar R_2=R_2+\ell_0$
and  define
$$
\eqalign{
&\DD(m,\underline\ell)\equiv\cr
&\quad
\cup_{1\le i\le \bar N}\left( {\cal R}^{\d,\z}( 
\eta(\ell_i-R_2,m),[\ell_i-R_2-L_1-2\ell_0,\ell_i-R_2])
\cap
{\cal O}^{\d,\z_5}([\ell_i-\bar R_2-L_1,\ell_i-\bar R_2])
\right)\cup\cr
&\quad
\cup_{1\le i\le \bar N}\left({\cal R}^{\d,\z}( 
\eta(\ell_i+R_2,m),[\ell_i+R_2,\ell_i+R_2+2\ell_0+L_1])
\cap
{\cal O}^{\d,\z_5}([\ell_i+\bar R_2,\ell_i+ \bar R_2+L_1])\right)
.\cr
}\Eq(up203)
$$
 The $ \DD(m,\underline\ell) $ is the set of configurations which
are simultaneously $\z$ close and $\z_5$ distant,  (recall $\z>\z_5$), 
from the equilibrium values in the interval
$[\ell_i-R_2-L_1-2\ell_0,\ell_i-R_2] \cup  [\ell_i+\bar R_2,\ell_i+ \bar
R_2+L_1]$ where $\ell_i$ are chosen as in \eqv (OPP.121).
Recalling \eqv(up2021) and  Proposition  \eqv(P1) one gets 
$$
\sum_{\underline \ell \subset [\underline a,\underline b]}
{\cal Z} (n'_0,n'_{\bar N+1},\underline \ell;\DD(m,\underline\ell))
\le   \bar N ^2 e^{\bar N \log \frac{\r}{\g}}   e^{\frac{\b}{\g}(8\d^*+4\z)}
e^{-\frac{\b}{\g}L_1\frac{\k(\b,\th)}{8}\d\z^3_5}.
\Eq(errup)
$$
To get \eqv (errup) one cuts at the points  $\ell_i+R_2$ and  $\ell_i+R_2+2\ell_0+L_1$ for the set 
${\cal R}^{\d,\z}( 
\eta(\ell_i+R_2,m),[\ell_i+R_2,\ell_i+R_2+2\ell_0+L_1])\cap
{\cal O}^{\d,\z_5}([\ell_i+\bar R_2,\ell_i+ \bar R_2+L_1]) $, 
and at the  points $\ell_i-R_2-L_1-2\ell_0$ and  $\ell_i-R_2$ for    the 
set  $ {\cal R}^{\d,\z}( 
\eta(\ell_i-R_2,m),[\ell_i-R_2-L_1-2\ell_0,\ell_i-R_2])
\cap
{\cal O}^{\d,\z_5}([\ell_i-\bar R_2-L_1,\ell_i-\bar R_2])$. 
 Notice that    we cut at points  $\eta^{\d,\z}\ne 0$
and  each time we make the error $ e^{\frac{\b}{\g}(2\z+4\d^*)}$.  This 
is the only place where making an error so large does not cause a
problem.
Namely we can choose $L_1$ suitable in \eqv(errup) so that
$L_1 \frac{\k(\b,\th)}{16}\d\z_5^3 > (8\d^*+4\z)$. Furthermore  denote 
$$
\BB(\underline\ell)
\equiv
\cap_{1\le i\le \bar N} 
\left({\cal O}^{\d,\z_5}([\ell_i-\bar R_2-L_1,\ell_i-\bar R_2])\right)^c
\cap 
\left({\cal O}^{\d,\z_5}([\ell_i+\bar R_2,\ell_i+ \bar R_2+L_1])\right)^c. 
\Eq(buoni)
$$
Since  for each $ \underline \ell$, $\AA(m, \underline \ell) \cap
\DD(m,\underline\ell)^c   \subset  \AA(m, \underline \ell) \cap  \BB(\underline\ell)$
we are left  to estimate
$$
\sum_{\underline \ell \subset [\underline a,\underline b]}
{\cal Z} \left(n'_0,n'_{\bar N+1},\underline \ell;\BB(\underline\ell)\right).
$$
 On each $ \AA(m, \underline \ell) \cap \left({\cal
O}^{\d,\z_5}([\ell_i-\bar R_2-L_1,\ell_i-\bar R_2])\right)^c$,  $ 1 \le i  \le \bar N$,  
there exists at least one block, say $[n'_i, n'_i+1)$  contained in 
$[ \ell_i-\bar R_2-L_1,\ell_i-\bar R_2) $ with
$\eta^{\d,\z_5}(n'_i)=\eta(a_i,m)$.  
Making the same on the right of $\ell_i$ and indexing  $n''_i$ the
corresponding block where $\eta^{\d,\z_5}(n''_i)=\eta(b_i,m) $, one gets
$$
\eqalign{
&\sum_{\underline \ell \subset [\underline a,\underline b]}
{\cal Z} \left(n'_0,n'_{\bar N+1},\underline \ell,\BB(\underline\ell)\right)\le\cr
&\quad 
\sum_{\underline \ell \subset [\underline a,\underline b]}\,\,
\sum_{\scriptstyle \underline n'\subset [\underline \ell-\bar R_2-L_1,\underline
\ell-\bar R_2]\atop\scriptstyle
\underline n''\subset [\underline \ell+\bar R_2,\underline
\ell+\bar R_2+L_1]}
{\cal Z} \left(n'_0,n'_{\bar N+1},\underline \ell;\,\cap_{1\le i\le \bar N}
\{ \eta^{\d,\z_5}(n'_i)=\eta(a_i,m),\eta^{\d,\z_5}(n''_i)=\eta(b_i,m)
\}\right).\cr
}\Eq(buoni2)
$$
The number of terms in the second sum of  \eqv(buoni2) does not exceed
$\exp(2\bar N(\log L_1))$.
Consider now a generic term in \eqv(buoni2),
$$
{\cal Z} \left(n'_0,n'_{\bar N+1},\underline \ell;
\cap_{1\le i\le \bar N}
\{ \eta^{\d,\z_5}(n'_i)=\eta(a_i,m),\eta^{\d,\z_5}(n''_i)=\eta(b_i,m)
\}\right).
\Eq(gene)
$$
Recalling \eqv(up2021), 
we cut the numerator of the partition function, as in \eqv (P.21),
at   the points $\underline n'$ and $\underline
n''$   to get an upper bound. 
Each time    we  cut we get the  error term
$e^{\frac{\b}{\g}(2\d^*+\z_5)}$.   
In the denominator, see \eqv(up202),  restrict the configurations to
be in 
$$
 \AA(m^*, \underline \ell)  \cap_{1\le i\le \bar N}
\left \{ \eta^{\d,\z_5}(n'_i)=\eta(a_i,m^*),\eta^{\d,\z_5}(n''_i)=\eta(b_i,m^*)\right \}
\Eq(buoni12)
$$ 
and then cut at all the points $\underline n'$ and $\underline
n''$.   In this way we obtain  an upper bound for \eqv(gene). 
 We use  notation   \eqv(OPP.2)  (case    {1} and   {3}) and \eqv(OPP.3) (case   {2})  
after cutting at $n'_i$ and $n''_i$.  Note that  $\eta (n'_i+1)=
\eta(a_i,m)$ and
$ \eta (n''_i-1)=  \eta(b_i,m)$ therefore we have 
in the  case  {1} and  {3},  see  \eqv(OPP.2),    
$$
\eqalign{
&\PP^{\r}_{[n'_i+1,n''_i-1]}(m,\ell_i,i)= \cr
&\quad \RR^{\d,\z}(\eta(a_i,m),[n'_i+1,\ell_i-R_2-1]) \cap  
\WW_1(\ell_i,m) \cap
\RR^{\d,\z}(\eta(b_i,m),[\ell_i+R_2+1,n''_i-1]),
}\Eq(OPP.1bis)
$$
in the case 2, see \eqv(OPP.3), 
$$
\eqalign{
&\PP^{\r}_{[n'_i+1,n''_i-1]}(m,\ell_i,i)=
\RR^{\d,\z}(\eta(a_i,m),[n'_i+1,\ell_i-R_2-1]) \cap\cr  
&\quad \RR^{\d,\z}(\eta(a_i,m),[\ell_i-R_2,\ell_i+R_2])
\cap
\RR^{\d,\z}(\eta(b_i,m),[\ell_i+R_2+1,n''_i-1]).
}\Eq(OPP.120bis)
$$
For the remaining parts   corresponding  to runs    between  two phase changes, {\it i.e} the intervals
$[n''_i,n'_{i+1}]$,      $ n''_i \in [a_i,b_i]$ and  $ n'_{i+1}\in  [a_{i+1},b_{i+1}]$, for $  i
\in \{ 1, \dots,  \bar N \} $, we denote 
$$
\PP^\r_{[n''_i,n'_{i+1}]} (m,\z_5)\equiv
\RR^{\d,\z}(\eta(b_i,m),[n''_i+1,n'_{i+1}-1])
\cap \{ \eta^{\d,\z_5}(n''_i)= \eta^{\d,\z_5}(n'_{i+1})=\eta(b_{i},m)\}.
\Eq(lbp)
$$
Similarly   in  the intervals $[n'_0,n'_1]$, and $[n''_{\bar N},n'_{\bar N+1}]$,
recalling \eqv(prop1) and \eqv(prop2), we 
have 
$$\eqalign{
&
\PP^\r_{[n'_0,n'_1]}(m,\z_5)\equiv \RR^{\d,\z}(\eta(v_1,m),[n'_0,n'_1])
\cap
\{\eta^{\d,\z_5}(n'_0)= \eta^{\d,\z_5}(n'_1)=\eta(v_1,m)\} \cr & 
=
\PP^\r_{[n'_0,n'_1]}(m^*,\z_5), }   
\Eq(xzero)$$
$$
\eqalign{
&\PP^\r_{[n''_{\bar N},n'_{\bar N+1}]}(m,\z_5 ) \equiv  \RR^{\d,\z}(\eta(v_2,m),[n''_{\bar N},n'_{\bar N+1}])
\cap \{\eta^{\d,\z_5}(n''_{\bar N})= \eta^{\d,\z_5}(n'_{\bar N+1})=\eta(v_2,m)\}  \cr & = 
\PP^\r_{[n''_{\bar N},n'_{\bar N+1}]}(m^*,\z_5).  
}\Eq(xbarN)
$$
As a result we have 
$$
\eqalign{
&{\cal Z} \left(n'_0,n'_{\bar N+1},\underline \ell;
\cap_{1\le i\le \bar N}
\{ \eta^{\d,\z_5}(n'_i)=\eta(a_i,m),\eta^{\d,\z_5}(n''_i)=\eta(b_i,m)
\}\right)\le \cr
&\quad   
 e^{+ \bar N  \frac {\b}{ \g} (4\z_5+8\d^*) }
\frac 
{Z_{[n'_0, n'_1]}^{0,0} \left (\PP^\r_{[n'_0,n'_1]}(m,\z_5) \right)} 
{ Z_{[n'_0 , n'_1]}^{0,0} 
\left ( \PP^\r_{[n'_0,n'_1]} (m^*,\z_5)\right)  }
\times
 \cr 
& \quad
\prod_{i=1}^{ \bar N-1} 
\left ( 
\frac 
{Z_{[n'_i+1,n''_i-1]}^{m,m} 
\left( \PP^\r_{  [n'_i+1,n''_i-1]} (m,\ell_i,i) \right ) }
{ Z_{[n'_i+1,n''_i-1]}^{m^*,m^*}
\left  ( \PP^\r_{ [n'_i+1,n''_i-1]} (m^*,\ell_i,i) \right )  }
 \frac 
{ Z_{[n''_i, n'_{i+1}]}^{0 ,0} \left (\PP^\r_{ [n''_i, n'_{i+1}]} (m,\z_5) \right)}
{ Z_{[n''_i, n'_{i+1}]}^{0 ,0} \left  ( 
 \PP^\r_{  [n''_i, n'_{i+1}]} (m^*,\z_5 )  \right )  }
\right ) \times
\cr & \quad 
\frac {  Z_{[n'_{\bar N}+1,n''_{\bar N}-1]}^{m,m} 
\left  ( 
 \PP^\r_{ [n'_{\bar N}+1,n''_{\bar N }-1]} (m,\ell_{\bar N},\bar N)  
\right )  }
{Z_{[n'_{\bar N}+1,n''_{\bar N}-1]}^{m^*,m^*}
\left ( 
 \PP^\r_{ [n'_{\bar N}+1,n''_{\bar N }-1]} (m^*,\ell_{\bar N},\bar N)  
\right )  }
\frac 
{ Z_{[n''_{\bar N}, n'_{\bar N+1}]}^{0 ,0} \left( 
 \PP^\r_{ [n''_{\bar N },n'_{\bar N+1}]} (m,\z_5)\right)}
{Z_{[n''_{\bar N},n'_{\bar N+1} ]}^{0 ,0} \left ( 
 \PP^\r_{ [n''_{\bar N }, n'_{\bar N+1}]} (m^*,\z_5)\right )} 
 .}\Eq(gene1)
$$
\vskip .5truecm
\noindent
\Remark(bc)
 Note that the boundary conditions of restricted partition functions
as $$Z_{[n'_i+1,n''_i-1]}^{m,m} ( \PP^\r_{  [n'_i+1,n''_i-1]}
(m,\ell_i,i))$$ in \eqv(gene1) are related on the left to
$\eta(a_i,m)$ and on the right to $\eta(b_i,m)$, see \eqv(OPP.1bis)
and \eqv(OPP.120bis).  

\vskip .5truecm
\noindent
Now the goal is to estimate separately all the ratios in the right hand
side of \eqv(gene1).
It follows from \eqv(prop1), \eqv(xzero), and \eqv(xbarN) that 
$$
\frac 
{Z_{[n'_0, n'_1]}^{0,0} \left (\PP^\r_{[n'_0,n'_1]}(m,\z_5) \right)} 
{ Z_{[n'_0 , n'_1]}^{0,0} 
\left ( \PP^\r_{[n'_0,n'_1]} (m^*,\z_5)\right)  }
=\frac 
{ Z_{[n''_{\bar N}, n'_{\bar N+1}]}^{0 ,0} \left(  
 \PP^\r_{ [n''_{\bar N },n'_{\bar N+1}]} (m,\z_5)\right)}
{Z_{[n''_{\bar N},n'_{\bar N+1} ]}^{0 ,0} \left ( 
 \PP^\r_{ [n''_{\bar N }, n'_{\bar N+1}]} (m^*,\z_5)\right )} 
=1.
$$
The remaining  ratios are estimated 
in Lemma \eqv (33a), Corollary \eqv(diff5) and 
Lemma \eqv (phases)   given below. 
\vskip .3truecm 
\noindent{ \bf Collecting}
We  insert the results of  Lemma
\eqv(33a),   Corollary  \eqv (diff5) and  Lemma
\eqv(phases) in \eqv(gene1). 
To  write in a unifying way  the contributions of   the jumps  
we note that for \eqv(tipo21bis) 
$$
-\FF^*= -\frac{\FF^*}{2\tilde m_\b}\sum_{a_i\le s\le b_i} |D\tilde m ( s)|=
-\frac{\FF^*}{2\tilde m_\b}
\sum_{a_i\le s\le b_i} \left (  |D\tilde m ( s)|-|D\tilde m^*(s)|\right) 
\Eq(PP.45)
$$
since  in the case 1, see \eqv(PP.26), $\sum_{a_i\le s\le b_i} |D\tilde m^*(s)|=0$. 
For \eqv(tipo22)
$$
+\FF^*=\frac{\FF^*}{2\tilde m_\b}\sum_{a_i\le s\le b_i} |D\tilde m^*(s)|=
-\frac{\FF^*}{2\tilde m_\b}\sum_{ a_i\le s\le b_i}  \left (   |D\tilde m( s)|-|D\tilde m^*(s)|\right) 
\Eq(PP.46)
$$
since  in the case 2, see \eqv(PP.27),  $\sum_{a_i\le s\le b_i} |D\tilde m ( s)|=0$.
Moreover, since neither $\tilde m$ nor $\tilde m^*$ have jump in $[b_i+1,a_{i+1}]$
for $i\in \{1,\dots,\bar N\}$, in $[v_1,a_1-1]$,  and in $[b_{\bar N}+1,v_2]$,
one gets simply 
$$
\prod_{i=1}^{ \bar N} 
e^{-\frac{\b}{\g}
\frac{\FF^*}{2\tilde m_\b}\sum_{a_i\le s\le b_i} 
\left[|D\tilde m( s)|-|D\tilde m^*( s)|\right]}
=
e^{-\frac{\b}{\g} \frac{\FF^*}{2\tilde m_\b}
\sum_{-L\le r \le L} \left[|D\tilde u(r)|-|D\tilde u^*_\g(r)|\right]}.
\Eq(PP.47)
$$
 Using  \eqv(PP.24bis), the random terms give a contribution 
$$
e^{ \frac \b \g 
\sum_{i=1}^{\bar N} \frac{\tilde u(r_i)-\tilde u^*_\g(r_i)}{2\tilde m_\b}
\left[\sum_{\a:\,\e\a\in [r_i, r_{i+1})}\chi(\a)\right]}.
\Eq(randompart)
$$
It remains to collect the error terms, see \eqv (errup), \eqv(buoni2), \eqv(gene1),
\eqv(PP.24bis), and Lemma \eqv(phases).
Denote    
$$
\EE_1 \equiv \bar N \left [ 4\z_5+8\d^*+ \g \log\frac{\r}{\g}+\g\log L_1+
\frac{20 V(\b,\th)}{ (g(\d^*/\g))^{1/4(2+a)}}+32\th(R_2+\ell_0+L_1)\sqrt{\frac{\g}{\d^*}}\right ],  
\Eq(errors)
$$
$$
-\AA_2\equiv
2\frac \g \b \log  \bar N+ \frac \g \b \bar N \log \frac{\r}{\g}  + 8\d^*+4\z 
- L_1\frac{\k(\b,\th)}{8}\d\z^3_5, 
\Eq(A2)
$$
and 
$$
 \AA \equiv   \frac{\FF^*}{2\tilde m_\b}
\sum_{-L\le r \le L} \left[|D\tilde u(r)|-|D\tilde u^*_\g(r)|\right] 
-
\sum_{i=1}^{\bar N} \frac{\tilde u(r_i)-\tilde u^*_\g(r_i)}{2\tilde m_\b}
\left[\sum_{\a:\,\e\a\in [r_i, r_{i+1})}\chi(\a)\right].   \Eq(penu2)
$$
We  have  proved  
$$ {\cal Z} (n'_0,n'_{\bar N+1} ) \le
e^{-\frac{\b}{\g}   \AA }  e^{ \frac{\b}{\g}  \EE_1  } +e^{-\frac \b \g \AA_2} 
\Eq(antepenu1)
$$
that entails  \eqv (up1). \eop 

\noindent   Next we state  the lemmas used for estimating the  different ratios in \eqv (gene1).  

\vskip0.5cm   
\noindent{\bf \Lemma (33a)}    {\it Under the same hypothesis of   Proposition \eqv(31) and on the 
probability space $\O_1\setminus \O_3$, for all $1\le i\le \bar N-1$,
for all $n''_i,n'_i$, we have 
$$ 
\frac 
{ Z_{[n''_i, n'_{i+1}]}^{0 ,0}
\left (\PP^\r_{ [n''_i, n'_{i+1}]} (m,\z_5) \right)}
{ Z_{[n''_i, n'_{i+1}]}^{0 ,0}
\left (\PP^\r_{  [n''_i, n'_{i+1}]} (m^*,\z_5 ) \right)  }  
= \left \{ \eqalign { & 1 \qquad \hbox {when } \eta(b_i,m)=\eta(b_i,m^*);  \cr & 
 \eqalign { & e^{\pm \frac {\b}{\g} \frac{ 1}{4c^2(\b,\th) g(\d^*/\g)}}
e^{ \frac \b \g 
\frac{\tilde m(b_i)-\tilde m^*(b_i)}{2\tilde m_\b}
\left[\sum_{\a:\,\e\a\in\g[n''_i+1, n'_{i+1}-1]}\chi(\a)\right]}  \cr
& \qquad  {\hbox {when } }\quad  \eta(b_i,m)=-\eta(b_i,m^*),  }  } \right.  
\Eq(PP.24)
$$
where in the last term  we have an upper bound for $\pm=+$ and a lower bound for $\pm=-$. }
\vskip0.3cm 
\noindent 
\proof
 When   $\eta(b_i,m)=\eta(b_i,m^*)$  \eqv (PP.24) is immediate, see definition  \eqv(lbp).   
When $\eta(b_i,m)=-\eta(b_i,m^*)$  the estimate is a direct
consequence of Lemma \eqv (33). \eop
\vskip0.5cm \noindent 
\noindent
 The  r.h.s of \eqv(PP.24)     gives when $ \eta(b_i,m)=-\eta(b_i,m^*)$ a term which should give   
the second sum in the right hand side of \eqv(L.66).
However in \eqv(PP.24) one has $\sum_{\a:\,\e\a\in\g[n''_i+1,
n'_{i+1}-1]}\chi(\a)$ instead of 
$\sum_{\a:\,\e\a \in [r_i, r_{i+1})}
\chi(\a)$ in \eqv(L.66).   To obtain the term in \eqv(L.66) we  add   the 
missed random field  to   reconstruct  $\sum_{\a:\,\e\a \in [r_i, r_{i+1})}
\chi(\a)$.  We  therefore  need to subtract the same term we added
as a result one has
\vskip0.5cm \noindent 
\noindent  {\bf \Corollary(diff5)} { \it Under the same hypothesis of   Proposition \eqv(31) and on the 
probability space $\O_1\setminus (\O_3\cup \O_4)$,
with $\P(\O_4)\le e^{-(\log g(\d^*/\g))\left(1-\frac{1}{\log\log
g(\d^*/\g)}\right)}$, 
for all $1\le i\le \bar N-1$,
for all $n''_i,n'_i$,  when  $ \eta(b_i,m)=-\eta(b_i,m^*)$  we have 
$$
\frac 
{ Z_{[n''_i, n'_{i+1}]}^{0 ,0}
\left (\PP^\r_{ [n''_i, n'_{i+1}]} (m,\z_5) \right)}
{ Z_{[n''_i, n'_{i+1}]}^{0 ,0}
\left (\PP^\r_{  [n''_i, n'_{i+1}]} (m^*,\z_5 ) \right)  }
= 
e^{\pm \frac {\b}{\g} \frac{20 V(\b,\th)}{ (g(\d^*/\g))^{1/4(2+a)}}}
e^{ \frac \b \g 
\frac{\tilde u(r_i)-\tilde u^*_\g(r_i)}{2\tilde m_\b}
\left[\sum_{\a:\,\e\a\in [r_i, r_{i+1})}\chi(\a)\right]}
\Eq(PP.24bis)
$$
where we have an upper bound for $\pm=+$ and a lower bound for $\pm=-$. }
\vskip0.3cm 
\noindent 
\proof
Recalling  that $ \g x_i=  r_i$, see   \eqv (Ma.2),   we need to  estimate  
$$
\KK(i,n''_i,n'_{i+1})=\sum_{\a:\,\e\a\in\g[x_i,n''_i]}\chi(\a)+
\sum_{\a:\,\e\a\in\g[n'_{i+1},x_{i+1}]}\chi(\a)
\Eq(diff)
$$
  uniformly  with respect to $1\le i\le \bar N$, $x_i$, $n''_i,
n'_{i+1}$. It is enough to  estimate 
$$
\max_{1\le i\le \bar N}\max_{x_i}\max_{1\le \ell \le L_1+\r\g^{-1}}
\Big|\sum_{\a:\e\a\in\g[x_i, \ell +x_i]}\chi(\a)\Big|. 
\Eq(diff1)
$$
However since the  point $x_i$ might be random and depending on  $\chi(\a)$, a little
care is needed.  An  upper bound for \eqv(diff1) is clearly   
$$
\KK(Q,L_1,\r,\e)\equiv \max_{\a_0: \e\a_0\in
 [-Q ,+Q   ]}\,\,
\max_{1\le \e  \bar \alpha \le  \g L_1   +\r }  
\Big|\sum_{\a=\a_0}^{\bar \a }\chi(\a)\Big|.
\Eq(diff2)
$$
Using, Levy inequality,  \eqv(PP.10) and exponential Markov inequality, one has
$$
\eqalign{
&\P\left  [\KK(Q,L_1,\r,\e)\ge 2V(\b,\th) \sqrt{2(\g L_1+\r)\log (g^5(\frac {\d^*} \g))}\right ]
\cr
&\quad \le \frac{2Q+1}{\e}\P\left  [\max_{1\le \e  \bar \alpha \le  \g L_1   +\r }   
\Big|\sum_{\a=\a_0}^{\bar \a }\chi(\a)\Big|\ge 2V(\b,\th) \sqrt{2(\g L_1+\r)\log (g^5(\frac {\d^*} \g))}
\right ] \cr
&\quad \le \frac{2Q+1}{\e} \frac{1}{g^5(\frac {\d^*} \g)}\le
e^{-\log g(\frac {\d^*} \g) \left(1-\frac{1}{\log\log
g(\frac {\d^*} \g)}\right)},
}\Eq(diff3)
$$
where we have used \eqv(epsilon) and \eqv(queue) at the last step.
Recalling that $0\le L_1+\ell_0\le \r\g^{-1}$ and \eqv(rho), one has
$$
2(\g L_1+\r)\log (g^5(\frac {\d^*} \g)) \le 4\r \log (g^5(\frac {\d^*} \g)) \le
4 \left (\frac{5} {g(\frac {\d^*} \g)}\right)^{1/(2+a)}\log (g^5(\frac {\d^*} \g))\le 
5^2\left (\frac{1} {g(\frac {\d^*} \g)}\right)^{1/2(2+a)}.\Eq(diff4)
$$
that entails \eqv(PP.24bis) after an easy computation. \eop 
 
Next we estimate  the remaining type of ratio in \eqv (gene1).
Recall that $n'_i\in [\ell_i-\bar R_2-L_1, \ell_i -\bar R_2]$ and $n''_i\in
[\ell_i+\bar R_2,\ell_i+\bar R_2+L_1]$ with $\bar R_2=R_2+\ell_0$, where $\ell_0$
is defined in \eqv(shrink2).
 
\vskip0.5cm 
\noindent {\bf \Lemma(phases)}
{\it On $\O_1\setminus \O_3 $, choosing the parameters as in
Subsection 2.5, for all $1\le i\le \bar N$,
in the case 1, we have
$$
\frac 
{Z_{[n'_i+1,n''_i-1]}^{m,m} 
\left( \PP^\r_{  [n'_i+1,n''_i-1]} (m,\ell_i,i) \right ) }
{ Z_{[n'_i+1,n''_i-1]}^{m^*,m^*}
\left  ( \PP^\r_{ [n'_i+1,n''_i-1]} (m^*,\ell_i,i) \right )  }
=e^{-\frac{\b}{\g}(\FF^*\pm 32\th(R_2+\ell_0+L_1)\sqrt{\frac{\g}{\d^*}})}.
\Eq(tipo21bis)
$$
In the case 2, we have
$$
\frac 
{Z_{[n'_i+1,n''_i-1]}^{m,m} 
\left( \PP^\r_{  [n'_i+1,n''_i-1]} (m,\ell_i,i) \right ) }
{ Z_{[n'_i+1,n''_i-1]}^{m^*,m^*}
\left  ( \PP^\r_{ [n'_i+1,n''_i-1]} (m^*,\ell_i,i) \right )  }
=e^{+\frac{\b}{\g}(\FF^*\pm 32\th(R_2+\ell_0+L_1)\sqrt{\frac{\g}{\d^*}})}.
\Eq(tipo22)
$$
In the case 3, we have
$$
\frac 
{Z_{[n'_i+1,n''_i-1]}^{m,m} 
\left( \PP^\r_{  [n'_i+1,n''_i-1]} (m,\ell_i,i) \right ) }
{ Z_{[n'_i+1,n''_i-1]}^{m^*,m^*}
\left  ( \PP^\r_{ [n'_i+1,n''_i-1]} (m^*,\ell_i,i) \right )  }
=e^{\pm \frac{\b}{\g}(32\th(R_2+\ell_0+L_1)\sqrt{\frac{\g}{\d^*}})}.
\Eq(tipo23)
$$
}
\vskip .5truecm
\noindent

\Remark(elleunter)  Note that here one needs to have
$L_1\sqrt{\frac{\g}{\d^*}}\downarrow 0$.

\vskip0.5cm 
\noindent
The proof of   \eqv(tipo21bis)  and \eqv (tipo22) follows from  Lemma
\eqv (T1).
The \eqv(tipo23) is a  consequence of
\eqv(tipo21bis) and \eqv(tipo22). 
Next we  estimate from below the r.h.s. of  \eqv (PPPP.1). 
 \smallskip
{\noindent \bf \Lemma(32DD)  }
{\it Under the same hypothesis of   Proposition \eqv(31) and on the   probability
space $\O_1\setminus   ( \O_3 \cup \O_4) $,
for  $\z_5$  as in \eqv (zeta51), 
we have 
 $$ \eqalign { & 
\frac{Z^{0,0}_{[v_1-1,v_2+1]}  
\Big(\PP^\r_{[v_1,v_2]} (m), 
\eta^{\d,\z_5}(v_1-1)=\eta(v_1-1,m^*), 
\eta^{\d,\z_5}(v_2+1)=\eta(v_2+1,m^*)\Big)}
{Z^{0,0}_{[v_1-1,v_2+1]}
\Big(
\PP^\r_{[v_1,v_2]} (m^*), 
\eta^{\d,\z_5}(v_1-1)=\eta(v_1-1,m^*), 
\eta^{\d,\z_5}(v_2+1)=\eta(v_2+1,m^*)\Big)}
\cr & \ge \left(e^{\frac{\b}{\g} (\AA+\EE_1)}+e^{-\frac{\b}{\g}\AA_2}\right)^{-1}
}\Eq(up2) 
$$
where $\AA$, $\EE_1$, and $\AA_2$ are defined in \eqv(penu2),
\eqv(errors),  and \eqv(A2) respectively.
}

\proof 
  Obviously  one can get     the lower bound  simply proving an upper bound for the inverse
  of  l.h.s.  of  \eqv (up2), i.e.  
$$
\frac
{Z^{0,0}_{{[v_1-1,v_2+1]}}
\Big(\PP^\r_{{[v_1,v_2]}} (m^*),\,   
\eta^{\d,\z_5}(v_1-1)=\eta([v_1-1,m^*), \,
\eta^{\d,\z_5}(v_2+1)=\eta(v_2+1,m^*)\Big)}
{Z^{0,0}_{[v_1-1,v_2+1]}
\Big(\PP^\r_{[v_1,v_2]} (m),\,   
\eta^{\d,\z_5}(v_1-1)=\eta(v_1-1,m^*), \,
\eta^{\d,\z_5}(v_2+1)=\eta(v_2+1,m^*)\Big)}.
\Eq(low)
$$
Note that  $\eta(v_1-1,m^*)= \eta(v_1-1,m)$ and $\eta(v_2+1,m^*)= \eta(v_2+1,m)$
and in the proof of the upper bound, see Lemma \eqv (32D),   we never 
used that  $m^*$ in the denominator  is the one given in Theorem \eqv (copv). 
Then \eqv (low) is equal to 
$$
\frac
{Z^{0,0}_{{[v_1-1,v_2+1]}}
\Big(\PP^\r_{{[v_1,v_2]}} (m^*),\,   
\eta^{\d,\z_5}(v_1-1)=\eta([v_1-1,m), \,
\eta^{\d,\z_5}(v_2+1)=\eta(v_2+1,m)\Big)}
{Z^{0,0}_{[v_1-1,v_2+1]}
\Big(\PP^\r_{[v_1,v_2]} (m),\,   
\eta^{\d,\z_5}(v_1-1)=\eta(v_1-1,m), \,
\eta^{\d,\z_5}(v_2+1)=\eta(v_2+1,m)\Big)}.
\Eq(low1)
$$
Then by Lemma \eqv (32D)   we obtain \eqv (up2).    \eop 
 
 \vskip0.5cm
 \noindent { \bf Proof of Proposition \eqv (31)} 
To prove \eqv(L.66), we use Lemma \eqv (32), then  Lemma \eqv (32DD)
to get a  lower
bound and Lemma \eqv (32D)  and Corollary \eqv (diff5) to get an  upper bound. 
For the lower bound we get applying  \eqv (PPPP.1) and \eqv(up2)  
 $$
\mu_{\b,\th,\g}\big(\PP^\r_{[q_1,q_2]}(m)\big) \ge 
e^{-\frac{\b}{\g}(4\z_5+8\d^*)}
\left({1-2K(Q) e^{-\frac{\b}{\g} \frac{1}{g(\d^*/\g)}}}-
2e^{-\frac{\b}{\g}\frac{\k(\b,\th)}{16}\d\z_5^3}\right) 
 \left( e^{\frac{\b}{\g} \AA}     e^{ \frac{\b}{\g} \EE_1}+e^{-\frac{\b}{\g}\AA_2}\right)^{-1}.
\Eq  (merc7)$$

 For the upper bound we get 
$$
\mu_{\b,\th,\g}\big(\PP^\r_{[q_1,q_2]}(m)\big)
\le
e^{-\frac{\b}{\g} \AA}  e^{+\frac{\b}{\g}\EE_1}
+ 2e^{-\frac \b\g \AA_2}
\Eq(ulti)
$$
where $\AA_2$ is defined in \eqv(A2).
To get \eqv(L.66) from \eqv(ulti), one needs  $\AA_2>\AA$, this will
be a consequence of an upper bound on $\AA$ and a lower bound on $\AA_2$.
We start estimating  the terms of $\AA$.
We easily obtain   
$$  \frac{\FF^*}{4\tilde m_\b}
\sum_{-L\le r \le L} \left[\|D  u(r)\|_1-\|D  u^*_\g(r)\|_1\right]
\le   \FF^*  \left [  N_{[-L,L]}(u) +  N_{[-Q,Q]}(
u^*_\g)\right].
\Eq(saltieffestar)
$$
We  use  that  
$N_{[-Q,Q]}( u^*_\g) \le K(Q)$, see \eqv (RAP.1),  where $K(Q) $ is given in  \eqv (M.1). 
If $ L$ is  finite for all $\g$, then  $ N_{[-L,L]}(u)$ is bounded since $ u \in BV_{loc}$. 
When $L$ diverges as $Q$  when $ \g \downarrow  0$ from the assumption   \eqv(condu) we have that 
$$   \bar N \le N_{[-L,L]}(u) +  N_{[-Q,Q]}( u^*_\g) \le [ f(Q)+1]K(Q) 
\Eq (merc3)
$$
where we set $$
f(Q)=e^{(\frac{1}{8+4a}-b)(\log Q)(\log\log Q)}. 
\Eq(effequeue)
$$
The second term  of   $\AA$  can be estimated as
  $$
\eqalign { 
\left | 
\sum_{i=1}^{\bar N} \frac{\tilde u(r_i)-\tilde u^*_\g(r_i)}{2\tilde m_\b}
\left[\sum_{\a:\,\e\a\in [r_i, r_{i+1})}\chi(\a)\right] \right |
&\le \bar N \max_{\{-\frac Q \e \le \a_0  \le \frac Q \e\} } 
\max_{\{ \a_0\le \bar \a \le \frac Q \e\}}
\left | \sum_{\a= \a_0} ^{\bar \a  }\chi(\a)\right |
\cr &  \le
2 \bar N  \max_{\{ - \frac Q \e \le \bar \a \le \frac Q \e\}  }  
\left | \sum_{\a=- \frac Q \e} ^{\bar \a  }\chi(\a)\right |. 
}
\Eq(diff2b)
$$
To  estimate the last term in \eqv(diff2b),   we use  Levy inequality, \eqv (PP.10)
and exponential Markov inequality  to get 
$$
\P\Big  [  \max _{\{ - \frac Q \e \le \bar \a \le \frac Q \e\}  }  
\left | \sum_{\a=- \frac Q \e} ^{\bar \a  }\chi(\a)\right |  \ge \sqrt
3 V(\b,\th)
\sqrt{[2 Q+1] \log (g(\frac {\d^*} \g))} \Big  ]
\le 4e^{-\log (g (\frac {\d^*} \g))}= \frac 4 {g (\frac {\d^*} \g)}. 
\Eq(diff3a)
$$
Denote $\O_5$ the  probability space for which    \eqv (diff3a) holds.
Then for  $\o \in    \O_1\setminus (\O_3 \cup \O_4 \cup \O_5) $ and
$\g_0$ small enough, one has
$$
\AA\le 
2 [ f(Q)+1]  K(Q)
V(\b,\th) \sqrt{(2 Q+1) \log (g (\frac {\d^*} \g))}\le \bar c(\b,\th) f(Q)Q^3
\Eq (merc2)
$$
for some $\bar c(\b,\th)$. The last inequality in \eqv (merc2) is obtained from  the choice of $f(Q)$ in  \eqv(effequeue), 
the one of $K(Q)$ in \eqv (M.1) and the choice of $Q$ in \eqv(queue). 
Namely, from   \eqv(queue)    $Q^2g(\d^*/\g)\le g^2(\d^*/\g)$.
Notice that in $\AA_2$, see \eqv (A2), $L_1$  enters. We   make  the following choice of $L_1$ 
$$
L_1=\left (g(\frac{\d^*}{\g})\right)^{19/2}.
\Eq(elle1def)
$$
This choice   satisfies the
requirement in  Proposition \eqv (32), i.e.  $L_1 < \frac \r \g$,  see
\eqv (rho). 
Furthermore as in [\rcite{COPV}] we make the choice 
$$
\z_5=\frac{1}{2^{18} c^6(\b,\th)}\frac{1}{g^3(\d^*/\g)}
\Eq(zeta5)
$$
for some constant $c(\b,\th)$.  Obviously \eqv (zeta5) satisfies  requirement \eqv (zeta51) provided $\z$ is chosen
according \eqv (TE.1).  Since  $  Q = g(\d^*/\g)^{\frac 1 { \log\log g(\d^*/\g)}} $,  see \eqv(queue), we have   
$\log g(\d^*/\g)=(\log Q) (\log\log g(\d^*/\g))$.   Iterating this
equation, for  $\g_0$  small enough to have  $\log\log\log
g(\d^*/\g)>0$, one gets 
$$
\log g(\d^*/\g)=(\log Q)(\log\log Q)\Big(1+\frac{\log\log\log
g(\d^*/\g)}{\log\log g(\d^*/\g)-\log\log\log g(\d^*/\g)} \Big)\ge (\log Q)(\log\log Q).
\Eq(invert)
$$
Therefore,  recalling \eqv(5.97037) and using \eqv(effequeue) one can
check that 
$$
L_1  \frac{\k(\b,\th)}{16}\d\z_5^3 \ge c(\b,\th) f(Q)Q^3.
\Eq(anelle1)
$$
It is not difficult to check that \eqv(anelle1) implies 
$$
L_1 \frac{\k(\b,\th)}{16}\d\z_5^3
> 2\g \log  \bar N+ \g \bar N \log \frac{\r}{\g}  + 8\d^*+4\z. 
\Eq(mart2)
$$
Therefore, recalling \eqv(A2), \eqv(anelle1) entails  $\AA_2 > \AA$ and finally 
 one gets 
$$
\mu_{\b,\th,\g}\big(\PP^\r_{[q_1,q_2]}(m)\big)
\le
e^{-\frac{\b}{\g} \AA} e^{ \frac{\b}{\g} \EE_1}
\left(1+2e^{-\frac{\b}{\g} \left\{ L_1
\frac{\k(\b,\th)}{16}\d\z_5^3\right\}}\right).
\Eq(postulti)
$$
It remains to check that $\EE_1\downarrow 0$. Recalling
\eqv(rho), one has $\g \log (\r/\g)\le (g(\d^*/\g))^{-1}$. Recalling
\eqv(erre2) one has $( R_2+\ell_0)\sqrt{\g/\d^*}\le
(g(\d^*/\g))^{-1}$. Therefore, using \eqv(invert), \eqv(effequeue) and
recalling that $0<b<1/(8+4a) $, see Proposition \eqv(31),  one has
$$
\EE_1 \le
K(Q)(f(Q)+1) \big[\z_5
+32 \th L_1\sqrt{\frac{\g}{\d^*}}
+\frac{c(\b,\th)}{(g(\d^*/\g))^{1/(8+4a)}}\big]\le (g(\d^*/\g))^{-b}.  
\Eq(mart4)
$$
So one gets the  upper bound in \eqv(L.66). 
Recalling \eqv(merc7),  it is easy to get the corresponding lower bound. \eop

\vskip0.5cm   
\noindent{\bf \Lemma (T1)} { \it  For $\o \in \O_1$ and choosing the parameters as in Subsection 2.5 we have 
  $$ e^{-\frac{\b}{\g}(\FF^* + 32\th R_2\sqrt{\frac{\g}{\d^*}})} \le 
\frac
{Z^{m,m}_{[\ell,\ell+2R_2]}\left(\WW_1([\ell,\ell+2R_2],R_2,\z)\right)}
{Z^{m^*,m^*}_{[\ell,\ell+2R_2]}\left(\PP^\r_{[\ell,\ell+2R_2]}(m^*)\right)}
\le  
e^{-\frac{\b}{\g}(\FF^*- 32\th R_2\sqrt{\frac{\g}{\d^*}})}. 
\Eq(PP.32)
$$}
The Lemma is proven   in [\rcite{COPV}], see Lemma 7.3. 
The random field is estimated as in  Lemma 3.3 (the rough
estimates).    The upper bound can be easily obtained  since $\FF^*$
is the minimum amount of energy needed to
go   from one phase to the other, see \eqv (min1).   More care should
be taken to show the lower
bound, see formula (7.34) of [\rcite{COPV}].

Next we   summarize  in Lemma \eqv (33) the estimates  needed to prove
Lemma \eqv (33a).
Let $\II$ be the interval that gives rise to an elongation. 
Denote by $\sign(\II)=  1 $, if $\II$ gives rise to 
a positive elongation, $\sign(\II)=  -1 $ in the other case.

\vskip0.5cm   
\noindent{\bf \Lemma (33)} { \it Let  $\O_1$ be the probability space of Theorem  \eqv (copv), 
let $ \II \subset \D_Q$  be an interval 
that gives rise  to  an elongation.  Then for any interval $ I \subset 
\g^{-1}\II  $ we have} 
$$ 
\frac 
{Z_{I}^{0,0} \left (\eta^{\d,\z_1} (\ell)=- \sign(\II),\, 
\forall \ell \in I \right)} 
{ Z_{I}^{0,0} \left  ( 
 \eta^{\d,\z_1} (\ell)= \sign(\II),\, \forall \ell \in I \right)  }= 
e^{-  \sign(\II)\frac {\b}{\g}  \sum_{\a \in \e^{-1} \g I}  \chi (\a)} \,\,
\frac {Z_{I,0}^{0,0} \left (\eta^{\d,\z_1} (\ell)=- \sign(\II),\, 
\forall \ell \in I \right)} 
{ Z_{I,0}^{0,0} \left  ( 
 \eta^{\d,\z_1} (\ell)= \sign(\II),\, \forall \ell \in I \right)  }.
\Eq (es.1) $$
{\it On $\O_1$, 
the last ratio satisfies: For all $I\subset \g^{-1} \II\subset \g^{-1}[-Q,+Q]$}
$$
\left|\log 
\frac {Z_{I,0}^{0,0} \left (\eta^{\d,\z_1} (\ell)=- \sign(\II),\, 
\forall \ell \in I \right)} 
{ Z_{I,0}^{0,0} \left  ( 
 \eta^{\d,\z_1} (\ell)= \sign(\II),\, \forall \ell \in I \right)  }\right|
\le \frac {\b}{\g} \frac{ 1}{4c^2(\b,\th) g(\d^*/\g)}
\Eq(es.10)
$$
{\it where $g$ is the  function given in Subsection 2.5 and $c(\b,\th)$  is some positive constant
that depends only on $\b,\th$. }
\smallbreak

The proof  has been  done in [\rcite{COPV}], see the proofs of Lemma 6.3   and  
of Proposition 4.8 there.  
It consists essentially in extracting the leading stochastic part 
and estimating the remaining term by using a classical deviation inequality for Lipschitz
functions of Bernoulli random variables. The corresponding Lipschitz norms
are estimated using the cluster expansion. The proof is however long
and tedious.

\vskip 1truecm
\chap {5    Probability estimates  and  Proof of Theorem \eqv (2)}5
\numsec= 5
\numfor= 1
\numtheo=1

\vskip0.5cm 
  In this section we  prove the probability estimates   needed 
for  proving   the main results stated in Section 2. The proof of
Theorem \eqv (2) is given after  Lemma \eqv(45).
This section is rather long and we  divided into  several subsections.
In the  first by using  a simple and direct application of the Donsker invariance
principle in the Skorohod space,  we prove that 
the main random  contribution identified in \eqv(4.800) 
suitably rescaled,  converges in law to  a
bilateral Brownian  process  (BBM), see \eqv (DPP.1).

In the second   subsection we recall the construction  done by 
Neveu-Pitman,   [\rcite{NP}],  to determine the  $h$-extrema for a bilateral
Brownian motion  and then we  adapt  it  to  the random walk
corresponding to the previous random contribution. 
In Subsection 5.3   we state definitions and main properties of the
maximal $b$--elongations  with excess $f$  introduced  in   [\rcite{COPV}].
In Subsection  5.4,  which is the most involved,  we  identify   them
with the $h$--extrema of Neveu--Pitman  by  restricting   suitably the
probability space we are working on.
Here   $b$, $f$, and $h$ are positive constant which will be specified. 
In the last section  we
present rough estimates on the number of renewals up to time $R$,
needed to prove the Theorem \eqv  (copv).

\medskip
\noindent{\bf 5.1. Convergence  to a  Bilateral Brownian Motion}
\medskip

Let
$ \e \equiv \e (\g)$, $\lim_{\g \to 0}  \e (\g) =
0$,
$\frac \e {\g^2} > \frac { \d^*} \g $,  so that
each block of length  $\frac
\e {\g^2}$   contains at least one block
$A(x)$ (see section 2.2 )
; to avoid rounding problems it is assumed that
$\e/\g\d^*\in \N$, and that the basic initial partition $A(x)\colon\, x \in
C_{\d^*}(\R)$ is a refinement of the present one,
 see \eqv(epsilon)  for the actual choice of $\e$.
Denote    by  $\{\hat W^\e  (t); t\in \R \} $  the  following
 continuous time   random walk: 
$$
\hat W^\e   (t) \equiv   \left \{ \eqalign { &  V^\e_1   (t)=\frac 1
{\sqrt { c(\b,\th, \g /\d^*)} }  \sum_{\a=1}^{ [  \frac t \e ]   }
\chi (\a),
\quad \,\,  t\ge  \e \, ; \cr &
0,   \qquad -\e \le  t\le  \e\,;  \cr & 
 V^\e_2   (-t)=\frac 1
{\sqrt { c(\b,\th, \g /\d^*)} }  \sum_{ \a=-[\frac t \e ]}^{\a=-1}
\chi (\a), \,  t\le -\e.\cr } \right.
\Eq(N.9)
$$
Here  $[x]$ is  the integer part of $x$ and $\chi(\a)$ was defined in
\eqv(3.1) for all $\a\in \Z$.
Definition \eqv (N.9)  allows to  see  $\hat W^\e  (\cdot) $ as a trajectory in the
space  of real functions   on  the line  that are right continuous and have
left limit, {\it i.e}  in the Skorohod space $D(\R,\R)$ endowed with
the Skorohod topology.
To define  a metric that makes it separable and complete,
let us denote $\L_{\rm Lip}$ the set of strictly
increasing Lipschitz continous function $\l$ mapping $\R$ onto $\R$ such
that
$$
\|\l\|=\sup_{s\ne t}\Big| \log \frac{\l(t)-\l(s)}{t-s}\Big|<\infty.
\Eq(Sk1)
$$
For $v \in D(\R,\R)$ and $T\ge 0$,  let us define
$$
v^T(t)=\cases {v(t\wedge T),&if $t\ge 0$;\cr
v(t\vee (-T)),&if $t<0$.\cr}
\Eq(Sk2)
$$
Define for  $v$ and $w$ in $D(\R,\R)$ 
$$
d(v,w)\equiv \inf _{\l\in \L_{\rm Lip}}\left [\|\l\| \vee
\int_{0}^{\infty} e^{-T} 
\sup_{t\in \R}(1\wedge (|v^T(t)-w^T(\l(t))|)\,{\rm d}T\right]. 
\Eq(Sk4)
$$
Note that for a given $T\in R$, the quantity
$|v^T(t)-w^T(\l(t))|$ is constant for $t>T\vee \l(T)$ and for
$t<(-T)\wedge (\l(-T))$, therefore the previous supremum over $t\in
\R$ is merely over $(-T)\wedge (\l(-T))\le t\le T\vee \l(T)$.
See [\rcite{Bill}] chapter 3 or [\rcite{EK}] chapter 3  where the case
of $D[0,\infty)$ is considered  with  all the needed details. 
Let us define the bilateral  Brownian motion $W =(W(t); t\in \R)$ by
$$
W(t)\equiv W_t= \left\{ \eqalign { &  B_1(t) \qquad \qquad t \ge 0 \cr &  B_2(-t) \qquad \qquad t \le 0. }
\right.
\Eq(BM.1)
$$
with $(B_1(t), t\ge 0)$ and $(B_2(-t), t\le  0)$ two  independent standard
Brownian motions. 
Note that $E[(W(t))^2]=|t|$ for all $t\in \R$, in
particular $W(0)=0$, and when $s\le 0\le t$,  $E[(W(t)-W(s))^2]=t-s$.
 Since  $\chi(\a)$ depends on $\e=\e(\g)$,
we need the following generalization of the Donsker Invariance Principle 
that can be proved following step by step the proof of Billingsley 
[\rcite{Bill}] pg 137. 
\medskip
\noindent{\bf Theorem (Invariance Principle) } {\it 
Let $\e \equiv \e (\g)>0$ so that    $\frac \e {\g^2} > \frac { \d^*} \g $,  
 $\lim_{\g \to 0}  \e (\g) = 0$.  Let   
$\PP^\e$ be the measure induced by $\{\hat  W^\e(t), t\in \R \} $
on
$  D \left (\R, \BB(\R)\right )  $
Then as $\g\downarrow 0$,
$ \PP^{\e} $ converges weakly to the Wiener measure
$\PP $, under which the coordinate
mapping process $W(t)$, $t \in \R$ is a
bilateral Brownian motion.
} 
\medskip
\noindent 
\remark One can wonder about the coherence between
the fact that  $\hat W^\e  (t) =0$,  for $-\e \le  t\le  \e$ in  \eqv(N.9)
and  $\chi(0)\neq 0$.  However we have  
$
 \chi(0)\equiv \g \sum_{x:\d^*x \in \tilde A_{\e/\g}(0)}
X(x)\1_{\{p(x)\le (2\g/\d^*)^{1/4}\}}.  $
\medskip 
\noindent{\bf \Lemma (52)  }   
$$ \lim_{\g \to 0} \P [ \chi(0) =0]  = 1
\Eq(chizero)
$$ 
\medskip
\noindent 
\proof
Using \eqv(PP.10),  one gets immediately for    $c>0$, 
$
\P [  |\chi(0) | \ge c] \le 2 e^{- \frac {c^2}  {3 \e V^2} }
$
which  implies \eqv(chizero), since $\e=\e(\g)\downarrow 0$. \eop
\medskip
\noindent
The  following result is immediate 
\vskip0.5cm 
\noindent{\bf \Lemma (52b)  }   {  \it Set   $\eta = \pm 1 $,
$I= [\frac a \g, \frac b \g ]$ (macro scale),  $a$ and $b$ in $\R$.   Then,  see  \eqv(4.800), we obtain } 
$$ \lim_{\g \to 0}   \left [ -\eta \g \D^\eta \GG(m^{\d^*}_{\b,I})\right] \mathrel{\mathop=^{\rm Law}}
V(\b,\th) [W(b)-W(a)]. \Eq (DPP.1)$$ 
\vskip0.5cm  \noindent 
\proof
 Recalling  \eqv(N.9), for $\eta = \pm 1 $ and $I= [\frac a \g, \frac b \g ]$
in macroscopic  scale with   $0<a<b$ or $a<b<0$
 one gets the following 
$$
 \g \D^{\eta} \GG(m^{\d^*}_{\b,I})=-\eta \g 
\sum_{x\in \CC_{\d^*}(I) } X(x)=  -\eta \sum_{\a =[\frac a \e]
}^{[\frac b\e ]}
\chi(\a)  = -\eta
\sqrt { c(\b,\th, \g/\d^*)}
\left [  \hat W^\e (b) - \hat W^\e (a) \right ].
\Eq (N.2)$$
 When  $0\in [a,b]$, we get
$$
\g \D^{\eta} \GG(m^{\d^*}_{\b,I})=-\eta 
\sqrt { c(\b,\th, \g/\d^*)}
\left [  \hat W^\e (b) - \hat W^\e (a) \right ]
-\eta \chi(0).
\Eq(N.2bis)
$$
 Therefore, using Lemma \eqv(52) to take care of the $\chi(0)$ term
and \eqv (2.21) we obtain \eqv  (DPP.1). \eop 
\medskip
\noindent
\remark Note that $I=[\frac a \g, \frac b \g ]$ corresponds to $\g I= [a,b]$ in the
Brownian scale, according to the notation in Subsection 2.2.
The \eqv(DPP.1) is the main reason to
have introduced the notion of   ``Brownian'' scale.
In this scale the main random contribution identified in \eqv(4.800)
becomes a functional of a bilateral Brownian motion. 

\medskip
\noindent
{\bf 5.2.  The Neveu-Pitman construction  of the  $h$--extrema for  the random walk  $ \{\hat W^\e\}$}
\medskip

We shortly recall the  Neveu-Pitman construction  [\rcite{NP}], 
used to determine the $h$--extrema  for   the bilateral  Brownian Motion  
$(W_t, t\in \R)$.  Realize it  as the coordinates of
the set $\O$ of real valued functions $\o$ on $\R$ which vanishes at the origin.
Denote by  $(\th_t, t \in \R)$, the flow of translation : $[\th_t
\o(\cdot)=\o(t+\cdot)-\o(t)]$ and by $\r$  the time reversal $\r\o(t)=\o(-t)$.
For $h>0$,  the  trajectory
$\o$ of the  BBM
 admits an $h$--minimum at the origin if $W_t(\o)\ge W_0(\o)=0$ for $t\in
[-T_h(\rho \o), T_h(\o)]$ where $T_h(\o)=\inf[t : t>0, W_t(\o)>h]$,
 and $-T_h(\r\o)=-\inf[t>0 : W_{-t}(\o)>h]\equiv\sup[t<0, W_{-t}(\o)>h]$.
The  trajectory
$\o$ of the  BBM   admits an
$h$--minimum (resp. $h$ maximum) at $t_0\in \R$ if $W\, o\, \th_{t_0}$
(resp. $-W\, o\, \th_{t_0}$)  admits an $h$ minimum at $0$. 

To define the point process of $h$--extrema for the BBM,
Neveu-Pitman consider first the one sided Brownian
motion $(W_t, t\ge 0, W_0=0)$, {\it i.e} the part on the right of the
origin of the BBM.    Denote  its running maximum by
 $$M_t =\big (\max(W_s\,;\,
0\le s\le t), t\ge 0\big) \Eq (E.1) $$ and define
$$
\eqalign{
&\t=\min(t\,;\,t\ge 0, M_t-W_t=h),\cr
&\b=M_{\t},\cr
&\s=\max(s\,;\,0\le s\le \t\,,\, W_s=
\b).\cr
}\Eq(NP1)
$$
The stopping time  $\t$ is the first time  that the Brownian motion achieves a
drawdown of size $h$,  see [\rcite{T},\rcite{Wi}]. Its Laplace transform is given by 
 $\E[\exp(-\l\t)]=(\cosh (h\sqrt{2\l}))^{-1}$, $\l>0$.  This is  consequence  of  the celebrated  L\'evy
Theorem [\rcite{KS}]     which states  that
$(M_t-W_t; 0\le t<\infty)$ and $(|W_t|; 0\le t<\infty)$ have the same
law.    Therefore $\t$ has the same law as the
first time a reflected Brownian motion reaches  $h$. The Laplace
transform of this last one  is obtained applying  the optional sampling
theorem to the martingale $\cosh(\sqrt{2\l}|W_t|)\exp(-\l t)$.
Further Neveu and Pitman   proved that 
$(\b,\s)$ and $\t-\s$ are independent and give the corresponding
Laplace transforms. In particular one has
$$
\E[e^{-\l\s}]=(h\sqrt{2\l})^{-1}\tanh(h\sqrt{2\l}).
\Eq(la)
$$
Now  call $\t_0=\t,\,\b_0=\b,\,\s_0=\s$ and 
define  recursively   $\t_n,\b_n,\s_n\, (n\ge 1)$, so that
$(\t_{n+1}-\t_n, \b_{n+1},\s_{n+1}-\t_n)$ is the $(\t,\b,\s)$--triplet
associated to the Brownian motion
$\big((-1)^{n-1}(W_{\t_n+t}-W_{\t_n}),\, t\ge 0\big)$.
By construction, for  $n\ge 1$,  $\s_{2n}$ is the time of  an  $h$-maximum
and for $n\ge 0$, $\s_{2n+1}$ is the time of a $h$-minimum.
Note that since we have considered just the part on the right of the
origin, in general $\s_0$ is not   an $h$ maximum.
The definition only requires $W_t \le W_{\s_0}$ for $t \in
[0, \s_0)$,  therefore $W_{\s_0}=W_{\s_0}-B_0$ could be smaller than $h$.   
The trajectory of  the BBM on the left of  the origin   will determine whether $\s_0$ is or is
not  an $h$--maximum. 
From the above mentioned fact that $(\b,\s)$ and $\t-\s$ are
independent, it follows that the variables $\s_{n+1}-\s_{n}$ for $n\ge
1 $ are independent with  Laplace transform 
$(\cosh (h\sqrt{2\l}))^{-1}$.  In this way  Neveu and Pitman  define   a renewal process
on $\R_+$, with a delay distribution, {\it i.e.} the one of $\s_0$,  that have Laplace transform
\eqv(la). 
 
Since the times of $h$-extrema
for the   BBM depend  only on its increments,
these times should  form a stationary process on $\R$. 
The above one side construction does not provide   stationary   on the positive real axis $\R^+$ 
since the delay distribution is not the one of the limiting
distribution of the residual life as it should  be,  see [\rcite{As}] Theorem
3.1.   In fact the Laplace transform of limiting distribution of the residual life is given
by
\eqv(M.6) which is different from
\eqv(la).

There is a standard way to get a stationary renewal process.
Pick up an $r_0>0$, translate the origin
to $-r_0$ 
 and repeat  for $(W_{t+r_0},\, t>-r_0)$ the above construction.
One gets $\s_0(r_0)$ and  the sequence of point of $h$-extrema  
$(\s_n(r_0),n\ge 1)$. 
Let  $\nu(r_0)\equiv \inf (n>0:\, \s_n(r_0)>0)$ be  the number of renewals
up to time $0$ (starting at $-r_0$). In this way, 
$\s_{\nu(r_0)}(r_0)$ is  the residual life at ``time''
zero for the Brownian motion  starting at $-r_0$. 
So taking $r_0 \uparrow \infty$, the distribution of $\s_{\nu(r_0)}(r_0)$
will converge to the one of the residual life and using
[\rcite{As}], Theorem 3.1, one gets a stationary renewal process on
$\R^+$. 
So conditionally on $\s_1(r_0)<0$, define $S_i(r_0)=\s_{\nu(r_0)+i-1}(r_0)$ for
all $i\ge 1$. 
 Then since the event
$\{\s_1(r_0)<0\}$ has a probability that goes to 1 as $r_0\uparrow \infty$, one gets,   
as $r_0\uparrow \infty$,  a stationary renewal process on $\R_+$   as well  
on $\R$.
Since the
Laplace transform of the inter--arrival time distribution is 
$(\cosh (h\sqrt{2\l}))^{-1}$, one gets easily  that 
the Laplace transform of the distribution of $S_1$ (and also of $S_0$)
is \eqv(M.6). 

\vskip .5truecm

With this in mind we start the construction for the random walk
$ \{\hat W^\e\}$.  Denote
 $$
V^\e (t)   = \left \{ \eqalign { &  V^\e_1   (t)=\frac 1
{\sqrt { c(\b,\th, \g /\d^*)} }  \sum_{\a=1}^{ [  \frac t \e ]
} \chi (\a)   \qquad t\ge  \e, \cr &
0   \qquad 0  \le  t\le  \e  } \right.
\Eq (N.9a)
$$
and $\hat \FF^+_t$, $t\ge 0$  the associated $\s$- algebra.
 Define  
 the rescaled running maximum for $ V^\e (t)$, $t\ge0$ 
$$ \sqrt \e \hat M(n)= \max_{ 0 \le k\le n}  V^\e (k \e ).\Eq (30.1) $$
The $ \sqrt \e$ multiplying $\hat M(n)$ comes from the observation, see \eqv (2.20),  that
 $ \E\left [ \left ( \frac 1 {\sqrt \e}  V^\e (k
\e )\right )^2 \right] =k$.   
  For any $h>0$, define the $ \hat \FF^+_t$ stopping time 
$$
\hat \t_0 (\e)\equiv\hat \t_0 =\min \{ n \ge 0: \sqrt \e \hat M(n)- V^\e (n \e ) \ge h  
\}, \Eq (3.2) $$
$$\sqrt \e \hat \b_0(\e)\equiv \sqrt \e \hat \b_0 = \max \{  V^\e (k
\e ): 1
\le k\le \hat \t_0 \}  \Eq (3.3a)$$
and $$ \hat \s_0(\e) \equiv \hat \s_0=\max \{k: 1 \le k \le \hat \t_0;
V^\e (k \e )= \sqrt \e \hat \b_0 \}. \Eq (3.3) $$
By construction  
$$ \sqrt \e \hat \b_0  \equiv   \sqrt \e \hat M(\hat \t_0)=
\max_{ 0 \le k\le \hat \t_0}  V^\e (k \e )= V^\e (\hat \s_0 \e )
\ge V^\e(\hat
\t_0 \e ) +h . \Eq (E.1bis)$$
It follows from   the invariance principle and 
the continuous mapping theorem,
Theorem 5.2 of [\rcite{Bill}], that the joint distribution of
$$ \left [ \sqrt \e \hat M([\frac t \e]), \e \hat \t_0 (\e),
\sqrt \e \hat \b_0(\e),\e\hat \s_0(\e) \right ] $$
converges as $\e \to 0$, to the joint distribution of
the  respective quantities 
defined for a Brownian motion, see \eqv(NP1) {\it i.e}
$$ \left [  M_t, \t_0, \b_0, \s_0 \right ].  $$

Since $\hat \t_0 $ is a  $\hat \FF^+_t$ stopping time for
$( V^\e (t),t\ge0 )$ , the translated   and reflected  motion 
$ (-1)[V^\e (\e\t_0  +t) -  V^\e (\e\t_0 )]$,  for $t \ge 0$,  is a new
random walk independent of $( V^\e (t), 0 \le t \le
\e \t_0)$ on which we will iterate the previous construction. We have 
 $$ \eqalign { & 
\hat \t_1 (\e)\equiv\hat \t_1 =\min \{ n \ge \hat \t_0:
 \max_{ \t_0 \le k\le n} [-  V^\e (k \e )] + V^\e (n \e )
\ge h  
\} \cr &= 
\min \{ n \ge \hat \t_0:  -   \min_{ \t_0 \le k\le n}
V^\e (k \e )  + V^\e (n \e )
\ge h  
\} \cr &= 
\min \{ n \ge \hat \t_0:     \min_{ \t_0 \le k\le n}
V^\e (k \e )  - V^\e (n \e )
\le - h  
\}  } 
.\Eq (3.2a) $$
$$\eqalign { & \sqrt \e \hat \b_1(\e)\equiv \sqrt \e \hat \b_1 =
\max \{  (-V^\e (k \e )):   \hat \t_0 \le k \le \hat
\t_1
\}\cr & =  - \min \{  V^\e (k \e ):  \hat \t_0 \le k \le  \hat \t_1
\} }
\Eq (3.4b)$$
$$
\hat \s_1(\e) \equiv \hat \s_1=\max \{k: \hat \t_0 \le k \le \hat
\t_1;     -V^\e (k \e )= \sqrt \e \hat \b_1 \} \Eq
(3.4c) $$
 Now for any  $i\in \N$, we can iterate  the above procedure
to get as   Neveu and Pitman the  family 
$$
\left[\sqrt\e\tilde M([\frac{t}{\e}]),\e \hat\t_0(\e),\sqrt{\e}\hat\b_0(\e),\e\hat\s_0(\e),\dots,
\e \hat\t_i(\e),\sqrt{\e}\hat\b_i(\e),\e\hat\s_i(\e)\right].
\Eq(RF.1)
$$
Using again the invariance principle and the continuous mapping
theorem one gets that
$$
\lim_{\e\downarrow 0}\left[\e
\hat\t_i(\e),\sqrt{\e}\hat\b_i(\e),\e\hat\s_i(\e), i\ge 0,  i\in
\N\right]\mathrel{\mathop=^{\rm Law}} \left[\t_i,\b_i,\s_i, i\ge 0, i\in \N\right] ,
\Eq(RF.2)
$$
where the quantity in the right hand side of \eqv(RF.2) are the ones 
defined after \eqv(la).
 Let us note the following properties of the previous points of $h$--extrema.
By construction the random walk  satisfies,
in  the interval   $ [ \hat \s_0,
\hat \s_1 ] $,      the following :
 
\medskip  
  \noindent {\bf Property (4.A) }  In the interval $[ \hat \s_0, \hat
\s_1 ] $ we have 
$$  V^\e (\hat \s_1 \e)- V^\e (\hat \s_0 \e)  \le -h, \,
\qquad   V^\e (k \e )  - V^\e (k' \e )     < h  
\qquad \forall k'<k \in  [ \hat \s_0,
\hat \s_1 ],  
\Eq(3.5a) $$
$$    V^\e (\hat \s_1 \e)  \le V^\e (k \e ) \le  V^\e (\hat \s_0 \e) 
\qquad  \hat \s_0  < k < \hat \s_1  . \Eq(3.5b) $$

\vskip0.5cm  \noindent 
The first property in \eqv (3.5a) is easily obtained.
Namely   adding and subtracting $ V^\e ( \e \hat \t_0 )$ 
one has 
$$ [V^\e (\e \hat \s_1  )-  V^\e (\e \hat \t_0 )] +
[ V^\e (\e \hat \t_0 )- V^\e (\e\hat \s_0 )]  \le  -h $$
since $ [V^\e (\e\hat \s_1  )-  V^\e (\e\hat \t_0 )]\le 0$ and
by construction $ V^\e (\e\hat \s_0)-  V^\e (\e\hat \t_0) 
\ge  h $. The other properties are easily checked.
Properties similar to 
\eqv (3.5a) and  \eqv (3.5b) hold in the interval
$[\hat \s_{2i}, \hat \s_{2i+1}]$,  for $ i>0$.   
Namely   by construction $\hat \s_{2i}$ is a point of $h$-maximum and 
$\hat \s_{2i+1}$ is a point of $h$-minimum.  
Further, since by construction  $\hat \s_{2i-1}$ is a point of $h$-minimum and 
$\hat \s_{2i}$ is a point of $h$-maximum  in   the interval $[ \hat \s_{2i-1}, \hat \s_{2i} ] $,   
$i\ge1$,   we have the following:  
\medskip  
  \noindent {\bf Property (4.B) }  In the interval  $[ \hat
\s_{2i-1}, \hat \s_{2i} ] $,   $i\ge1$,  we have
$$ V^\e (\hat \s_{2i} \e)   -V^\e(\hat \s_{2i-1} \e) \ge h, \,
\qquad   V^\e(k \e) - V^\e ( k'\e)    > -h  
\qquad \forall k'<k\in  [ \hat \s_{2i-1}, \hat \s_{2i}],    
\Eq(3.6a) $$
$$   V^\e(\hat \s_{2i-1} \e)  \le   V^\e( k \e)  \le  V^\e(\hat \s_{2i} \e)
\qquad  \hat \s_{2i-1} < k < \hat \s_{2i}.
\Eq(3.6b) $$

\vskip0.5cm  \noindent 

Following  the Neveu--Pitman construction,  
one   translates the origin  of the random walk  $ \{V^\e\}$
to $-r_0$, being $r_0$ positive and large enough
and repeats the previous construction.  To obtain the
$h$--extrema  as in Neveu-Pitman  we should let  first  
$\e\to 0$, obtaining  by the Donsker invariance principle that
 $$   V_{ r_0}^\e (\cdot ) \equiv  V^\e (\cdot + r_0) \Eq (NN.1)
$$   converges  in law to the standard 
BM  translated by 
$-r_0$, then  $ r_0 \to \infty$.
However we cannot proceed in this way since to
  control  some 
probability estimates we need to have $\e$
small but  different from zero.
For the moment, the picture to have in mind is merely to take a
suitable $r_0=r_0(\g)$ that diverges when $\g\downarrow 0$. 
We denote by
$(\hat \s_{i}(r_0) =\hat \s_{i}(\e,r_0), i\ge 1, i\in \N) $  the points  of $h$--extrema for
$V_{r_0}^\e (\cdot )$.

\vskip0.5cm

\medskip
\noindent{\bf 5.3.  The maximal $b$  elongations  with  excess $f$  as defined in [\rcite{COPV}]}
\medskip
In this subsection we recall  definitions of the  maximal  elongations from
[\rcite{COPV}].
We  extract them from the first 5 pages of Section 5 of
[\rcite{COPV}],   with  different conventions that will be pointed out. 
This subsection is not completely self--contained since an involved
probability estimate done in [\rcite{COPV}], see \eqv(Zorro.4b) is  just
recalled. However if one accepts it, the rest is self--contained.  
In [\rcite{COPV}],  formula (5.3) we introduced  the following
$$
{\YY}(\a) \equiv \left \{ \eqalign {
& \sum_{\tilde \a \in [0, \a]}\chi(\tilde \a), \text {if} \a \ge 1;
\cr &0 \qquad   \text {if} \a=0 ;\cr & -\sum_{\tilde \a \in (\a, -1) }\chi (\tilde \a),  \text{if} \a < -1
}
\right.  \quad \a \in \Z.  \Eq(3.C6)
$$

\noindent{\bf \Definition (P39):}
{\it Given $b>f $ positive real numbers,
we say that an interval $[\a_1,\a_2]$ gives rise to a negative 
$b$--elongation with excess $f$, for $\YY(\a), \a \in \Z$ if
$$
 \YY  ( \a_2) - \YY (\a_{1})  \le -b -f;\,  \quad
\YY  ( y) - \YY  ( x)    \le b -f,  
\quad \forall x<y \in  [  \a_1, \a_{2}].
\Eq(3.3aabb) $$
We say that $[\a_1,\a_2]$  gives rise to a positive $b$--elongation
with excess $f$ if
$$ \YY  ( \a_{2}) - \YY (\a_{1})  \ge b +f;\,
\quad   \YY  ( y) - \YY  ( x)    \ge -b +f,  
\quad \forall x<y \in  [  \a_{1}, \a_{2}].    
\Eq(3.3ban) $$
In the first case we say that the sign of the $b$--elongation 
with excess $f$ is $-$; in the second case, $+$. 
}
\medskip
\noindent
\Remark (v1) To decide  if a given  interval $[\a_1,\a_2]$ gives rise
to a $b$--elongation with excess $f$ depends only on the variables
$\chi(\a)$ with $\a_1\le \a\le\a_2$, i.e. it is a local procedure. 

\medskip
\noindent
To our aim we need to determine the   $b$--elongations  with
excess $f$ which are ``maximal'', i.e the intervals  of maximum length which give  rise  to a positive or
negative 
$b$--elongations 
with excess $f$.

\noindent  {\bf  \Definition  (P41)    (The maximal $b$--elongations  with
excess $f$).} 
{\it Given $b>f$ positive real numbers, 
the  $\YY (\a)$, $\a \in \Z $,
 have  maximal $b$--elongations  with excess  $f$ if there exists an increasing sequence
$\{\a^*_i, i\in \Z\}$ such that
in each of the intervals
$[\a^*_{i},\a^*_{i+1}]$ we have either (1) or (2) below:  

(1) In the interval $ [  \a^*_{i}, \a^*_{i+1}]$ (negative maximal
elongation):
$$ \YY  ( \a^*_{i+1}) - \YY (\a^*_{i})  \le -b -f;\,  \quad
\YY  ( y) - \YY  ( x)    \le b -f,  
\quad \forall x<y \in  [  \a^*_{i}, \a^*_{i+1}];    
\Eq(3.3aa) $$
$$   \YY ( \a^*_{i+1}) \le  \YY ( \a)  \le  \YY (  \a^*_{i}),  
\quad  \a^*_{i}  \le  \a \le \a^*_{i+1}.  \Eq(3.3ab) $$

(2) In the interval $ [  \a^*_{i}, \a^*_{i+1}]$ (positive maximal elongation):
$$ \YY  ( \a^*_{i+1}) - \YY (\a^*_{i})  \ge b +f;\,
\quad   \YY  ( y) - \YY  ( x)    \ge -b +f,  
\quad \forall x<y \in  [  \a^*_{i}, \a^*_{i+1}];    
\Eq(3.3ba) $$
$$   \YY ( \a^*_{i})  \le  \YY ( \a)  \le  \YY (  \a^*_{i+1}),  
\quad  \a^*_{i} \le  \a \le \a^*_{i+1}.
\Eq(3.3bb) $$
Moreover, if in the  interval $  [  \a^*_{i}, \a^*_{i+1}]$ we have
\eqv(3.3aa) and \eqv(3.3ab) (resp. \eqv(3.3ba) and \eqv(3.3bb)) then in 
the adjacent interval $[\a^*_{i+1},\a^*_{i+2}]$ we have   \eqv(3.3ba)
and \eqv(3.3bb) (resp. \eqv(3.3aa) and \eqv(3.3ab)).
At last, we make the convention 
$$
\a^*_0\le 0 < \a^*_1.
\Eq(grr)
$$} 
\medskip
\noindent
\Remark (v2) In [\rcite{COPV}]  the convention $\a^*_{-1}\le 0< \a^*_0$
was assumed.

\medskip
\noindent
We say that the interval $ [  \a^*_{i}, \a^*_{i+1}] $ gives rise to 
negative  {\it maximal} $b$ elongations with excess  $f$ in the first
case and  the interval $ [  \a^*_{i}, \a^*_{i+1}] $ gives rise to  positive
{\it maximal} $b$ elongations with excess  $f$ in the second case.
\medskip
\noindent
\Remark (a1)  Note that if $\{\a^*_i, i\in \Z\}$  gives rise   to 
{\it maximal} $b$ elongations with excess  $f>0$, then $\{\a^*_i, i\in \Z\}$  gives rise  to 
{\it maximal} $b$ elongations with excess $f'$ with  $0 \le f'\le f $.  

\medskip
\noindent
The $\a^*_{i}$ are in fact $\a^*_{i}  \equiv \a^*_{i}  (\g,\e, b, f,
\o)$.
We will write  explicitly the dependence on one, some or all the
parameters only   when   needed.  Since we are considering  a random walk and $\a^*_i$ are points
of local extrema, see \eqv(3.3ab) and  \eqv(3.3bb), for a given
realization of the random walk, various sequences $\{\a^*_i, i\in
\Z\}$ could have the properties listed above.
This because  a random walk can
have locally and globally multiple maximizers or minimizers.
 Almost surely this  does not happen  for
the  Brownian motion.
In [\rcite {COPV}], 
we have chosen to take the first minimum time
or the first maximum time instead of the last one as in \eqv(NP1).
In the Brownian motion case the last and first maximum (resp. minimum) time
are almost surely equal. However we  could  have taken the last minimum
time or the last maximum time without any substantial change. From now on,  we make this last choice. 
With this choice and the convention \eqv(grr) the  points $\a^*_i$ are unambiguously
defined. 
The interval $[\a^*_0,\a^*_1]$ is called 
maximal $b$--elongation with excess $f$ that contains the origin.

\medskip

 \noindent \Remark (v3)   Obviously the construction of  maximal $b$--elongation with excess
$f$ cannot be  a local procedure.  So  to  determine $[\a^*_0,\a^*_1]$, for example, implies to know
that  the intervals
   adjacent to $[\a^*_0,\a^*_1]$  give risen to  $b$--elongation with excess $f$  (not necessarily
maximal) of  sign opposite to the one associated to $[\a^*_0,\a^*_1]$.

\medskip
\noindent
In [\rcite{COPV}] we
  determined   the maximal
$b$-- elongation with excess $f$  containing the origin and  
estimated the $\P$-probability
of the occurrence of $[\a^*_0,\a^*_{1}] \subset [-Q/\e,+Q/\e]$
taking  care of  ambiguities  
mentioned above.   
 Namely, applying  5.8, 5.9 and   Corollary 5.2   of  [\rcite{COPV}],  
choosing   $\d^*$,  $Q$ and $\e$ as in  
Subsection 2.5, $b=2\FF^*$,
and  see (5.30) in [\rcite{COPV}],
$f=5/g(\d^*/\g)$, we have proved
$$
\eqalign { & 
\P\left [ \left ( [\a^*_0,\a^*_{1}] \subset [-Q/\e,+Q/\e]
\right )^c\right]\le 
3 e^{-\frac{Q}{2C_1}}+ \e^{\frac{a}{16(2+a)}}
+Q^2\e^{\frac{a}{8+2a}}+ Q e^{-\frac{1}{2\e^{3/4}V^2(\b,\th)}} \cr &
\le \e^{\frac{a}{32 (2+a)}}=\Big(\frac{5}{g(\d^*/\g)}\Big)^{\frac{a}{8(2+a)}}.
}\Eq(Zorro.4b)
$$
where $C_1 \equiv C_1(\b,\th)$ is an explicit constant, $V(\b,\th)$ as in
\eqv(C.1a)  and $a>0$.
 Estimate \eqv(Zorro.4b) is obtained in [\rcite{COPV}] estimating the probability
to have enough
$b$--elongation with excess $f$ (not necessarily maximal) within $[-Q/\e,Q/\e]$
 to 
be sure that there exists a  maximal one   containing  the origin. 
Here we have a slightly different point of view, we want to be
able to construct {\it all} the maximal $b$--elongations with excess $f$
that are within $[-Q/\e,Q/\e]$. After a moment of reflection, one  
realizes that  the simultaneous occurrence of the events  that
two $b$--elongations with excess $f$ not necessarily maximal with 
opposite sign on the right of $[-Q/\e,Q/\e]$ and the same  on
its left should allow to  construct all the maximal
$b$--elongations with excess $f$
that are within $[-Q/\e,Q/\e]$. 
There is a simple device used constantly in  [\rcite{COPV}], to
estimate the $\P$--probability of the simultaneous occurrence
of such events
on $[Q/\e,(Q+L)/\e]$ and on  $[-(Q+L)/\e,-Q/\e]$ for some
$L>0$. Let us call these events  $\O^-_L(Q,f,b)$ and $\O^+_L(Q,f,b)$. 
Since it is rather long to introduce this device and it will be 
used  for other purposes, we postpone to the Subsection 5.5 the 
proof that choosing  the parameter as in Subsection 2.5, taking  
$L=cte \log (Q^2g(1/\g))$, one gets 
$$
\P[\O^-_L([Q,f,b)\cap \O^+_L(Q,f,b)] \ge 1-2\e^{\frac{a}{32(2+a)}}
\Eq(poml)
$$ 
for some $a>0$, see after \eqv(4.480B). 
Let us call 
$$
\O_L([-Q,+Q], f,b,0)\equiv \O^-_L([Q,f,b)\cap\left\{[\a^*_0,\a^*_{1}] \subset [-Q/\e,+Q/\e]\right\}
\cap \O^+_L(Q,f,b)
\Eq(oml1)
$$
where $0$ in  the argument of $\O_L( \cdot )$ is to  recall that $\YY(0)=0$.
The space  $\O_L([-Q,+Q], f,b,0) $ depends on the variables $\chi(\a)$ for $\e\a\in
[-Q-L,Q+L]$.  
Collecting \eqv(Zorro.4b) and \eqv(poml) one gets
$$
\P[\O_L([-Q,+Q],f,b,0)] \ge 1-3\e^{\frac{a}{32(2+a)}}.
\Eq(poml2)
$$ 

On $\O_L([-Q,Q],f,b,0)$ we have
$$
-\frac{Q}{\e}<\a^*_{\k^*(-Q)+1}\le \dots\le 
\a^*_{0}<0< \a^*_1 \le \dots \a^*_{\k^*(Q)-1}<\frac{Q}{\e}.
\Eq(palpha)
$$
  where $\k^*(\pm Q)$  are defined in \eqv (num). 
 The construction done in    [\rcite{COPV}],  just described  is a {\it bilateral} construction.
We considered the process $\YY(\cdot)$, $\YY(0)=0$, see \eqv (3.C6) and we determined to the right and
to the left of the origin the  $b-$ elongations with excess $f$.
The Neveu--Pitman construction, recalled in Subsection 5.2 is a {\it one side} construction.  The
determination of the points $S_i$  is achieved  moving the origin of the BM to  $-r_0$, and then
letting $ r_0 \to \infty$.      To be able to compare what we just  recalled   with the
Neveu--Pitman construction for the random walk we  translate   the origin to
$-r_0=-4pQ$, for some
$p\in \N$, and called
$\YY_{r_0}$ the new random walk with $\YY_{r_0}(r_0)=0$.
We want to  construct in the  interval   $[-Q/\e,+Q/\e]$ the    maximal
$b$--elongations with excess $f$ for  the process  $ (\YY_{r_0} (\a),
\a\in \Z)$   considering   as above,   extra elongations   on the left and on
the right of  $[-Q/\e,+Q/\e]$.  In this way we are able to compare in the same probability space
the      construction   done in   [\rcite{COPV}]  with the one by  Neveu-Pitman specialized for the
random walk   $\YY_{r_0}$.  

Repeating step by step the
construction of maximal $b$--elongations with excess $f$ given  in [\rcite{COPV}], 
and recalling \eqv(oml1), on a  subset  $\O_L([-Q,+Q],f,b,r_0)$ that depends only of the
variables $\chi(\a)$ for $\e\a\in [-Q-L,Q+L]$   (and in particular does not depends
of the variables $\chi(\a)$ for $\e\a \in [-Q-L-r_0, -Q-L-1]$), 
we can construct
all the maximal elongations that are within $[-Q,+Q]$ for the process
$ (\YY_{r_0} (\a),\a \in \Z)$.  
By translation invariance, using \eqv(poml2) we have
$$
\P[\O_L([-Q,+Q],f,b,r_0)]=
\P[\O_L([-Q,+Q],f,b,0)] \ge 1-3\e^{\frac{a}{32(2+a)}}.
\Eq(pomlrzero)
$$
 Similarly to  \eqv (palpha) we have on   $\O_L([-Q,+Q],f,b,r_0)$
$$
-\frac{Q}{\e}< \a^*_{\k^*(-Q,r_0)}(r_0)+1\le \dots\le  
\a^*_{0}(r_0)<0< \a^*_1(r_0) \le \dots   \a^*_{
\k^*(Q,r_0)-1}(r_0)<\frac{Q}{\e}
\Eq(tildealpha)
$$
where  
$$
\k^*(-Q, r_0)=\sup(i\ge1\,: \, \e  \a^*_i(r_0)<-Q)
\Eq(tildekappa01)
$$
and 
$$
\k^*(Q,r_0)=\inf (i\ge \k^*(-Q, r_0) :\, \e  \a^*_i(r_0)>Q)     
\Eq(tildekappa02)
$$
with the usual convention $\inf(\emptyset )=+\infty$. 

Since the previous construction depends only on the increments of
$\YY(\a)$ and  is exactly the one  used to construct 
$(\a^*_i ,i \,:-Q<\e\a_i^*<Q)$,
we have
$$
\eqalign{
&\left(  \a^*_i(r_0),\, \forall i\in \Z\,:\,  \k^*(-Q,r_0)<
i< \k^*(Q, r_0)\right)\,\,{\rm  on}\,\,  \O_L([-Q,+Q],f,b,r_0)\cr
&\mathrel{\mathop=^{\rm Law}}
\left( \a_i^*,\, \forall  i \in \Z\,:\, \k^*(-Q) < i <\k^*(Q) \right)\,\,{\rm
on}\,\, \O_L([-Q,Q],f,b,0).\cr
}\Eq(tildealpha2)
$$

Here $X \,\,{\rm on}\,\, \O_1 \mathrel{\mathop=^{\rm Law}}Y\,\,{\rm
on}\,\,\O_2$ means that the respective conditional distributions are the
same. 
 Note that  we have $  \a^*_0(r_0) <0<  \a^*_1(r_0)$ and
$\a^*_0<0<\a^*_1$.   
In particular
\eqv(tildealpha2) implies  that $ \a^*_1(r_0)$ on
$\O_L([-Q,+Q],f,b,r_0)$ and 
$\a^*_1$ on  $\O_L([-Q,Q],f,b,0)$ have the same law.

\medskip
\noindent{\bf 5.4.  Relation between $h$--extrema and maximal  $b$--elongation with  excess $f$ }
\medskip

Recalling \eqv(N.9), we have 
$${\YY}(\a)  =\sqrt { c(\b,\th, \g /\d^*)}  \hat W^\e  (\a \e),\quad
\forall \a\in \Z.  \Eq (E.2) $$
Furthermore   taking into account that    $(\hat \s_i(r_0), i\ge 1)$ are the times of $h$--extrema
for the random walk $V^\e_{r_0}$  starting a $-r_0=-4pQ$,
see    the end of Subsection 5.2,    
and the  properties (4.A) and (4.B) satisfied by   $(\hat \s_i(r_0), i\ge 1)$  one   recognizes immediately  
 that  the intervals   $[\hat \s_{i}(r_0),\hat \s_{i+1}(r_0))$
for $i\ge 1, i\in \N$   give rise to maximal 
$b=h\sqrt{c(\b,\th,\g/\d^*)}$ elongations with excess $f=0$, for any $b>0$. 
Let us define 
$$  \hat \k(-Q, r_0)=\sup\left(i \ge 1: \e\hat \s_i (r_0)<-Q\right). 
\Eq(kappa1)   
$$
We impose  $i\ge 1$ in \eqv(kappa1)  so that  
 $\hat\s_{\hat\k(-Q,r_0)}(r_0)$ is
a  time of a $h$--extremum. Recall that     
$\hat\s_0(r_0)$ may not be  a point of $h$--extrema.   
Furthermore  we   define  
$$
\hat \nu(r_0)=\inf\left(i\ge \hat \k(-Q,r_0): \e\hat\s_i(r_0)>0\right)
\Eq(nukappa)
$$
and 
$$
\hat \k(Q,r_0)=\inf\left(i\ge \hat \nu(r_0): \e\hat   \s_i
(r_0)>Q\right).
\Eq(kappa2)
$$ 
Note that on   $\{\hat \k(Q,r_0)<\infty\} $,
there are $\hat \k(Q,r_0)-\hat \k(-Q,r_0)+1$ points of $h$--extremum
within $[-Q,+Q]$.
So let 
$$
\O_0(Q,r_0)\equiv\left\{\o\in\O,\, \hat \k(-Q,r_0)<\nu(r_0)<\hat
\k(Q,r_0)<\infty,
\hat \k(Q,r_0)-\hat \k(-Q,r_0)\ge 1\right\}
\Eq (NN.2) $$
be the  set of realizations such  that there exists at least one
interval $ [\e \hat\s_{i}(r_0), \e \hat\s_{i+1}(r_0)] \subset
[-Q,Q]$, for some
$i\in \Z$ with $\hat\s_i(r_0)$ and  $ \hat\s_{i+1}(r_0)$ that are  $h$--extrema of
$V^\e_{r_0}(\cdot)$. 
On $ \O_0(Q, r_0)$   we have  
$$
-\frac Q \e <   \hat \s_{\hat\k(-Q,r_0)+1}(r_0) < ...< \hat \s_{\hat \nu(r_0)-1}(r_0)
<0<\hat \s_{\hat \nu(r_0)}(r_0) < ...<\hat \s_{\hat \k(Q,r_0)-1}(r_0) <\frac Q \e.
\Eq (E.5)
$$
Note that $\O_0(Q,r_0)\supset\O_L([-Q,+Q],b,f,r_0)$. Namely,  see Remark \eqv (a1),   if $[\e 
\a^*_i (f,r_0),\e
\a^*_{i+1} (f,r_0))$
gives rise to a maximal $b$--elongation with excess $f$, then
  it   gives rise
to a maximal $b$--elongation with excess $f=0$. Therefore
$\e \a^*_i (f,r_0)$ and $\e \a^*_{i+1} (f,r_0)$
are points of $h=b/\sqrt{c(\b,\th,\d^*/\g)}$ extrema.
Of course, it could exist a pair of  points of $h$--extrema, 
$h=b/\sqrt{c(\b,\th,\d^*/\g)}$, $\e\hat\s_{i}(r_0),\e\hat\s_{i+1}(r_0)$
for $\hat\k(-Q,r_0)\le i<\hat\k(Q,r_0)$ such that
$[\e\hat\s_{i}(r_0),\e\hat\s_{i+1}(r_0))$
gives rise to a maximal $b$--elongation with excess $f=0$ without
giving rise to a maximal $b$--elongation with excess $f>0$. 
That is, a priori on $\O_0(Q,r_0)\cap\O_L([-Q,+Q],b,f,r_0)=\O_L([-Q,+Q],b,f,r_0)$,
we have $ \hat\k(Q,r_0)-\hat\k(-Q,r_0)> \k^*(Q,r_0)- \k^*(-Q,r_0)$.

\medskip
\noindent {\bf \Lemma (46P)} {\it 
   Set
 $ b=2 \FF^*$,  $h= \frac { 2 \FF^*} {\sqrt{c(\b,\th,\d^*/\g)}}$, all the remaining parameters as in
Subsection 2.5, $L=cte\log (Q^2 g (\frac {\d^*} {\g} ))$ and $f= \frac 5 {g (\frac {\d^*} {\g}  )}$. Set  
$$ 
\O (f,r_0)=\O_L([-Q,+Q],b,f, r_0)\cap \left\{\hat\k(Q,r_0)-\hat\k(-Q,r_0) > \k^*(Q,r_0)-
\k^*(-Q,r_0)\right\}. 
\Eq(om1)
$$
 We have 
$$
\P [  \O (f,r_0)] \le
3 e^{-\frac{Q}{2C_1}}+ \e^{\frac{a}{16(2+a)}}
+Q^2\e^{\frac{a}{8+2a}}+ Q e^{-\frac{1}{2\e^{3/4}V^2(\b,\th)}}
\le \e^{\frac{a}{32 (2+a)}}.
\Eq(PPPP.1bis) $$
  where $C_1 \equiv C_1(\b,\th)$ is an explicit constant, $V(\b,\th)$ as in
\eqv(C.1a)  and $a>0$.}    

\proof 
Denote 
$$ \eqalign { \O'= & 
\Bigl \{ \o:\,
-\frac Q \e <   \hat \s_{\hat\k(-Q,r_0)+1}(r_0)<  \dots < \hat \s_{\hat \nu(r_0)-1}(r_0)
<0<\hat \s_{\nu(r_0)}(r_0) < \dots < \hat \s_{\hat \k(Q,r_0)-1}(r_0) <\frac Q \e
; \cr
&\qquad \exists i,\, \hat \k(-Q,r_0)+1\le i\le \hat\k(Q,r_0)-2\, \hbox {such that}\, 
[\hat \s_i(r_0),\hat\s_{i+1}(r_0))\, \hbox
{does not  satisfy (1)  and    (2) of}\cr
&\qquad  \hbox {Definition \eqv (P41) but does 
satisfy \eqv (3.5a) and \eqv (3.5b)  or \eqv (3.6a) and  \eqv
(3.6b) }
\Bigr\} } 
\Eq (Ebis.1) $$ 
Note that 
$$  
 \O(f,r_0) 
  \subset 
 \O' \cap\O_L([-Q,+Q],f,b,r_0).
 \Eq (E.9P)$$ 
To estimate the $\P$--probability of the event in the right hand side
of \eqv(E.9P), let $i$,  $\hat \k(-Q,r_0)+1 \le i\le \hat\k(Q,r_0)-2$ be such 
that $[\hat \s_i(r_0),\hat\s_{i+1}(r_0)]$  does not
satisfy (1) and   (2) of Definition
\eqv (P41) but does satisfy \eqv (3.5a) and  \eqv (3.5b) {or}
\eqv (3.6a) and  \eqv (3.6b).

It is enough
to consider the case where 
$[\hat\s_{i}(r_0),\hat\s_{i+1}(r_0)]$  does not  satisfy (1) of Definition
\eqv (P41) but does satisfy \eqv (3.5a) and  \eqv (3.5b). 
There are two cases: 
\noindent \item{$\bullet$} first case 
$$
-b -f \le \YY  ( \hat\s_{i+1}(r_0)) - \YY (\hat\s_{i} (r_0))  \le -b, \,
\quad   \YY  ( y) - \YY  ( x)    \le b -f  
\quad \forall x,y :\, x<y \in  [  \hat\s_{i}(r_0), \hat\s_{i+1}(r_0)]    
\Eq(S.1) $$
$$   \YY ( \hat\s_{i+1}(r_0))  < \YY ( \a)  \le  \YY (  \hat\s_{i}(r_0))  
\qquad  \forall \a :\,\hat\s_{i}(r_0) <  \a \le \hat\s_{i+1}(r_0)  \Eq(PS.2) $$

\noindent \item{$\bullet$} second case 

$$  \YY  ( \hat\s_{i+1}(r_0)) - \YY (\hat\s_{i} (r_0))  \le -b-f, \,
\exists x_0, y_0,\,  x_0<y_0 \in  [  \hat\s_{i}(r_0), \hat\s_{i+1}(r_0)]: \,
\, b \ge \YY  ( y_0)-\YY  ( x_0)  \ge b -f,    
\Eq(S.12)
$$
$$
\YY ( \hat\s_{i+1}(r_0))  < \YY ( \a)  \le  \YY (  \hat\s_{i}(r_0))  
\qquad  \hat\s_{i}(r_0) <  \a \le \hat\s_{i+1}(r_0).
\Eq(S.2a)
$$
Let us denote
 $$
\YY^*(\underline \a,  \a_1,  \a_2)
\equiv
\max_{ \a_1 \le \tilde \a\le  \a_2}
\sum_{\a=\underline \a}^ {\tilde \a} \chi(\a)
\Eq(4.52P)
$$
and
$$
\YY_*(\underline \a,  \a_1,  \a_2)
\equiv
\min_{ \a_1 \le \tilde \a\le  \a_2}
\sum_{\a=\underline \a}^ {\tilde \a} \chi(\a) \Eq(4.52PP)
$$
where  $\e\underline \a=-Q$.
To estimate both the cases 
we follow an argument already  used in the
proof of Theorem 5.1 in [\rcite{COPV}].
Take   $\rho'=(9f)^{1/(2+a)}$, for some $a>0$. Divide the interval $[-Q, Q]$  
into blocks of length $ \r'$  and  consider the event   
$$
\tilde\DD(Q,\rho',\e)  
\equiv \Big\{
\exists
\ell, \ell',\,-Q/\rho' \le \ell<\ell'\le (Q-1)/\rho' ;
|\YY^*(\underline \a, \sfrac{\rho' \ell}{\e}, \sfrac{\rho' (\ell+1)}{\e})
  -
 \YY_*(\underline \a,  \sfrac{\rho' \ell'}{\e},\sfrac{\rho'
(\ell'+1)}{\e})-b| \le 9f \Big\}.
$$
Simple observations show that those $\o$ that belong to 
$\{\max_{\a\in [-Q/\e,Q/\e]}|\chi(\a)|\le f\}$ and are such that
there exists $i$, $\hat\k(-Q, r_0)+1\le i\le \hat\k(Q,r_0)-2$    such that 
\eqv (S.1) and \eqv (PS.2) hold, belong also to  $\tilde\DD(Q,\rho',\e)$.

For the second case, we can assume that
$x_0$ is a local minimum and $y_0$ a local maximum, therefore
those $\o$ that belong to $\{\max_{\a\in [-Q/\e,Q/\e]}|\chi(\a)|\le
f\}$
and are such that there exists $i$, 
$\hat\k(-Q, r_0)+1\le i\le \hat\k(Q, r_0)-2$ such that 
\eqv (S.12) and \eqv (S.2a) hold, belong also to  $\tilde\DD(Q,\rho',\e)$.
Therefore we obtain that
$$  \O' \cap   \{ \max_{\a\in [-Q/\e,Q/\e]}|\chi(\a)|\le
f\}  \subset   \tilde \DD(Q,\rho',\e). $$
The estimate of    
 $  \P \left [\tilde\DD(Q,\rho',\e)  \cap  \O_L([-Q,+Q],f,b,r_0) \right]$  is done  
  in [\rcite{COPV}] where the   same set  $\tilde\DD(Q,\rho',\e) $, see pag 834 there,  was
considered.   It is based on Lemma 5.11 and Lemma 5.12 of [\rcite{COPV}].   
Here we recall the final estimate
$$ \eqalign { &   \P \left [\tilde\DD(Q,\rho',\e)  \cap  \O_L([-Q,+Q],f,b,r_0) \right]  \le  \cr & 
 8(2(Q+L)+1)^2\frac{2\sqrt{2\pi}}{V(\b,\th)} (9f)^{a/(2+a)}  +
(2(Q+L)+1)\frac{1296}{V(\b,\th)}  \frac{9f+(2+V(\b,\th))
\sqrt{\e\log \frac{C_1}{\e}}}{(9f)^{3/(4+2a)}} \cr
&\quad +\frac{4(Q+L)}{\e} e^{-\sfrac{f}{4\e
V^2(\b,\th)}}.} \Eq(3.1411)
$$
 Furthermore by  Chebyshev inequality  we obtain that 
$$
\P \left [ \{\max_{\a\in [-Q/\e,Q/\e]}|\chi(\a)|\ge
f\}\right ]  \le  \frac  {\E \left [ \{\max_{\a\in
[-Q/\e,Q/\e]}|\chi(\a)| \}\right ] } f  
\le 
2 \left (  \e V^2_+ \log \{ \frac {2 Q} \e \}  \right )^{\frac 12}  ( 1+  \frac 1 {  \log \{ \frac {2 Q}
\e
\}} )
$$
For the last inequality,  see  formula 5.38 in [\rcite{COPV}]. 
Choosing the parameters as in Subsection 2.5 we obtain the thesis.  
\eop

On $\tilde \O_L([-Q,+Q],f,b,r_0)=\O_L([-Q,+Q],f,b,r_0)\setminus \O (f,r_0)$,
\eqv(tildealpha) and \eqv(E.5) hold: a point is a beginning or an  ending   of an interval of 
maximal $b$--elongations  with excess $f$ if and only if it is a  point  of $h$--extremum. 
Relabel  the variables $\hat \s_i(r_0)$
in \eqv(E.5) as in Neveu and Pitman, that is  
define 
$$
\hat S_i(r_0)=\hat \s_{\nu(r_0)+i-1}(r_0),\,
\forall i\in \Z\,:\,\hat \k(-Q, r_0) \le \nu(r_0)+i-1<\hat \k(Q, r_0).
\Eq(mefi)
$$
Therefore, on $\tilde \O_L([-Q,+Q],f,b,r_0)$, we have
$$
\hat S_i(r_0)=  \a^*_i(r_0),\, \forall i\in \Z\, :\, -\frac{Q}{\e} \le
\hat S_i(r_0)\le \frac{Q}{\e}. 
\Eq(yemanja)
$$   

\medskip
\noindent {\bf \Lemma (45)} {\it 
  Take 
$$ b=2 \FF^*, \qquad  \qquad  h= \frac { 2 \FF^*} { V(\b,\th) },  $$ 
 all the remaining parameter as in Subsection  2.5,   
 $L=cte\log (Q^2 g (\frac {\d^*} {\g}))$ and $f= \frac 5 {g (\frac
{\d^*} {\g}) }$.

Let  $ \O_L([-Q,Q],f,b,0) $ be the probability space defined in
\eqv(oml1)
with $\P[ \O_L([-Q,Q],f,b,0)] \ge 1 -3\e^{\frac{a}{32(2+a)}}$   for
some $a>0$. Let 
$$ -\frac Q \e <    \a^*_{\k^*(-Q)+1}<....<  \a^*_{-1}
<\a^*_0  <0 < \a^*_1 <...<\a^*_{\k^*(Q)-1} <\frac Q \e
$$
be the maximal  $b$--elongations with excess $f$, see   \eqv (palpha), 
and 
$\{ S_i, i\in \Z \}$   the point process of $h$-extrema of the
BBM defined in Neveu-Pitman [\rcite{NP}]. 
We have  
$$ \lim_{\g \to 0} \e (\g) \a^*_i (\e(\g), f(\g))  \mathrel{\mathop=^{\rm Law}}  S_i  \qquad  i
\in \Z.  \Eq (D.2)   $$ 
} 

\medskip
\noindent \proof
This  is an immediate consequence of  \eqv(tildealpha2),
Lemma \eqv (46P), \eqv(yemanja), \eqv (2.21)  and the continuous mapping theorem.
\eop

\vskip0.5cm \noindent
{\bf Proof of Theorem \eqv (2)} The \eqv (law) is proved in Lemma \eqv(45).
The properties of $S_i$
are recalled in Subsection 5.2  and  \eqv (Laplace1) is proved in  [\rcite{NP}].     To
derive  \eqv (M.6)   let  $X= S_2-S_1$   be the interarrival times of
the renewal process $\{ S_i,  i\in \Z \}$.    
 Then using 
$$
\int_{0}^{\infty} \l e^{-\l z} \1_{\{x\ge z\}}\, dz=1-e^{-\l x },
$$ 
one gets
$$ \E [ e^{-\l S_1} ] = \frac 1 {h^2 \l} [1- \E [e^{-\l X} ] ]  = \frac 1 {h^2 \l}  [1-
\frac 1 {\cosh (h   \sqrt {2 \l})} ] 
\qquad {\rm for}\,\, \l \ge 0 .  $$ 
The distribution \eqv (6)  has been obtained in   [\rcite{Ke}],
 applying  the Mittag-Leffler representation for $(cosh z )^{-1}$.
Since $ \P[S_1>z] =\frac{1}{h^2}\int_{z}^{\infty} \P[X>x]dx$,
one obtains  differentiating  \eqv (6) the distribution in \eqv (D.1). \eop

 \medskip
\noindent 
{\bf Proof of Corollary \eqv (2bis)}
Since we already proved the convergence of finite dimensional
distributions see \eqv(law), to get \eqv(law2)
it is enough to prove that for any  
subsequence $\{u^*_\g, 0<\g<\g_0\}\in BV_{\rm
loc}(\R,\{m_\b,Tm_\b\})$, with $\g\downarrow 0$,  one can
extract a subsequence $\{u^*_{\g_n}, 0<\g_n<\g_0\}$ that convergences in Law.
In fact, since  $ BV_{\rm loc}(\R,\{m_\b,Tm_\b\}) $ is endowed with the topology induced
by the metric $d(\cdot,\cdot)$ defined in \eqv(Sk4), this  implies
that the points of jumps of  $\{u^*_{\g_n}, 0<\g_n<\g_0\}$ will converge
in Law to some points that by \eqv(law) are  necessarily the $(S_i, i\in
\Z)$, this will imply \eqv(law2). 

 So let $\g \downarrow 0$ be any subsequence that goes to $0$.
We will prove that for any chosen $\e_1$, 
it is possible to extract a  subsequence $\g_n\downarrow 0$ and 
to  construct  a probability subset
$\KK_\e \subset \O$ with 
$$
\P[\KK_{\e_1} ]\ge 1-\e_1
\Eq(compact21)
$$
so that  on  $\KK_{\e_1}  $, the subsequence 
$\{u^*_{\g_n},\,0<\g_n\le \g_0\}$ is a compact subset of
$BV_{\rm loc}(\R,\{m_\b,Tm_\b\})$.

To construct $\KK_{\e_1}$ and the subsequence $\g_n$, we use the following
probability estimates.   
Let $b=2\FF^*$
 and $\O_L([-Q,+Q],f,b,0)$ the probability subspace defined in   
\eqv(oml1),  $\P[\O_L([-Q,+Q],f,b,0)]\ge  1-3\e^{\frac{a}{32(2+a)}}$,
see \eqv (poml2).  On  $\O_L([-Q,+Q],f,b,0)$    
$u^*_\g(\cdot)$ jumps at the points $\{ \e\a_i^*, \k^*(-Q)+1 \le i\le
\k^*(Q)-1 \} $.  
It was proved in  Proposition 5.3  of [\rcite{COPV}]  that 
for   $i\in \Z$ and for   $0\le x\le (\FF^*)^2/(V^2(\b,\th)18\log 2)$  
$$
\P[\e\a^*_{i+1}-\e\a^*_{i}<x ]\le 2e^{-\frac{(\FF^*)^2}{18xV^2(\b,\th)}}.
\Eq(NewPPOP.1)
$$ 
By Lemma \eqv(36),  on the probability subspace $\O_{\rm urt}$, with 
$P[\O_{\rm urt}]\ge 
1-(\frac{5}{g(\d^*/\g)})^{\frac{a}{8(2+a)}}
$ for some $a>0$, the  
number of jumps within $[-Q,+Q]$ is smaller than 
$4+ 
\frac{8V^2_+}{(\FF^*)^2} Q\log\left[ Q^2 g(\d^*/\g)\right].$
Therefore, calling
$$
\O_Q(x,\g)
\equiv\left\{\o \in \O_{\rm urt};   \forall i:   \e\a^*_i \in [-Q,+Q], \e\a^*_{i+1}-\e\a^*_{i}>x\right\}
$$
one has 
$$
\P[\O_Q(x,\g)]\ge 1-
4(\frac{5}{g(\d^*/\g)})^{\frac{a}{32(2+a)}}
-\left (4+  \frac{8V^2_+}{(\FF^*)^2} Q\log\left[ Q^2 g(\d^*/\g)\right]\right)
 2e^{-\frac{(\FF^*)^2}{18xV^2(\b,\th)}}.
\Eq(Compact1)
$$

For any subsequence  $\g\downarrow 0$, one can pick up a subsequence $\{\g_n \}$    such that
$$
\sum_{n\ge 1}
\left(\frac{5}{g(\d^*(\g_n)/\g_n)}\right)^{\frac{a}{32(2+a)}}<\infty
\Eq(gammaenne)
$$
and recalling that $Q=Q(\g)\uparrow \infty$ when $\g\downarrow 0$, one can take $x=x(\g_n)>0$ such that
$$
\sum_{n\ge 1} \left (4+  \frac{8V^2_+}{(\FF^*)^2} Q(\g_n)\log\left[ Q^2(\g_n) g(\d^*(\g_n)/\g_n)\right]\right)
 2e^{-\frac{(\FF^*)^2}{18x(\g_n)V^2(\b,\th)}}<\infty.
\Eq(xgammaenne)
$$
Now using \eqv(Compact1), \eqv(gammaenne) and \eqv(xgammaenne), 
given $\e_1>0$, one can choose $n_0=n_0(\e_1)$ such that
$$
\P[\bigcap_{n\ge n_0}\O_{Q(\g_n)}(x(\g_n),\g_n)
]\ge 1-\e_1.
\Eq(compact21a)    
$$
Denote  $\KK_{\e_1}\equiv \bigcap_{n\ge n_0}\O_{Q(\g_n)}(x(\g_n),\g_n)
$ and we have proven \eqv (compact21).

Let $ \o \in \KK_\e $ and 
$\{ u^*_{\g_n}= u^*_{\g_n} (\o), n\ge n_0 \}$  the above constructed  subsequence.   
Sufficient  and necessary conditions  for the compactness of  $\{
u^*_{\g_n}, n\ge n_0 \}$
is   to exhibit for all $\tilde \e $ say, $\tilde \e<1/2$
and for some numerical constant $c$ a  {\sl finite}  $c \tilde\e$--net for $\{u^*_{\g_n},
n\ge n_0( \e)\}$, see [\rcite{Bill}] pg. 217. 
One can also assume that $n_0=n_0(\e,\tilde \e)$ is such that 
$$
e^{-Q(\frac{\d^*(\g_{n_0})}{\g_{n_0}})}\le \tilde \e
\Eq(condition)
$$
Set $y^2\equiv y^2_{\g_n} =\frac {\tilde \e x (\g_n)}{ 4(1+\tilde
\e)}$, {\bf let} $k_Q \in \Z $ and $k_{-Q}\in  \Z$ 
so that  $k_Q y_n^2\le Q<(k_Q+1)
y^2$ and respectively  $k_{-Q} y^2 \le -Q<(k_{-Q}+1)y^2$.  Denote 
 $\BB(y^2, Q) \subset  BV_{\rm loc}$   the finite subset 
 $$\BB(y^2, Q) = \left \{ \eqalign { & u_0 \in BV_{\rm loc}: u_0
\quad \hbox {constant on} \quad  [ky^2, (k+1)y^2),
k\in [k_{-Q}, k_Q] \cap   \Z,   \cr &  \forall r \ge Q,   u_0(r)= 
u_0(k_Q); \quad  \forall r\le -Q,   u_0(r)=u_0(k_{-Q}) } \right  \} $$
 Let  $ \o \in \KK_\e$  and $k_i^*\equiv  k_i^*(\o, \g_n) \in \Z $
such that $k_i^*y^2\le \e(\g_n)\a^*_i(\o, \g_n)<(k^*_i+1)y^2$, for
all $i$ such that $\e\a^*_{i-1} \in [-Q,+Q]$.    Let $ u_0 \in
\BB(y^2,Q)$ such that $u_0(k_i^*y^2)=u^*_{\g_n}(\e\a^*_i)$.
It remains to  check that  $d(u^*_{\g_n}, u_0)\le c\tilde \e$
for some numerical constant $c$, where $d(\cdot, \cdot)$ is defined in \eqv(Sk4).
Let us define the following $\l_{\g_n}(.) \in \L_{\rm Lip}$ by
$\l_{\g_n}(k_i^*y^2)=\e\a_i^*$ and linear between $k_i^*y^2$ and $(k_{i}^*+1) y^2$ for 
$r>Q$ take $\l_{\g_n}(r)=\l_{\g_n}(Q)+t-Q$ and for $r\le -Q$ take
$\l_{\g_n}(r)=\l_{\g_n}(-Q)+t+Q$.
For all $i$ such that  $\e\a^*_{i-1} \in [-Q,+Q]$,
one has 
$$
\left|\l_{\g_n}(k_i^*y^2)-\l_{\g_n}(k_{i-1}^*)-(k_i^*-k_{i-1}^* )y^2
\right|=
\left|\e\a^*_{\ell+1}-\e\a^*_{\ell}-(k_i^*-k_{i-1}^*) y^2 \right|\le 2y^2.
\Eq(compact2)
$$
On the other hand on $\KK_\e$ one has
$\e\a^*_i-\e\a^*_{i-1}\ge x(\g_n)$ and therefore $(k_i^*-k_{i-1}^*)y^2>x(\g_n)-2y^2$.
Using $2y^2\le \tilde \e(x(\g_n)-2y^2)$ and \eqv(compact2), one gets
$$
\left|\l_{\g_n}(k_i^*y^2)-\l_{\g_n}(k_{i-1}^*)-(k_i^* -k_{i-1}^*) y^2 \right|
\le2y^2\le \tilde \e(x(\g_n)-2y^2)\le \tilde \e (k_i^*y^2-k_{i-1}^*y^2).
\Eq(compact4)
$$
Since $\l$ is piecewise linear one has also, for $s<t \in [k^*_{i-1}y^2, k^*_iy^2)$  
$$
\left|\l_{\g_n}(t)-\l_{\g_n}(s)-(t-s)\right|\le \tilde \e (t-s).
\Eq(compact5)
$$
Since $\l_{\g_n}$ has a slope 1 outside $[-Q,+Q]$, one gets for all 
$s<t \in \R$
$$
\log(1-\tilde \e)\le \log \frac{\l_{\g_n}(t)-\l_{\g_n}(s)}{t-s}\le
\log(1+\tilde \e).
\Eq(compact6)
$$
Therefore, recalling \eqv(Sk1), \eqv(compact6) entails
$\|\l_{\g_n}\|\le 4 \frac {\tilde \e} 3$ and 
using \eqv(condition) to control $\int_Q^\infty e^{-T}\,dT$ in \eqv(Sk4)
, one gets after an easy computation 
$d(u^*_{\g_n}, u_0)\le 3 \tilde \e$.  \eop

\medskip
\noindent{\bf 5.5. Probability estimates}
\medskip

We recall the already mentioned  device constantly used in
[\rcite{COPV}].     Lemma \eqv(35), stated  below,  gives   
 lower and upper bound on  the $\a^*_i$, $ i \in \Z$,  in term of
suitable
stopping times.
We set
$\hat T_0=0$, and define, for $k\ge 1$:
$$
\eqalign{
\hat T_k&=\inf \{ t>\hat T_{k-1}\colon
|\sum_{\a=\hat T_{k-1}+1}^{t} \chi(\a)| \ge \FF^*+\frac f 2 \},\cr
\hat T_{-k}&=
 \sup \{ t<\hat T_{-(k-1)}\colon
|\sum_{\a=t+1}^{\hat T_{-(k-1)}}\chi(\a)| \ge \FF^*+\frac f 2 \}.}
\Eq(4.20B)
$$
Clearly, the random variables $\D\hat T_{k+1}:= \hat T_{k+1}-\hat T_{k}$, $k\in
\Z$, are
independent and identically distributed.
(Note  that, by definition,  $\D\hat T_1=\hat T_1$.)

\remark Note that $(\hat T_i, i\in \Z)$ was called
$(\tau_i, i \in \Z)$ in [\rcite{COPV}], we change their names to avoid
ambiguities with the
$\t$ defined in \eqv(NP1) and the ones defined after \eqv(la). 

We define,
 $$ \tilde  S_k= {\rm sgn} \Big (
\sum_{j=\hat T_{k-1}+1}^{\hat T_{k}}
\chi(j)\Big); \qquad
 \tilde  S_{-k}= {\rm sgn} \Big ( \sum_{j=\hat T_{-k }+1}^{\hat T_{-k+1}}
\chi(j)\Big)\qquad \text{for  } k \ge 1.
\Eq(4.23B)
$$
The following lemma  estimates the probability to  detect
at least   one $b= 2 \FF^*$  elongation,
with excess $f$, not necessarily maximal.   The proof
is done in  [\rcite{COPV}], Lemma 5.9 there.

\smallskip
\noindent{\bf \Lemma (lem6)}
{\it 
There exists an $\e_0$ such that 
for  all $0<\e <\e_0$,
all integer $k\ge 1$, and all $s>0$ we have
$$
\P\left[\hat T_k \le \frac { k (s +\log 2)C_1 }{\e}
; \exists i\in\{1,\dots,k-1\}, \tilde  S_i=\tilde  S_{i+1}\right]
\ge \left(1- e^{- sk}\right)(1-\frac 1{2^{k-1}}).
\Eq(4.43B)
$$
for some  $C_1=C_1(\b,\th)$.}
\smallskip
\noindent
To  detect  elongations with alternating sign,  we introduce 
on  the right of the origin
$$
\eqalign{ i^*_1&\equiv \inf \left\{i\ge 1: \tilde  S_i=\tilde  S_{i+1}\right\}\cr
i^*_{j+1}&\equiv\inf\left\{i \ge (i^*_j+2): \tilde  S_i=\tilde
S_{i+1}=-\tilde  S_{i^*_j}
\right\} \quad \qquad  j \ge 1,}
\Eq(4.46B)
$$
and on  the left
$$
i^*_{-1} \equiv
\cases{
-1 & if
$\tilde  S_{-1}=\tilde  S_1=-\tilde  S_{i^*_1}$,\cr 
 \sup \left\{i\le -2: \tilde  S_i=\tilde  S_{i+1}=-\tilde S_{i^*_1}\right\} & if
$\tilde  S_{-1}\neq \tilde S_1$ or $\tilde S_1=-\tilde S_{i^*_1}$,\cr }
$$
$$
i^*_{-j-1}\equiv\sup\left\{i \le i^*_j-2: \tilde S_i=\tilde S_{i+1}=-\tilde S_{i^*_j} 
 \right\} \quad  \qquad  j \ge 1.
\Eq(4.46'B)
$$
The corresponding estimates are given by the following Lemma which  was
proved in [\rcite{COPV}], see Lemma 5.9 there.

\noindent{\bf \Lemma(lem7)} {\it 
There exists an $\e_0$ such that for  all $0<\e <\e_0$,
all $k$ and  $L$ positive integers,  $L$ even,
(just for simplicity of writing) and all $s>0$ we have: 
$$
\P\left[ \hat T_{kL-1} \le \frac{ (kL-1)(s
+\log2)C_1}
{\e },\forall_{1\le j \le k}\,\, i^*_j < jL\right]
\ge \left(1-e^{-s (kL-1)}\right)\left(1-\sfrac{1}{2^{L-1}}\right)
\left( 1-\left(\sfrac34 \right)^{L/2}\right)^{k-1}    
\Eq(4.48B)
$$
  and 
$$
\eqalign{
&
\P\left[ \hat T_{-kL} \ge \frac{ -kL(s
+\log2)C_1}
{\e},\, 
\hat T_{L-1} \le \frac{ (L-1)(s
+\log2)C_1}
{\e}, \, i^*_1 <L, \,
\forall_{1\le j \le k}\,\, i^*_{-j} > -jL\right] 
\cr&\quad 
\ge \left(1-e^{-s (kL-1)}\right)\left(1-\sfrac{1}{2^{L-1}}\right)
\left( 1-\left(\sfrac34 \right)^{L/2}\right)^{k}.
    \cr
}\Eq(4.480B)
$$
where $C_1=C_1(\b,\th) $ is a constant.}
\vskip0.5cm \noindent 
Applying  Lemma \eqv(lem7) with $L=cte \log(Q^2g(\frac {\d^*} \g ))$,
taking  the parameters as in 
Subsection 2.5,    one gets
\eqv(poml) by  a  short computation.  
The basic fact that was used constantly in [\rcite{COPV}] even if it
was not formulated in its whole generality is the following.  

\noindent{\bf \Lemma(35)} { \it  On $\O_L([-Q,+Q],f,b,0)$, see \eqv(oml1),
we have 
$$
\hat T_i\le \a^*_{i+1}, \, \forall i: 1\le i< \k^*(Q)
\Eq(PPP.10)
$$
and
$$
 \a^*_{i}\le \hat T_{i^*_{i+1}},\, \forall i: 1\le i<\k^*(Q),
\Eq(PPP.11)
$$
where $\k^*(Q)$ is defined in  \eqv(num). 
}

\proof 
Recall that on $\O_L([-Q,+Q],f,b,0)$ 
we have assumed that $\a^*_0\le 0< \a^*_1$. 
To prove \eqv(PPP.10) we  start  proving 
 that $\hat T_1\le \a^*_2$. 
Suppose that 
$\a_2^* < \hat T_1$.  
Then,  from \eqv (4.20B), since $\a^*_1<\a^*_2<\hat T_1$ we have
$$
 |\YY(\a^*_1)|<\FF^*+\frac{f}{2} \qquad \hbox { and } \qquad 
|\YY(\a^*_2)|< \FF^*+\frac{f}{2}
 \Eq(APP.1)
$$
which is a contradiction since by assumption  $[\e\a^*_1, \e\a^*_2]$ is a maximal
$2\FF^*$ elongation with excess $f$, see Definition
\eqv (P41). Similar arguments apply  for $ i \ge 2$.     
 Now we prove \eqv(PPP.11). 
Assume that $[\a^*_0,\a^*_1]$ gives rise
to a positive elongation. The case of a negative elongation is  similar. 
Let us check that  $\a^*_1\le \hat T_{i^*_2}$.
By definition of ${i^*_1}, {i^*_2}$  we have 
that $[\hat T_{i^*_1-1},\hat T_{i^*_1+1}]$ is within an elongation
with a sign, say $\hat S_{i^*_1}$
and $[\hat T_{i^*_2-1},\hat T_{i^*_2+1}]$ is within an elongation with 
opposite sign, $\hat S_{i^*_2}=-\hat S_{i^*_1}$. Therefore, either $\hat S_{i^*_1}$ or
$\hat S_{i^*_2}$ is negative, which implies that $\a^*_1\le \hat T_{i^*_2}$.
The general case is clearly the same. \eop

Given an integer $R>0$, we denote as in  \eqv(num) 
$
\k^*(R)=\inf\{i\ge 1:\e\a^*_i\ge R\}.
$
We   define the stopping time  $\tilde k(R)=\inf\{i\ge 0: 
\e\hat T_i\ge R\}$. By definition
$$
\e\hat T_{\tilde k(R)-1}<R\le \e\hat T_{\tilde k(R)}
\Eq(APP.3)
$$
Using \eqv(PPP.10), we get that
$$
R\le \e\hat T_{\tilde k(R)} \le \e\a^*_{\tilde k(R)+1}
\Eq(APP.4)
$$
therefore 
$$
\k^*(R)\le \tilde k(R)+1.
\Eq(APP.5)
$$  
\vskip0.5cm 
\noindent{\bf \Lemma(36)}
{\it  There exists    $\O_{urt}$,
$\P [ \O_{urt}] \ge 1-(\frac{5}{g(\d^*/\g)})^{\frac{a}{8(2+a)}}
$, where $a>0$,
 so that  for all $R>1$
$$
\k^*(R) \le 1+\tilde  k(R)\le
2+ 
\frac{4V^2_+}{(\FF^*)^2} R\log\left[ R^2 g(\d^*/\g)\right]
\Eq(RAP.1)
$$
and
$$
\e\a^*_{\k^*(R)+1} \le  
 \frac{24C_1 V^2_+\log 2}{(\FF^*)^2\log (4/3)} R  
\left[\log(R^2g(\d^*/\g))\right]^2, 
\Eq(RAP.2)
$$
  where $ V_+= 
  V(\b,\th) 
\left[ 1+ ( \g/\d^*)^{\frac 1 5 }\right]$  and  $C_1=C_1(\b,\th)$ is a positive constant }. 
\medskip
\noindent
\remark It is well known that, almost surely,   
$\lim_{R\uparrow \infty} \tilde k(R)/R=(\E[\hat T_1])^{-1}$, 
see [\rcite{As}]  Proposition 4.1.4. The estimate \eqv(RAP.1)
allows us to have an upper  bound  valid uniformly with respect to $R\ge 1$ with an 
explicit bound on the probability.  This is the main reason to have
a $\log [R^2 g(d^*/\g)]$ in the right hand side of \eqv(RAP.1). 
\medskip
\noindent
\proof
We can assume that we are on  $\O_L([-Q,+Q],f,b,0)$.
Suppose first that $\tilde k(R)>1$.  Since \eqv(APP.3), we get
$$
\frac{ \e\hat T_{\tilde k(R)-1}}{\tilde k(R)-1}
<\frac{R}{\tilde k(R)-1}\le \frac{\e\hat T_{\tilde k(R)}}
{\tilde k(R)-1}.
\Eq(APP.6)
$$
Applying   Lemma 5.7 of  [\rcite{COPV}],   setting  $ b=\FF^*+(f/2)$ and $ V_+= 
  V(\b,\th) 
\left[ 1+ ( \g/\d^*)^{\frac 1 5 }\right]$,   we obtain 
that for all
$s$ with  
$0<s<(\FF^*+(f/2))^2[4(\log 2)V_+^2]^{-1}$, for all positive integer $n$
$$
\P\left[\e\hat T_n\le {ns}\right]\le e^{-n\frac{(\FF^*)^2}{4sV^2_+}}.
\Eq(APP.7)
$$ 
Therefore 
$$
\P\left[\exists n\ge 1: \frac{\e\hat T_n}{n}\le {s}\right]\le 
\frac{e^{-\frac{(\FF^*)^2}{4sV^2_+}}}{1-e^{-\frac{(\FF^*)^2}{4sV^2_+}}}.
\Eq(APP.8)
$$
Applying \eqv(APP.6), we get that for $\tilde k(R)>1$
$$
\P\left[\tilde k(R) \le 1+\frac{R}{s}\right]\ge 
\frac{1-2e^{-\frac{(\FF^*)^2}{4sV^2_+}}}{1-e^{-\frac{(\FF^*)^2}{4sV^2_+}}}.
\Eq(APP.9)
$$
When   $ \tilde k(R) =0$ or  $ \tilde k(R) =1$,  \eqv(APP.9) is certainly true,
therefore  \eqv(APP.9) holds  for all $\tilde k(R)\ge 0$.
Choosing in \eqv (APP.9) 
$$
s_0^{-1}=\frac{4V^2_+}{(\FF^*)^2} 
\left[ \log R^2 g(\d^*/\g)\right]
\Eq(APPP.10)
$$
we get
$$
\P\left[\forall{R\ge 1},\, \tilde k(R) \le 1+\frac{R}{s_0}\right]\ge 
1-\sum_{R\ge 1} \frac{\frac{ 2}{g(\d^*/\g)R^{2}}}{1- 
\frac{ 2}{g(\d^*/\g)R^{2}}}
\ge 1-\frac{3}{g(\d^*/\g)}. 
\Eq(APPP.11)
$$
Recalling  \eqv (APP.5),   for all $R\ge 1$, 
$$
\k^*(R) \le 1+\tilde  k(R)\le 2+ 
\frac{4V^2_+}{(\FF^*)^2} R\left[ \log R^2 g(\d^*/\g)\right]
\Eq(APPP.110)
$$
which is \eqv(RAP.1).  Next we  prove \eqv(RAP.2). 
Applying  \eqv(PPP.11)  and \eqv(APP.5)   we have
$$
\e\a^*_{k(L)+1} \le \e \hat T_{i^*_{\tilde k(L)+2}}.
\Eq(APP.11)
$$
Using \eqv(4.48B) with
$$
\eqalign{
L&=  L_0 = 1+3 \frac{\log (R^2g(\d^*/\g) )}{\log(4/3)}\cr
k&=  k_0= 2+\frac{4V^2_+}{(\FF^*)^2} R\log[R^2g(\d^*/\g)].\cr
}\Eq(APPP.13)
$$
After an easy computation, given $R\ge 1$
with a $\P$--probability greater than 
$
1-c(\b,\th)\frac{\log (R^2g(\d^*/\g))}
{g(\d^*/\g)^{3/2}R^{2}}
$
we have
$$
\e\hat T_{(2+k_0)L_0}
\le 
 \frac{24C_1 V^2_+\log 2}{(\FF^*)^2\log (4/3)} R 
 \left[\log(R^2g(\d^*/\g))\right]^2,\quad
\forall j:\, 1\le j \le 
k_0,\,\, i^*_j < j L_0.
\Eq(APPP.14)
$$
Therefore, with a $\P$--probability greater than 
$$
1-c(\b,\th) \frac{\log g(\d^*/\g)}{g(\d^*/\g)^{3/2}} \ge 1- \frac 1  {g(\d^*/\g)}  
\Eq(APPP.15)
$$
for all $R\ge 1$, \eqv(APPP.14) holds. Using \eqv(APPP.110)
we have, for all $R\ge 1$, 
$$
1+\k^*(R)< \tilde k(R)+2\le 3+
\frac{4V^2_+}{(\FF^*)^2} R\left[ \log R^2 g(\d^*/\g)\right].
\Eq(APPP.16)
$$
Therefore collecting \eqv(APPP.14) and \eqv(APPP.16)
we obtain that for all $R\ge 1$
$$
i^*_{\tilde k(R)+2}\le (2+k_0)L_0.
\Eq(APPP.17)
$$
 From which using again \eqv(APPP.14) and recalling \eqv(APP.11),
 we get that  for all $R\ge 1$
 
$$
\eqalign{
\e\a^*_{\k^*(R)+1} & \le \e\hat T_{i^*_{\tilde k(R)+2}} 
  \le \e\hat T_{(2+k_0)L-0}\cr
&\le 
 \frac{24C_1 V^2_+\log 2}{(\FF^*)^2\log (4/3)} R  \left[\log(R^2g(\d^*/\g))\right]^2\cr
}\Eq(APPP.18)
$$
which is \eqv(RAP.2).    Denote  by $\O_{urt} $ the intersection of        
$\O_L([-Q,+Q],f,b,0)  $  with the probability subsets in  
\eqv (APPP.11) and
\eqv (APPP.15).  Recalling \eqv (poml) and the choice of $\e$, see \eqv(epsilon),
we get  the Lemma. 
\eop

\vskip0.5cm \noindent  {\bf {Proof of Theorem \eqv (copv):} }  

We  need to estimate the Gibbs probability of the set 
$\PP^\r_{\d,\g,\z,[-Q,+Q]}(u^*_\g(\o))$, see \eqv (copvprob). 
According to the  definition  \eqv (D.30) we need to prove that  
  on $\O_1$   the minimal distance between two
points of jump of $u^*_\g$ is larger than $8\r+8\d$.
 Define 
$$
\O_{1,1}=\left\{\o \in \O_{urt} : \forall  i,  \,\, -Q\le \e\a^*_i\le Q;
\,   \e\a^*_{i+1}-\e\a^*_{i}\ge 8\r+8\d \right\}.
\Eq(omega4)
$$
where $\O_{urt}$ is the probability subspace that occurs in Lemma
\eqv(36). 
On $\O_{urt}$, see Lemma
\eqv(36),   the total number of jumps of $u^*_\g$ within
$[-Q,+Q]$ is bounded by $2K(Q)+1$ with $K(Q)$ given in \eqv(M.1).  
Since the points of jumps of $u^*_\g$ are the $\e\a^*_i,
i\in \Z$,  from Proposition 5.3  in [\rcite{COPV}]  we have that 
for all $i\in \Z$, for all $0\le x\le (\FF^*)^2/(V^2(\b,\th)18\log 2)$  
$$
\P[\e\a^*_{i+1}-\e\a^*_{i}<x ]\le 2e^{-\frac{(\FF^*)^2}{18xV^2(\b,\th)}}.
\Eq(NewPPO.1)
$$
Then one gets  
$$
\P\left[\O_{1,1}\right]
\geq 1- 
\left(\frac{5}{g(\d^*/\g)}\right)^{\frac{a}{8(2+a)}}
-2K(Q)
e^{-\frac{(\FF^*)^2}{18(3\r+3\d)V^2(\b,\th)}}.
\Eq(NewPPO.2)
$$
Recalling \eqv(rho), \eqv (queue) and \eqv(5.97037) 
 one gets  
$$
\P\left[ \O_{1,1}\right]
\geq 1-\left(\frac{5}{g(\d^*/\g)}\right)^{\frac{a}{10(2+a)}}.
\Eq(NewPPO.3)
$$
Denote by
$$ \O_1= \O_{\g, K(Q)} \cap \O_{1,1} \Eq (omega)$$
where  $\O_{\g, K(Q)}$ is the probability subspace in 
 Theorem   2.2   of [\rcite{COPV}] and  $ K(Q)$ is given in \eqv
(M.1). From
the results stated in Theorem 2.1, 2.2 and 2.4 of [\rcite{COPV}] 
we obtain  \eqv (5.02)  and \eqv (copvprob).

\vskip 1truecm

\chap{6 Proof of Theorem \eqv (min) and \eqv (main1)}6
\numsec= 6
\numfor= 1
\numtheo=1
\medskip
 
\noindent{\bf 6.1. Proof of Theorem \eqv (min) }
\medskip
 Let $\{ W(r), r \in
\R \}$  be  a realization of  the Bilateral  Brownian motion.
Let $u^* (r)\equiv u^* (r, W)$, $r
\in
\R$, be  the function defined in \eqv (SS.1a) and  \eqv (SS.1).  
 As 
    consequence  that
$ \PP$ a.s the number of renewals in any finite interval is finite,
we have that $\PP$ a.s  $u^* \in BV_{{\rm loc} } $ .
To prove the theorem we  need to show
that  for $v  \in BV_{{\rm loc}}$,  $ \PP$ a.s. $  \G (v|u^*, W) \ge 0$, the equality holding 
only  when $v= u^*$.  Let $ S_{0}$ be a point
of $ h-$ minimum, $h=\frac   {2 \FF^*} {V(\b,\th)} $.     This,  by
definition, implies that  in the interval $ [  S_{0},S_{1})$ 
$$ W ( S_{1}) - W (S_{0})  \ge  \frac   {2 \FF^*} {V(\b,\th)},
\,  \qquad   W( y) - W ( x)  >  - \frac   {2 \FF^*} {V(\b,\th)},    
\qquad \forall x<y \in  [  S_{0}, S_{1})    
\Eq(3.3v) $$
$$   W( S_{0})  \le W ( x)   \le   W (  S_{1})  
\qquad  S_{0} \le  x \le  S_{1} .  \Eq(3.3b) $$
Suppose first that  $v$    differs from
$u^* (W)$ only in   intervals contained   in $[S_0,S_1)$. 
Since $u^*(r)= m_\b$,    for $r \in [S_0,S_1)$, see \eqv (SS.1a), set $v(r)= Tm_\b
\1_{[r_1,r_2)}$ for
$[r_1,r_2)\subset [S_0,S_1)$ and $v(r)= u^*(r)$  for $r \notin [r_1,r_2)$.
When  the interval $[r_1,r_2)$ is strictly contained in $ [S_0,S_1)$ the function $v$
has two jumps more than $u^*$.  
Then the  value of  $\G (v|u^*, W)$,  see \eqv
(D.5),   is 
$$   \G (v|u^*, W)=  \G_{[S_0, S_1)} (v|u^*, W)=  2 \FF^* + V(\b,\th) [W(r_2)-W(r_1)] > 0,  \Eq (7.1) $$
which is strictly positive  using   the second property in \eqv (3.3v).
 When $[r_1,r_2) \equiv [S_0,S_1)$ then the function $v$
has two jumps less than $u^*$.  Namely $u^*$  jumps   in $S_0$ and   in $S_1$ and $u$ does not.
The value of $  \G (v|u^*, W)$  in such   case is   
$$   \G (v|u^*, W)=   \G_{[S_0, S_1)} (v|u^*, W)+  \G_{[S_1, S_2)} (v|u^*, W)= 
 -  2 \FF^* + V(\b,\th) [W(S_1)-W(S_0)] \ge 0.  \Eq (7.2) $$
The last inequality holds since the first property in \eqv (3.3v).
 In the case in which $[r_1,r_2) \subset [S_0,S_1)$, $r_1=S_0$, $r_2 < S_1$ then the function
$v$ has the same number of jumps as   $u^*$. 
The value  of $  \G (v|u^*, W)$     
is   
$$   \G (v|u^*, W)=  \G_{[S_0, S_1)} (v|u^*, W)=  V(\b,\th) [W(r_2)-W(S_0)] \ge 0  \Eq (Ma.6)  $$
which is still positive  because of    \eqv (3.3b).
When  $[r_1,r_2)  \subset   [S_0,S_1)$, $r_1>S_0$, $r_2 = S_1$ then, as in  the
previous case,   the function
$v$ has the same number of jumps as   $u^*$ and again by     \eqv (3.3b), 
$$   \G (v|u^*, W)=  \G_{[S_0, S_1)} (v|u^*, W)+  \G_{[S_1, S_2)}
(v|u^*, W) =   V(\b,\th) [W(S_1)-W(r_1)] \ge 0.  \Eq (Ma.5) $$
 The case when  $v$ differs from  $u^*$, still only  in $[S_0,S_1)$,
but  in more than one interval     can
be   reduced to the previous cases.  For simplicity, suppose that 
 $v(r)= Tm_\b
\1_{[r_1,r_2) \cup [r_3,r_4) }$ for
$[r_1,r_2) $ and  $[r_3,r_4) $ both subsets of $[S_0,S_1)$ and  $v(r)= u^*(r)$  for $r
\notin [r_1,r_2) \cup [r_3,r_4) $.
We have that
 $$  \G (v|u^*, W)=
  \G (v_1|u^*, W) +  \G (v_2|u^*, W) \Eq (Ma.1) $$
where $v_1(r) = Tm_\b
\1_{[r_1,r_2)} + u^*
\1_{[r_1,r_2)^c}$ and $v_2(r) = Tm_\b
\1_{[r_3,r_4)} + u^*
\1_{[r_3,r_4)^c}$. Equality \eqv (Ma.1) is an obvious consequence of the linearity of the
integral and    the observation that $  \G (u^*|u^*, W)=0$.    Each  term in \eqv (Ma.1) can
  then be  treated as in  the previous cases.  By assumption
$v \in    BV_{{\rm loc} } $  and
    then the number of intervals   in $[S_i, S_{i+1}) $  where  $v$ might differ from $u^*$ is 
$\PP$ a.s  finite.   The conclusion is therefore that  if
$v
\neq u^*$ in
$ [S_0,S_1)$ 
$$  \G (v|u^*, W) \ge 0.  \Eq (Ma.8)$$
  When  
$v$ differs from $u^*$   in  $[S_1,S_2)$ and  $S_1$  is an $h-$ maximum one  repeats the
previous arguments   recalling  that by definition  in   $ [ S_{1},S_{2})$ 
$$ W( S_{2}) -W (S_{1})  \le -h \,  \qquad   W  ( y) - W ( x)    \le  h   
\qquad \forall x<y \in  [  S_{1}, S_{2})    
\Eq(7.3) $$
$$  W ( S_{2})  \le  W ( x) <  W (  S_{1})  
\qquad  S_{1}  \le  x < S_{2}  \Eq(7.4) $$
and  $u^* (r)= Tm_\b $, for $r \in   [S_1,S_2)$, see \eqv (SS.1). 
Then
repeating step by step  the previous scheme 
one concludes that  $\PP$ a.s. 
$$  \G (v|u^*, W) \ge 0. $$
In the general case   
 $$  \G (v|u^*, W)= \sum_{i \in \Z } \G_{ [S_i, S_{i+1})}(v|u^*, W)  \ge 0 . \Eq (MaM1)$$
 To  prove that  $u^*$ is  $ \P$ a.s. the unique minimizer of  $ \G (\cdot| u^*, W)$ it is enough    
to  show that each term     among  \eqv (7.2),    \eqv (Ma.6) and
\eqv (Ma.5)   is   strictly positive, so that  we    get     a strict
inequality in  \eqv (MaM1).  Since, see  [\rcite {RY}], page 108,  exercise (3.26),
$$\PP [  \exists r \in [S_0,S_1] :    [W(r)-W(S_1) ]=0] =0 ,$$ 
we obtain that \eqv (Ma.5),  \eqv (7.2) and by a simple argument
\eqv (Ma.6) are  strictly positive.  \eop
\medskip
 
\noindent{\bf 6.2. Proof of Theorem \eqv (main1) }  The proof of \eqv
(main)  is an immediate consequence of 
Proposition \eqv (31) and  Theorem \eqv(2).   \eop
\medskip
 
\centerline{\bf References}
\vskip.3truecm
\item{[\rtag{Ah}]} A. Aharony. (1978). {
Tricritical points in systems with random fields. }
 {\it  Phys. Rev. B} {\bf 18}, 3318--3327.
\item{[\rtag{AW}]}M. Aizenman and  J. Wehr. (1990).
{ Rounding of first or\-der pha\-se tran\-si\-tions 
in sys\-tems with quenched disorder.}
 {\it  Com. Math. Phys.} {\bf 130}, 489--528.
\item{[\rtag{AP}]}J.M.G.  Amaro de Matos and J. F. Perez. (1991).
{ Fluctuations in the Curie-Weiss version of the random field Ising model.}
{\it  J. Stat. Phys.} {\bf 62}, 587--608 .
\item{[\rtag{APZ}]}J.M.G. Amaro de Matos, A E Patrick, and  V A
Zagrebnov.  (1992).
{  Random infinite volume Gibbs states for the Curie-Weiss random field
Ising model.}
{\it  J. Stat. Phys.} {\bf 66}, 139--164.
\item{[\rtag{As}]} S. Asmussen. (1987).{\it Applied Probability and Queues},
J. Wiley and Sons .
\item{[\rtag{B}]}A. Berretti. (1985) {Some properties of random Ising 
models.} {\it  J. Stat. Phys,} {\bf 38}, 483-496 .
\item{[\rtag{Bill}]} P. Billingsley  (1968){\it Convergence of probability measures.}   
Wiley and Sons, NY .
\item{[\rtag{BRZ}]} P. Bleher, J. Ruiz, and V. Zagrebnov.  (1996).
{ One-dimensional random-field Ising model: Gibbs states and structure
of ground states.} {\it  J.Stat. Phys } {\bf 84} 1077--1093.
\item{[\rtag{Bo}]} T. Bodineau.  (1996). {Interface in a one-dimensional Ising spin
system.} {\it Stoch. Proc. Appl.} {\bf 61}, 1--23. 
 \item{[\rtag{BGP4}]} A. Bovier, V. Gayrard and P. Picco.  (1997). 
{ Distribution of profiles for the Kac-Hopfield model.}
 {\it Comm. Math. Phys.} {\bf 186} 323--379.
\item{[\rtag{BZ}]} A. Bovier, and M. Zahradnik.  (1997). 
{ The Low temperature phase of Kac-Ising models.}
{\it J. Stat. Phys.} {\bf 87}, 311--332.
\item{[\rtag{BK}]} J. Bricmont and A. Kupiainen.  (1988).
{Phase transition in the three-dimensional 
random field Ising model.}
 {\it Com. Math. Phys.},{\bf 116}, 539--572.
\item{[\rtag{Br}]} D.C. Brydges. (1986) {\it A short course on cluster expansions.}
in Critical phenomena, random systems, gauge theories, ed
K. Osterwalder, R. Stora, Les Houches XLlll, North Holland.
\item{ }
R. Kotecky and  D. Preiss.  (1996).
{ Cluster expansion for abstract polymer models.}  
{\it Comm. Math. Phys.} {\bf 103} 491--498.
\item{  } B. Simon. (1993)
{\it The statistical mechanics of lattice gases.}
 Vol I, Princeton University Press .
\item{[\rtag{COP1}]} M. Cassandro, E. Orlandi, and P.Picco.  (1999).
{ Typical configurations for one-dimensional random field Kac model. }
{\it  Ann. Prob.}  {\bf 27}, No 3, 1414-1467. 
\item{[\rtag{COP2}]}   M. Cassandro, E.Orlandi, and  P. Picco.  (2002).    
{  Uniqueness and global stability of the interface in a model of phase separation.}
{ \it Nonlinearity}   {\bf 15} No. 5, 1621-1651. 
\item{[\rtag{COPV}]} M. Cassandro, E. Orlandi, P. Picco and
M.E. Vares  (2005).
{ One-dimensional random field Kac's model: Localization of the Phases}
{\it  Electron. J. Probab.}, {\bf 10}, 786-864. 
\item{[\rtag{COP}]} M. Cassandro, E. Orlandi, and E. Presutti.  (1993).
{ Interfaces and typical Gibbs configurations for one-dimensional
Kac potentials.} { \it Prob. Theor. Rel. Fields} {\bf 96}, 57-96. 
 \item{[\rtag{CP}]} M. Cassandro, and  E. Presutti.   (1996).
{ Phase transitions in Ising systems with long but finite range
interactions.}
{\it  Markov Process.  Rel. Fields,}  {\bf 2}, 241-262. 
\item{[\rtag{DMM}]} G. Dal Maso and L. Modica. (1986).
{Non linear stochastic
homogenization.} {\it  Ann. Math. Pura. Appli.} {\bf 4} 144 347--389 
\item{[\rtag{EK}]} S.N. Ethier and T.G. Kurtz. (1986).
{\it Markov Processes
Characterization and convergence.} Wiley and Sons
\item{[\rtag{FFS}]} D.S. Fisher, J. Fr{\"o}hlich, and T. Spencer.  (1984).
{ The Ising model in a random magnetic field.} 
 {\it J. Stat. Phys.}  {\bf 34}, 863--870. 
\item{[\rtag{Go}]}A.O. Golosov.  (1986).
{ On limiting distribution for a random walk in random environment.} 
{\it  Russian Math. Surveys}   {\bf 41}, 	199--200.
\item{[\rtag{I}]}J. Imbrie. (1985). 
{ The ground states of the three-dimensional 
random field Ising model.}
 {\it  Com. Math. Phys.}  {\bf 98}, 145--176 . 
\item{[\rtag{IM}]}Y. Imry and  S.K. Ma.  (1975).
{ Random-field instability of the 
ordered state of continuous symmetry.}
 {\it Phys. Rev. Lett.} {\bf 35}, 1399--1401. 
\item{[\rtag{KUH}]} M. Kac, G. Uhlenbeck, and P.C. Hemmer. (1963){  On the van
der Waals theory of vapour-liquid equilibrium. I. Discussion of a
one-dimensional model.} {\it  J.  Math. Phys.} {\bf 4}, 216--228;
(1963){ II. Discussion of the distribution functions.} 
{\it J. Math. Phys.} {\bf 4}, 229--247;
(1964) { III. Discussion of the critical region.} 
{\it J. Math. Phys.} {\bf 5}, 60--74.
\item{[\rtag{Ke}]} H. Kesten.  (1986) {The limit distribution of Sinai's 
random walk in random environment.}
{\it Physica} 138A, 299-309. 
\item{[\rtag{KS}]} I. Karatzas and S.E. Shreve.  (1991). 
{\it Brownian motion and stochastic calculus.}
Second Edition  Springer Berlin-Heidel\-berg-New York. 
\item{[\rtag{Ku}]} C. K\"ulske.  (2001).
{ On the Gibbsian nature of the random field Kac model under 
block-averaging.}
{\it  J. Stat. Phys.} no 5/6, 991-1012 .
\item{[\rtag{LP}]} J. Lebowitz, and O. Penrose.  (1966).
{Rigorous treatment of the Van der Waals Maxwell theory of the liquid-vapour
transition.} {\it  J. Math. Phys.} {\bf 7}, 98--113.
\item{[\rtag{LMP}]} J. Lebowitz, A. Mazel, and  E. Presutti.   (1999).
{  Liquid-vapor phase transition for systems with finite-range interactions
system.} {\it  J. Stat. Phys.} {\bf 94}, no 5/6, 955-1025. 
 \item{[\rtag{LT}]} M. Ledoux and M. Talagrand.  (1991). {\it Probability in 
 Banach Spaces.} Springer, Berlin-Heidel\-berg-New York.
\item{[\rtag{NP}]} J. Neveu and J. Pitman. (1989).  
{ Renewal property of the extrema and tree property of the excursion
of a one dimensional Brownian motion}
S\'eminaire de Probabilit\'e XXlll, Lect. Notes in Math. {\bf 1372}
239--247,  Springer  Berlin-Heidel\-berg-New York.
\item{[\rtag{PL}]} O. Penrose and J.L. Lebowitz.  (1987).
{ Towards a rigorous molecular theory of metastability.} 
in Fluctuation Phenomena ( e.W. Montroll and J.L. Lebowitz
ed) North-Holland Physics Publishing. 
 \item{[\rtag{RY}]} D. Revuz  and M. Yor. (1991)
{\it Continuous Martingales and Brownian Motion.}
Springer Verlag,  Berlin-Heidel\-berg-New York.
\item{[\rtag{S}]} Y. Sinai. (1982). { The limiting behavior of a one-dimensional
random walk in random environment} {\it Theory of Prob. and its  Appl.} {\bf 27}, 256--268.
\item{[\rtag{SW}]} S.R. Salinas and W.F. Wreszinski.  (1985).
{ On the mean field Ising model in a random external field.} 
{\it  J. Stat. Phys. }{\bf 41}, 299--313.
\item{[\rtag{Sk}]} A.V. Skorohod. (1956) {  Limit Theorems for stochastic
processes} {\it  Theor. Probability Appl. } {\bf 1} 261-290.
\item{[\rtag{T}]} H. M. Taylor.  (1975)
{  A stopped Brownian Motion Formula } 
{\it Ann. Prob.} {\bf 3} 234--236. 
\item{[\rtag{Wi}]}  D. Willians. (1975).
{ On a stopped Brownian Motion Formula of H.M. Taylor} 
Sem. Prob.  Strasb. X, Springer Lectures Note in Math.  {\bf 511} 235--239.

\end